\documentclass[12pt]{article}

\usepackage{frobeniusgraphcalc1}
\usepackage{array}

\usepackage{color}
\usepackage{amscd,amsfonts,amssymb,amscd,amsmath,latexsym}
\usepackage[pdftex]{hyperref}
\usepackage[all]{xy}
\usepackage{mathrsfs}
\usepackage{graphicx}

\newtheorem{theorem}{Theorem}[section] 
\newtheorem{prop}[theorem]{Proposition}
\newtheorem{lem}[theorem]{Lemma}
\newtheorem{ddd}[theorem]{Definition}
\newtheorem{kor}[theorem]{Corollary}

\newtheorem{fact}[theorem]{Fact}
\newtheorem{rem}[theorem]{Remark}
\newtheorem{ex}[theorem]{Example}

\newcommand{\Ree}{\mathrm{Re}}
\newcommand{\cF}{\mathcal{F}}
\newcommand{\Gr}{\mathrm{Gr}}
\newcommand{\Bun}{{\tt Bun}}
\newcommand{\Aut}{{\tt Aut}}
\newcommand{\cB}{{\mathcal{B}}}

\newcommand{\bN}{{\mathbf{N}}}
\newcommand{\DD}{{\tt DD}}

\newcommand{\tmf}{{tmf}}
\newcommand{\bW}{{\mathbf{W}}}

\newcommand{\Thom}{{\tt Thom}}

\newcommand{\hocolim}{{\tt hocolim}}
\newcommand{\ev}{{\tt ev}}
\newcommand{\ord}{{\tt ord}}
\newcommand{\im}{{\tt im}}
\newcommand{\id}{{\tt id}}

\newcommand{\cA}{{\mathcal{A}}}

\newcommand{\cM}{{\mathcal{M}}}
\newcommand{\cG}{{\mathcal{G}}}
\newcommand{\hA}{{\mathbf{\hat A}}}
\renewcommand{\ker}{{\tt ker}}
\def\hB{\hspace*{\fill}$\Box$ \newline\noindent}
\newcommand{\proof}{{\it Proof.$\:\:\:\:$}}
\newcommand{\nat}{{\mathbb{N}}}
\newcommand{\Z}{{\mathbb{Z}}}
\newcommand{\Q}{{\mathbb{Q}}}
\newcommand{\R}{{\mathbb{R}}}
\newcommand{\C}{{\mathbb{C}}}
\newcommand{\tM}{{\tt M}}
\newcommand{\image}{{\tt image}}
\newcommand{\Hom}{{\tt Hom}}
\newcommand{\Td}{{\mathbf{Td}}}
\newcommand{\ch}{{\mathbf{ch}}}
\newcommand{\bV}{{\mathbf{V}}}
\newcommand{\bR}{{\mathbf{R}}}
\newcommand{\bE}{{\mathbf{E}}}
\newcommand{\bU}{{\mathbf{U}}}
\newcommand{\Dirac}{{D\hspace{-0.27cm}/}}
\newcommand{\Tr}{{\tt Tr}}
\newcommand{\sign}{{\tt sign}}
\newcommand{\ind}{{\tt index}}
\newcommand{\pr}{{\tt pr}}
\renewcommand{\lim}{{\tt lim}}
\newcommand{\colim}{{\tt colim}}
\newcommand{\Sym}{{\tt Sym}}

\title{On the topological contents of $\eta$-invariants}
\author{Ulrich Bunke\thanks{NWF I - Mathematik,
Universit{\"a}t Regensburg,
93040 Regensburg,
GERMANY, ulrich.bunke@mathematik.uni-regensburg.de}}

\begin{document}
\maketitle

\begin{abstract}
We discuss an universal bordism invariant obtained from the Atiyah-Patodi-Singer $\eta$-invariant from the analytic and homotopy theoretic point of view.
Classical invariants like the Adams $e$-invariant, $\rho$-invariants and $String$-bordism invariants are derived as special cases. The main results are a secondary index theorem about the coincidence of the analytic and topological constructions and intrinsic expressions for the bordism invariants.
\end{abstract}

\tableofcontents

\section{Introduction}\label{subsecintro}

The purpose of this work is to understand which topological information is encoded in  the $\eta$-invariant,
a spectral-geometric invariant introduced by 
Atiyah-Patodi-Singer  \cite{MR0397797} in the context of index theory for boundary value problems for Dirac-operators. We are in particular interested in bordism invariants derived from the 
  $\eta$ invariant. By now we know  many examples, see e.g. \cite{MR0397799},\cite{MR870805}, \cite{MR883375}, \cite{2009arXiv0912.4875B},  \cite{MR762355}, \cite{westburyjones}, \cite{2010arXiv1012.5237C}.  
   In the present paper we consider the universal structure behind these examples. We define a bordism invariant which we call the universal $\eta$-invariant. We use Section \ref{sect5} in order to review some of the known $\eta$-invariant based bordism invariants. We put the emphasis on the demonstration  how they can be interpreted as special cases of our universal construction.  
 Our construction also subsumes (by constructions similar to the one in  Subsection \ref{subseccg})  some, but not all, of the    Kreck-Stolz type or Eells-Kuiper type invariants, \cite{MR932303}, \cite{MR0372929}, or the generalised Rochlin invariants of \cite{MR889873}.  The universal $\eta$-invariant does not seem  to incorporate the invariant introduced by Kervaire-Milnor
in order to distinguish exotic spheres \cite{MR0148075}.

In the present paper we introduce and compare two versions of the  universal $\eta$-invariant.  The analytic version  $\eta^{an}$  given in Definition \ref{etaandef} is the bordism invariant which is derived from the appearance of the reduced  $\eta$-invariant in the local index theorem for the Atiyah-Patodi-Singer boundary value problem by cancelling out the dependence on geometric data. The ideas for this construction are 
more or less standard and have been used previously in many special situations.
  
The topological counterpart $\eta^{top}$ introduced in Definition \ref{etatopdefv} is constructed by a simple  homotopy-theoretic consideration using the interplay of $\Q/\Z$-bordism and $K$-theory.

While it is not so complicated to see  that $\eta^{an}$ is a bordism invariant, to understand its homotopy-theoretic meaning is slightly deeper. The bridge between analysis and topology is provided by our first main Theorem \ref{indthm}
stating that $$\eta^{an}=\eta^{top}\ 
.$$
Its proof uses standard methods in index theory like the analytic picture of  $K$-homology \cite{MR679698}, \cite{MR918241}, 
the Atiyah-Patodi-Singer index theorem \cite{MR0397797}, and some ideas from $\Z/l\Z$-index theory \cite{MR1144425}.

Bordism classes can be  represented geometrically by manifolds with additional structures, called cycles (see Subsection \ref{sec1400} for details). It is then an interesting question how one can calculate the universal $\eta$-invariant or its specialisations in terms of the cycle.  The definition of both, the topological or the analytical version of the universal $\eta$-invariant involves the choice of a zero bordism of some multiple of the cycle. In applications it is often complicated to find such a zero bordism.
It is a striking advantage of the analytic picture that it can be reorganised to an expression which only
involves structures on the cycle itself. In special cases this has been previously exploited in \cite{MR687857}, \cite{MR762355} (the case of the Adams $e$-invariant, see Subsection \ref{subsec51}),  and 
\cite{2009arXiv0912.4875B} (to calculate $String$-bordism invariants, see  Subsection \ref{subsec55}).

We consider the intrinsic formula  for $\eta^{an}$ given in Theorem \ref{them2} as one of the main original contributions of the present paper. This formula is based on a new object which we call a geometrisation.
If a map $f:M\to BG$  classifies a principal $G$-bundle $P$, then a geometrisation of $f$ is essentially given by a connection on $P$. In general,
the notion of a geometrisation  partially generalises
the notion of a connection on the in general non-existent principal bundle classified by the map $f:M\to B$ for an arbitrary space $B$. The details are slightly more complicated since
we will take structures on the normal bundle into account. 

\begin{rem}{\rm
In this paper we generally decided to work with complex $K$-theory.
We think that there is a real version of the whole theory which can be obtained by replacing
complex $K$-theory by real $KO$-theory, $BSpin^{c}$ by $BSpin$,
and taking the real structures on the spinor bundles into account properly on the analytic side.
The real version of the universal $\eta$-invariant would be slightly stronger than its complex counterpart
which loses some two-torsion classes. In order to recover the Adams $e$-invariant or the string bordism invariant \cite{2009arXiv0912.4875B} completely as special cases of the universal $\eta$-invariant  we would need the real version.}
\end{rem}

Let us now describe the contents of the paper.
In Section \ref{sec2222} we introduce the topological version $\eta^{top}$ of the universal $\eta$-invariant and study its properties. Most interesting is probably the relation with the Adams spectral sequence Proposition \ref{prop1} which asserts   that 
the universal $\eta$-invariant  detects the first non-trivial sub-quotient of the bordism theory
with respect to the $K$-theory based Adams filtration.

In Section \ref{sec333} we introduce the analytic version $\eta^{an}$ of the universal $\eta$-invariant and prove the secondary index Theorem \ref{indthm} stating that $\eta^{an}=\eta^{top}$. Before we can define $\eta^{an}$ we have to recall in Subsections \ref{sec1400} and \ref{sec1677} some preliminary technical details   concerning the relation of structures on the stable normal bundle as they come out of the Pontrjagin-Thom construction, and structures on the tangent bundle which will be used to do geometry and analysis.

Section \ref{sec444} is devoted to geometrizations (Definition \ref{def903})  and the intrinsic formula for $\eta^{an}$ (Theorem \ref{them2}).

In the last Section \ref{sect5} we discuss in detail various specialisations of the universal $\eta$-invariant. It contains mainly a review of known constructions and results with slight improvements or generalisations at some points (e.g. Corollary \ref{thm7v}). In the Propositions 
 \ref{lem77} and \ref{prop77} we show how the usual geometric structures of $Spin$- and $String$-geometry (see \cite{2009arXiv0906.0117W} for the latter) give rise to geometrisations which lead to the known intrinsic formulas for the corresponding bordism invariants. It has  been the initial motivation for this work
to understand the general principles behind the $String$-bordism invariants introduced in \cite{2009arXiv0912.4875B}. It should be easy to adapt the  arguments used here for the $String=MO\langle8\rangle$-bordism case
  to bordism theories $MO\langle n\rangle$ associated to higher connected covers $BO\langle n\rangle$ of $BO$.

\bigskip

  


{\it Acknowledgement: I thank Bernd Ammann, Sebastian Goette,  Diarmuid Crowley,   Niko Naumann for stimulating discussions. I am in particular grateful to M. V\"olkl for suggesting various inprovements. 
The pictures have been typeset using the {\tt frobeniusgraphcalc.sty}-package written by Clara L\"oh.}

\section{The topological construction}\label{sec2222}

\subsection{Motivation}

In this section we use complex $K$-theory in order to detect torsion elements in the homotopy groups $\pi_{*}(E)$ of a spectrum $E$. We introduce the topological version $\eta^{top}$ of the universal $\eta$-invariant  as a secondary version of the degree
$\epsilon:\pi_{*}(E)\to \pi_{*}(K\wedge E)$. The idea is to first lift
the torsion element to the $\Q/\Z$-version of $E$, then apply the degree, and finally to detect
the result via its evaluations against $K^{0}(E)$.  This construction is a generalisation of the 
construction of the Adams $e$-invariant in the case of the sphere spectrum $E=S$.

In Subsection \ref{sec141} we will give the construction of the invariant. In Subsection \ref{sec222} we
analyse its target in some detail. Finally, in Subsection \ref{secneu8} we understand completely in terms of the Adams spectral sequence which piece of
$\pi_{*}(E)_{tors}$ the universal $\eta$-invariant can detect.

\subsection{The definition of $\eta^{top}$}\label{sec141}

For an abelian group $A$ we let
 $\tM A$ denote the Moore spectrum 
which  is characterised by  its integral homology
$$\pi_{n}(H\Z\wedge \tM A)\cong \left\{\begin{array}{cc} 
A&n=0\\0&n\not=0\end{array}\right. \ .$$
More generally, for a spectrum $E$ we have short
exact sequences
\begin{equation}\label{torseq}0\to \pi_{n}(E)\otimes A\to \pi_{n}(E\wedge \tM A)\to \pi_{n-1}(E)*A\to 0\end{equation}
(\cite[(2.1)]{MR551009}) for all $n\in \Z$. In order to simplify the notation
  we abbreviate $EA:=E\wedge \tM A$. 
  The equivalences $\tM \Z\cong S$ and $\tM \Q\cong H\Q$ induce 
equivalences
$E\Z\cong E$ and $E\Q\cong E\wedge H\Q$. 

The starting point for the construction of  the topological version $\eta^{top}$ of the universal $\eta$-invariant  is the fibre sequence
\begin{equation}\label{eq631}\tM\Z\to \tM\Q\to \tM\Q/\Z\to \Sigma \tM\Z\end{equation}
of Moore spectra. Smashing this sequence with the map $\epsilon:E\to K\wedge E$ induced by the unit of the $K$-theory spectrum we get the diagram 
\begin{equation}\label{eq3}\xymatrix{\Sigma^{-1}E\Q \ar@{.>}[dr]^{u}\ar[d]\ar[r]^{\epsilon}& \Sigma^{-1} K\wedge E\Q \ar[d]\\ {}^{\hat x}\ \Sigma^{-1}E\Q/\Z \ar[r]^{\epsilon}\ar[d]&{}^{\tilde x}\   \Sigma^{-1}K\wedge E\Q/\Z \ar[d]\\{}^{x}\ E   \ar[d]\ar[r]^{\epsilon}&\ K\wedge E \ar[d]\\{}^{0}\ E\Q \ar[r]&K\wedge E \Q }\ ,\end{equation}
where we use the symbol $\epsilon$ also to denote other maps naturally induced by the unit.

\begin{rem}\label{proffinite1}{\rm
If $E$ is a spectrum and $X$ is a space or spectrum, then for $n\in \Z$ we consider the cohomology $E^{n}(X)$ as a topological group  where a basis of neighbourhoods of zero is given by the kernels of restrictions along maps $Y\to X$ from finite $CW$-complexes or finite cell spectra $Y$.
The topology on $E^{n}(X)$ is called the profinite topology. This should not be confused   with the notion of a profinite group in algebra. 
}
\end{rem}
 
We have an evaluation pairing 
 \begin{equation}\label{neu15neu}\langle-,-\rangle:\pi_{n}(\Sigma^{-1}K\wedge E\Q/\Z)\otimes K^{0}(E)\to \pi_{n+1}(K\Q/\Z)\end{equation}
which sends $x\otimes \phi$ to the composition
$$\Sigma^{n+1}S\stackrel{x}{\to} K\wedge E\wedge \tM\Q/\Z\xrightarrow{\id_{K}\wedge\phi\wedge \id_{\tM\Q/\Z}}K\wedge K\wedge \tM\Q/\Z \stackrel{\mu}{\to} K\wedge \tM \Q/\Z\ ,$$
where $\mu$  is the multiplication of the ring spectrum $K$.
Its adjoint
is a map
$$\pi_{n}(\Sigma^{-1}K\wedge E\Q/\Z)\to \Hom^{cont}(K^{0}(E), \pi_{n+1}(K\Q/\Z))\ ,$$
where $\Hom^{cont}(-,-)$ stands for continuous homomorphisms and $ \pi_{n+1}(K\Q/\Z)$ has the discrete topology.
We let  
\begin{equation}\label{sec31u}U\subseteq  \Hom^{cont}(K^{0}(E),\pi_{n+1}(K\Q/\Z))\end{equation} denote the subgroup
given by the pairings with the elements in the image
$$\image(u:\pi_{n}(\Sigma^{-1}E\Q)\to \pi_{n}(\Sigma^{-1}K\wedge E\Q/\Z))$$
of the map $u$ which can be read off from \eqref{eq3}.
\begin{ddd}\label{defqq} For every $n\in \Z$
we  define the abelian group
\begin{equation}\label{eq5}
Q_{n}(E):=\frac{ \Hom^{cont}(K^{0}(E),\pi_{n+1}(K\Q/\Z))}{U}\ .\end{equation}
\end{ddd}

The following definition uses a diagram chase of elements as indicated in \eqref{eq3}.
\begin{ddd}\label{etatopdefv}
The homotopy theoretic version of the universal $\eta$-invariant is  the homomorphism 
$$\eta^{top}:\pi_{n}(E)_{tors}\to Q_{n}(E)$$
defined by the following prescription:
If $x\in \pi_{n}(E)_{tors}$, then we choose a lift
$\hat x\in \pi_{n}(\Sigma^{-1}E\Q/\Z)$. 
We let $\eta^{top}(x)\in Q_{n}(E)$ be the class represented by the homomorphism 
given by the pairing against $\epsilon(\hat x)$.
\end{ddd}

We must show that $\eta^{top}$ is well-defined. Indeed, the choice of $\hat x$ is unique up to elements which come from $\pi_{n}(\Sigma^{-1}E\Q)$, but this ambiguity is taken care of by taking the quotient by $U$ in the definition  \eqref{eq5} of $Q_{n}(E)$.

 The following Lemma immediately follow from the definitions:
\begin{lem}\label{lem1000} A map of spectra $E\to F$   naturally induces a commutative diagram
  $$\xymatrix{\pi_{n}(E)_{tors}\ar[d]\ar[r]^{\eta^{top}}&Q_{n}(E)\ar[d]\\
\pi_{n}(F)_{tors}\ar[r]^{\eta^{top}}&Q_{n}(F)}\ .$$
\end{lem}

\subsection{Simplification of $Q_{n}(E)$}\label{sec222}

Because of its definition as a quotient it is difficult to define maps out of  $Q_{n}(E)$. In this subsection we analyze the structure of this group and explain how one can detect its elements.
We consider the ring spectrum
$$HP\Q:=H\Q[b,b^{-1}]\ ,$$
where $\deg(b)=-2$. Additively it has a decomposition 
\begin{equation}\label{eq8560}HP\Q\cong \bigvee_{i\in \Z}  \Sigma^{2i}H\Q\end{equation} and   for each $i\in \Z$ we consider the projection $p_{2i}:HP\Q\to \Sigma^{2i}H\Q$ to the corresponding component. We set $p_{n}:=0$ for odd $n$. 
The Chern character is an equivalence of ring spectra  $$\ch:K\Q\stackrel{\sim}{\to} HP\Q\ .$$

The composition of the Chern character  with  the projection  $p_{n+1}$ gives a map $$p_{n+1}\circ \ch :K^{0}(E)\to H\Q^{n+1}(E)$$
whose kernel will be denoted by 
\begin{equation}\label{eq131}V_{n}:=\ker\left(p_{n+1}\circ \ch \right)\ .\end{equation}
\begin{lem}\label{lem1} 
\begin{enumerate}
\item
The restriction to $V_{n}\subseteq K^{0}(E)$ induces a well-defined map
\begin{equation}\label{eq100}Q_{n}(E)\to \Hom^{cont}(V_{n},\pi_{n+1}(K\Q/\Z))\ .\end{equation}
\item This restriction map 
 is an isomorphism if we assume  that   $E$ is lower bounded and $\dim H\Q^{n+1}(E)$ is finite.
 \end{enumerate}
\end{lem}
\proof
First we show that the restriction is well-defined. We must show that if $\phi\in V_{n}$, then the pairing of $\phi$ with  $u(y)\in \pi_{n+1}(K\wedge E\Q/\Z)$  vanishes  for every
$y\in \pi_{n+1}(E\Q)$.
This follows from   the equality
\begin{equation}\label{t2013}\langle u(y),\phi\rangle=q(\langle \ch(\epsilon(y)), p_{n+1}(\ch(\phi))\rangle)\end{equation}
where $$q: \pi_{n+1}(HP\Q)\stackrel{\ch^{-1}}{\to} \pi_{n+1}(K\Q)\to \pi_{n+1}(K\Q/\Z)$$
and $\epsilon$ is as in \eqref{eq3}.


We now show 2.  We use the general fact that if $f:A\to V_{n}$ is a homomorphism of an abelian group into  a $\Q$-vector space such that its image is finitely generated as an abelian group, then there exists a splitting $A\cong \ker(f)\oplus A^{\prime}$. Indeed, in this case the image is free and hence projective.

Note that there exists an integer $N$ (only depending on $n$ and the lower bound of $E$) such that the image of 
$p_{n+1}\circ  \ch(\dots)$ is contained in the image of $\frac{1}{N} H\Z^{n+1}(E)\to  H\Q^{n+1}(E)$ and is therefore finitely generated as an abelian group since we assume that $H\Q^{n+1}(E)$ is finite dimensional. We conclude that $$K^{0}(E)\cong V_{n}\oplus V_{n}^{c}\ ,$$
where $V_{n}^{c}\cong \image(p_{n+1}\circ  \ch)$ is a free abelian group.
This immediately implies that (\ref{eq100}) is surjective.

Any homomorphism
$\phi\in \Hom^{cont}(K^{0}(E),\pi_{n+1}(K\Q/\Z))$
 can uniquely be decomposed as a sum of its restrictions to $V_{n}$  and $V_{n}^{c}$.
  We claim that $U\cong \Hom(V_{n}^{c},\pi_{n+1}(K\Q/\Z))$, where $U$ is as in (\ref{sec31u}). 
 The claim implies that (\ref{eq100}) is injective.  
 
 We can assume that $n$ is odd.
 The claim follows  from the surjectivity of the following composition
 $$\pi_{n}(\Sigma^{-1}E\Q)\cong \Hom^{cont}(H\Q^{n+1}(E),\Q)\cong \Hom^{cont}(V_{n}^{c},\pi_{n+1}(K\Q))\twoheadrightarrow \Hom^{cont}(V_{n}^{c},\pi_{n+1}(K\Q/\Z))\ , $$
 where we use the Chern character  for the second equivalence and the fact that $V_{n}^{c}$ is free
 in order to conclude the surjectivity of the last map. Note that in the present situation continuity of homomorphisms is automatic by our finiteness assumption.
\hB 

The definition of $\eta^{top}$ is based on first lifting the torsion element in the homotopy group  of $E$ to a $\Q/\Z$- homotopy class   which is then paired with elements of $K$-theory. The pairing with torsion $K$-theory elements can be expressed in a dual way as a pairing of the original homotopy class with $\Q/\Z$-lifts of the $K$-theory elements. We now explain the details.

Assume that $\phi\in K^{0}(E)$ satisfies $\ch(\phi)=0$.
Then $\phi\in V_{n}$ and we get an evaluation $$\ev_{\phi}:Q_{n}(E)\to \pi_{n+1}(K\Q/\Z)\ .$$
In view of the exact sequence
$$K\Q/\Z^{-1}(E)\stackrel{\partial}{\to} K^{0}(E)\stackrel{\ch}{\to} HP\Q^{0}(E)$$
we can choose 
$\hat \phi \in K\Q/\Z^{-1}(E)$ such that $\partial \hat \phi=\phi$.
If we want to calculate $\ev_{\phi}(\eta^{top}(x))$, then instead of lifting the class $x$ to a $\Q/\Z$ class we can instead evaluate the class $\epsilon(x)\in \pi_{n}(K\wedge E)$ against the lift $\hat \phi$.
The following assertion follows easily from the definition of $\eta^{top}$ and  commutativity of  the diagram
$$\xymatrix{\Sigma^{-1}E \wedge \tM\Q/\Z\ar[rr]^{\phi\wedge \id_{\tM\Q/\Z}}\ar[d]^{\partial}&&\Sigma^{-1}K\Q/\Z\ar@{=}[d]\\
E\ar[rr]^{\hat \phi}&&\Sigma^{-1}K\Q/\Z}\ .$$
\begin{lem}\label{lem900}
For $x\in \pi_{n}(E)_{tors}$ and $\hat \phi\in K\Q/\Z^{-1}(E)$   we have
$$\ev_{\phi}(\eta^{top}(x))=\langle  \hat \phi,\epsilon(x) \rangle\ ,$$
where  $\phi:=\partial \hat \phi$.
\end{lem}

\bigskip

In the present paper the spectrum $E$ will often be a Thom spectrum $MB$  associated to 
  a map of spaces $B\to BSpin^{c}$. The spectrum $MB$ is   $K$-oriented by
\begin{equation}\label{betadef}\beta:MB\to MSpin^{c}\xrightarrow{ABS} K\ ,\end{equation}
where $ABS$ is the Atiyah-Bott-Shapiro orientation.  
In this case  can use the  Thom isomorphisms
$$\Thom^{K}:K^{0}(B)\stackrel{\sim}{\to} K^{0}(MB)\ , \quad 
\Thom_{H\Q}:H\Q_{n+1}(MB)\stackrel{\sim}{\to}  H\Q_{n+1}(B) 
$$
in order to express $Q_{n}(MB)$ in terms of $B$.
We have an isomorphism
\begin{equation}\label{spacerep1}Q_{n}(MB)\cong \frac{\Hom^{cont}(K^{0} (B ),\pi_{n+1}(K\Q/\Z))}{U^{\prime}})\ ,\end{equation}
 where $U^{\prime}$ is obtained by precomposing the elements in $U$ (see \eqref{sec31u}) with $\Thom^{K}$.
We let $\Td\in HP\Q^{0}( BSpin^{c})$ be the universal Todd class and use the same symbol in order to denote the pull-back of this class to $B$.
Then using \eqref{t2013}, the  isomorphism 
$$\pi_{n+1}(MB\Q)\cong H\Q_{n+1}(MB)\stackrel{\Thom_{H\Q}}{\cong} H\Q_{n+1}(B)\ ,$$  and the Riemann-Roch theorem we can describe $U^{\prime}$ as the space homomorphisms given by
\begin{equation}\label{localthomformula}K^{0}(B)\ni \phi\mapsto q(\langle  y , p_{n+1}(\Td^{-1}\cup \ch(\phi))\rangle)\in \pi_{n+1}(K\Q/\Z)\end{equation}
for all $y\in H\Q_{n+1}(B)$.

\subsection{Relation with the Adams spectral sequence and $K$-localisation}\label{secneu8}

In this subsection we will see that the topological $\eta$-invariant essentially detects the first line in the 
$K$-based Adams spectral sequence $(E^{*,*}_{r}(E),d_{r})$  for $E$. We refer 
 to \cite{MR1324104} and \cite{MR860042} for details on the Adams spectral sequence.   We define  the spectrum $\bar K$ by the cofibre sequence
$$\Sigma^{-1}\bar K\to S\to K\to \bar K$$
 and form the Adams tower
$$\xymatrix{\Sigma^{-2}E\wedge \bar K\wedge \bar K\ar[dd]&&\\&&\\\Sigma^{-1}E\wedge \bar K\ar[dd]\ar[rr]&&\Sigma^{-1}E\wedge \bar K\wedge K\ar@{..>}[uull]\\&&\\E\ar[rr]&&E\wedge K\ar@{..>}[uull]}\ .$$
The horizontal maps are induced by the unit $S\to K$ while the dotted arrows are degree one-maps which turn the triangles into fibre sequences.
This tower gives rise to an Adams spectral sequence
$$(E_{r}^{*,*}(E),d_{r})_{r\ge 1}$$
 and a filtration $(\cF^{s}\pi_{*}(E))_{s\ge 0}$ of the homotopy groups of $E$, where
  $\cF^{s}\pi_{*}(E)\subseteq \pi_{*}(E)$ is defined as  the image of $\pi_{*}(\Sigma^{-s}E\wedge \bar K^{\wedge s})$.
 We have a natural injective map
 $$s:\Gr^{1}\pi_{*}(E)\to E_{2}^{1,*+1}(E)\ .$$
 
 \begin{prop}\label{prop1}
 \begin{enumerate}
\item We have $\cF^{1}\pi_{*}(E)\subseteq \pi_{*}(E)_{tors}$.
 \item 
 We have $\cF^{2}\pi_{*}(E)\subseteq \ker(\eta^{top})$. Hence $\eta^{top}$ induces a map (denoted by the same symbol)
 $$\eta^{top}:\Gr^{1}\pi_{*}(E)\to Q_{*}(E)\ .$$ 
 \item There exists a map
$$\kappa: E_{2}^{1,*+1}(E)\to Q_{*}(E)$$
such that
$$\eta^{top}=\kappa\circ s:\Gr^{1}\pi_{n}(E)\to Q_{n}(E)\ .$$
\item 
If $\pi_{*}(E\wedge K )$ is torsion-free, then $$\pi_{*}(E)_{tors}= \cF^{1}\pi_{*}(E)\ .$$ If $n$ is odd, then the restriction $$\kappa_{|E_{2}^{1,n+1}(E)_{tors}}:E_{2}^{1,n+1}(E)_{tors}\to Q_{n}(E)$$ of $\kappa$ to the torsion subgroup is injective
Consequently, the map
 $$\eta^{top}:\Gr^{1}\pi_{n}(E)\to Q_{n}(E)$$ induced by $\eta^{top}$ is injective.
\end{enumerate}
 \end{prop}
\proof
We have a map
$$f:K\to K\Q\cong \prod_{p\in \Z} \Sigma^{2p}H\Q\stackrel{\pr}{\to} H\Q$$
which we use to build the diagram
$$\xymatrix{\Sigma^{-2} E\wedge \bar K\wedge \bar K\ar[dd]&& \\&&&\\\Sigma^{-1} E\wedge \bar K\ar[rr]\ar[dd] \ar[dr]&&\Sigma^{-1}E\Q/\Z\ar[dd]\ar[dr]&\\ &\Sigma^{-1}E\wedge \bar K\wedge K\ar[rr]\ar@{..>}[uuul]&&\Sigma^{-1}E\Q/\Z\wedge K\\E\ar@{=}[rr]\ar[rd]&&E\ar[dr]&\\
&E\wedge K\ar@{..>}[uuul]\ar[rr]^{f}&&E\Q\ar@{..>}[uuul]}\ .$$
It connects a piece of the Adams tower of $E$ with the basic diagram used in the definition of $\eta^{top}$. 

Assertions 1. and 2. follow immediately from a diagram chases. We now show assertion 3.
The map
$$E_{1}^{1,*+1}(E)=\pi_{*}(\Sigma^{-1}E\wedge \bar K\wedge K)\to \pi_{*}(\Sigma^{-1}E\Q/\Z\wedge K)\to Q_{*}(E)$$
annihilates the image of $E\wedge K$, i.e the image of the boundary map of the spectral sequence and therefore factors through a map
$$\kappa:E_{2}^{1,*+1}(E)\to Q_{n}(E)\ .$$
Assertion 3.  again follows by again by a diagram chase.

We now show Assertion 4.  In view of 1. the first part is clear. As in \cite[Sec. 5.3]{2009arXiv0912.4875B}   our additional assumption implies that we have a short exact sequence
\begin{equation}\label{eq112}0\to   E_{2}^{0,n+1}(E )\otimes \Q/\Z\stackrel{j}{\to}    E_{2}^{0,n+1}(E\Q/\Z ) \stackrel{\delta}{\to}  E_{2}^{1,n+1}(E)_{tors}\to 0\ .\end{equation}
Let $\gamma\in E^{1,n+1}_{2}(E)_{tors}$ and assume  that
$\kappa(\gamma)=0$.
We write $\gamma=\delta (\beta)$ for some $\beta\in E_{2}^{0,n+1}(E\Q/\Z)$. 
We have $ \kappa(\gamma)=0$ if an only if there exists $y\in \pi_{n+1}(E\Q )$ which induces the same pairing with $K^{0}(E)$ as $a(\beta) $, where
$a:E_{2}^{0,n+1}(E\Q/\Z)\to \pi_{n+1}(E\Q/\Z\wedge K)$ is the inclusion.
By Pontrjagin duality   we have
$$\pi_{n+1}(E\Q/\Z \wedge K )\cong \Hom^{cont}(K^{n+1}(E),\Q/\Z) \cong \Hom^{cont}(K^{0}(E),\pi_{n+1}(K\Q/\Z))\ .$$ We conclude that $a(\beta)=a(j(q(\epsilon(y))))$, where 
$\epsilon(y) \in E_{2}^{0,n+1}(E\Q)$  and $$q:E_{2}^{0,n+1}(E\Q)  \cong E_{2}^{0,n+1}(E)\otimes \Q \to E_{2}^{0,n+1}(E)\otimes \Q/\Z \ .$$ Since $a$ is injective this implies
$\gamma=\delta(j(q(\epsilon(y))))=0$. \hB 

We let $E\to E_{K}$ denote the Bousfield localisation of the spectrum $E$ at the complex $K$-theory spectrum. The following Lemma shows that the universal $\eta$-variant factorizes over the $K$-localisation.

\begin{lem}
We have a commuting diagram
$$\xymatrix{\pi_{n}(E)_{tors}\ar[r]^{\eta^{top}}\ar[d]&Q_{n}(E)\ar[d]^{\cong}\\
\pi_{n}(E_{K})_{tors}\ar[r]^{\eta^{top}}&Q_{n}(E_{K})}\ .$$
 \end{lem}
\proof
Since $E\to E_{K}$ induces an isomorphism in $K$-theory we conclude that the right vertical map is an equivalence. The rest is Lemma \ref{lem1000}. \hB

\section{The spectral geometric construction}\label{sec333}

\subsection{Motivation}

In this section we define an analytic invariant  $\eta^{an}$ of torsion elements in the $B$-bordism theory. The analytic invariant will be derived from geometric and spectral geometric quantities
associated to geometric cycles  for bordism classes. The relation between the 
geometric  and homotopy theoretic picture of the bordism group is given by  Thom-Pontrjagin construction, see \cite[Ch IV.7]{MR1627486}.  
In Subsection \ref{sec1400} we give the details of the geometric picture of the $B$-bordism theory.
Subsection \ref{sec1677} is devoted to some technical details on the transfer of $Spin^{c}$-structures from the normal bundle to the tangent bundle.
A reader with some experience with the Thom-Pontrjagin construction and $Spin^{c}$-structures may immediately proceed to the construction of $\eta^{an}$ in Subsection \ref{sec81}. The final Subsection \ref{theproof} of this part contains the proof of the theorem about the equality of the analytic and topological universal $\eta$-invariant.



%

%

\subsection{Geometric cycles for $B$-bordism theory}\label{sec1400}

We assume that we have chosen representatives of the homotopy types of the classifying spaces of the classical Lie groups like $O(k)$, $Spin^{c}(k)$ etc. and representatives of the homotopy classes  of the usual maps connecting them like $BSpin^{c}(k)\to BO(k)$ or $\iota_{k}:BO(k)\to BO(k+1)$.
We further choose universal euclidean bundles $\xi_{k}$ on $BO(k)$ and isomorphisms
$\iota_{k}^{*}\xi_{k+1}\cong \xi_{k}\oplus BO(k)\times \R$ for all $k\in \nat$.

We consider a map of spaces  $B\to BSpin^{c}$ and the corresponding Thom spectrum $MB$.
Cycles for elements of the $B$-bordism group $\pi_{n}(MB )$ of $X$ are pairs $(M,f)$ consisting of      a closed $n$-dimensional Riemannian manifold   $M$ and a stable normal $B$-structure  on the map  $f:M\to B$.  The additional data of a stable normal $B$ structure is not written explicitly and fixes the relation between the map $f$ and the tangent bundle of $M$.
In the following we explain the details. Since $M$ is compact   there exists a factorisation
$$\xymatrix{&&BO(k)\ar[d]\\M\ar@{.>}[urr]^{\hat f}\ar[r]^{f}&B\ar[r]&BO}$$
up to homotopy for a suitable integer $k$.
 We require that $\hat f^{*} \xi_{k}$ is a complement of the tangent bundle of $M$, i.e. that there exists an isomorphism
 \begin{equation}\label{eq1360}TM\oplus \hat f^{*}\xi_{k}\cong M\times \R^{n+k}\ .\end{equation}
\begin{ddd} A normal $B$-structure on $f$ consists of the choice of $\hat f$ and the isomorphism \eqref{eq1360}.
 \end{ddd}
 There is an obvious notion of a stabilisation of a normal $B$-structure which allows to increase $k$. A stable normal $B$-structure is an equivalence class of normal $B$-structures under the relation generated by stabilisation.

%

The bordism group $\pi_{n}(MB)$ is the set of equivalence classes of cycles, where the equivalence relation is given by bordism, and the group structure is induced by the disjoint sum.
 A zero bordism of $(M,f)$ is given by a pair $(W,F)$ of similar data, where $W$ is a compact $n+1$-dimensional Riemannian manifold with boundary $\partial W\cong M$  and product structure and   $F:W\to B$ carries a stable normal $B$-structure  which extends the one on $f$. 
   In detail this means the following. First of all the map $\hat F$ extends $\hat f$. The   stable normal $B$-structure on $W$ is  represented by an  isomorphism
\begin{equation}\label{eq1361}TW\oplus \hat F^{*}\xi_{k}\cong W\times \R^{n+1+k}\ .\end{equation}  The outgoing normal field of $TW_{|\partial W}$ provides
an orthogonal decomposition
\begin{equation}\label{eq1361q}TW_{|M}\cong TM\oplus (M\times \R)\ . \end{equation} 
It is here where we use the additional datum of the Riemannian metrics in order to rigidify the choice of the unit normal vector field.
We require that the isomorphism 
$$TM\oplus (M\times \R)\oplus \hat f^{*}\xi_{k} \stackrel{\eqref{eq1361q}}{\cong}  TW_{|M}\oplus \hat f^{*}\xi_{k}\stackrel{\eqref{eq1361}_{|M}}{\cong} M\times \R^{n+1+k}$$ represents the stable normal $B$-structure on $f$.

%
\subsection{Normal and tangential $Spin^{c}$-structures}\label{sec1677}

Because of the factorisation $B\to BSpin^{c}\to BO$ 
a normal $B$-structure induces a normal $Spin^{c}$-structure. As we will do geometry on the tangent bundle we must transfer normal $Spin^{c}$-structures to tangential $Spin^{c}$-structures. The homotopy theoretic picture of this transition is explained  in \cite[Sec. 8]{2009arXiv0912.4875B} in the example of $String$-structures. In the following we describe its geometric counterpart.

Let $V\to M$ be an $m$-dimensional real vector bundle. Then a  $Spin^{c}$-structure on $V$ is pair $(P,\kappa)$, where
$P\to M$ is a $Spin^{c}(m)$-principal bundle and $\kappa$ is an isomorphism of real vector bundles
$$\kappa:P\times_{Spin^{c}(n)}\R^{m}\cong V\ .$$
With this definition a $Spin^{c}$-structure induces an euclidean metric and an orientation on $V$ so that the oriented orthonormal frame bundle is $SO(V):=P\times_{Spin^{c}(m)} SO(m)$.

The collection  of all $Spin^{c}$-structures on the vector bundle  $V$  naturally forms a groupoid 
 $Spin^{c}(V)$. 
The objects of the groupoid  $Spin^{c}(V)$ are the $Spin^{c}$-structures $(P,\kappa)$, and the morphisms
$(P,\kappa)\to (P^{\prime},\kappa^{\prime})$
are isomorphisms of $Spin^{c}(m)$-principal bundles $P\to P^{\prime}$ which are compatible with the isomorphisms $\kappa$ and $\kappa^{\prime}$. This in particular implies that automorphisms of $(P,\kappa)$  are given by the central action of $C^{\infty}(M,U(1))$ on $P$.

If we associate to any open subset
 $A\subseteq M$ the groupoid $Spin^{c}(V_{|A})$,  then we obtain a sheaf of groupoids $\underline{Spin^{c}(V)}$ which actually 
is an $U(1)$-banded gerbe. We refer to \cite{MR2362847},\cite{MR2681698}, or \cite{MR1876068} for an introduction to gerbes. 
Isomorphism classes of $U(1)$-banded gerbes $\cG$ are classified by their Dixmier-Douady classes
$\DD(\cG)\in H^{3}(M;\Z)$. In particular,  the Dixmier-Douday class of the $Spin^{c}$-gerbe  $\underline{Spin^{c}(V)}$
  is the class
$$\DD(\underline{Spin^{c}(V)})=W_{3}(V)=\beta (w_{2}(V))\in H^{3}(M;\Z)\ ,$$ where
$w_{2}(V)\in H^{2}(M;\Z/2\Z)$ is the second Stiefel-Whitney class and $\beta:H^{2}(M;\Z/2\Z)\to H^{3}(M;\Z)$ is the Bockstein operator \cite[Thm. D2]{MR1031992}. The groupoid $Spin^{c}(V)$ is non-empty  if and only if $W_{3}(V)=0$, i.e. the class $W_{3}(V)$ is the obstruction against the existence of a $Spin^{c}$-structure on $V$. 
 In the following we will simplify the notation and write $P$ for the $Spin^{c}$-structure $(P,\kappa)$.  
 
Let $\cB U(1)(M)$ denote the Picard groupoid (see \cite{del73}) of $U(1)$-principal bundles  on $M$.
Given an $U(1)$-principal bundle $E\in \cB U(1)(M)$ and a $Spin^{c}$-structure
$P\in Spin^{c}(V)$,  we can define a new $Spin^{c}$-structure
$E\otimes P\in Spin^{c}(V)$. 
A  formula for this tensor product is given by (\ref{actsglob}) specialised to the case $n=0$, see below.  This construction defines bifunctor
\begin{equation}\label{acts}\cB U(1)(M)\times Spin^{c}(V)\to Spin^{c}(V)\ .\end{equation}
If $Spin^{c}(V)$ is not empty, then the set of isomorphism classes of $Spin^{c}$-structures on $V$ is a torsor over the group of isomorphism classes in $\cB U(1)(M)$. Since the latter is canonically isomorphic to $H^{2}(M;\Z)$ we get a simply transitive
action of $H^{2}(M;\Z)$ on $Spin^{c}(V)/iso$.
Furthermore,  we get natural isomorphisms
$$C^{\infty}(M,U(1))\cong \Aut_{\cB U(1) (M)}(E)\cong \Aut_{Spin^{c}(V)}(E\otimes P)\cong \Aut_{Spin^{c}(V)}(P)\ .$$

The sum of two vector bundles with $Spin^{c}$-structures has a naturally induced  $Spin^{c}$-structure. This is formalised
 as the natural bifunctor
\begin{equation}\label{actsv} Spin^{c}(V)\times Spin^{c}(U)\to  Spin^{c}(V\oplus U)\ .\end{equation}
On the level of objects this bifunctor  is given by 
$$(P,Q)\mapsto P\otimes Q \ ,$$
 where the $Spin^{c}(n+m)$-principal bundle
\begin{equation}\label{actsglob}P\otimes Q:=(P\times_{M} Q)\times_{(Spin^{c}(n)\times Spin^{c}(m))}Spin^{c}(n+m)\end{equation}
is obtained from the 
$Spin^{c}(n)\times Spin^{c}(m)$-principal bundle $P\times_{M}Q$ by extension of structure groups along the upper horizontal map in the diagram
$$\xymatrix{Spin^{c}(n)\times Spin^{c}(m)\ar[d]\ar[r]&Spin^{c}(n+m)\ar[d]\\SO(n)\times SO(m)\ar[r]&SO(n+m)}\ .$$
Here $n=\dim(V)$ and $m=\dim(U)$, and  the compatibility with the lower part of this diagram is used to define the structure map
$\kappa_{P\otimes Q}$ from
$\kappa_{P}$ and $\kappa_{Q}$.
The bifunctor comes equipped with natural associativity constraints. We omit the details of the latter two aspects. 

We set $Spin^{c}(0):=U(1)$ and let $0_{M}$ denote the zero dimensional vector bundle on $M$.
 Then we get an identification of
$Spin^{c}(0_{M})\cong \cB U(1)(M)$, and for $n=0$  the bifunctor (\ref{actsv})
specialises to (\ref{acts}).
As a consequence of associativity the bifunctor (\ref{actsv}) is compatible with the action (\ref{acts}) of $\cB U(1)(M)$ in the sense that for $E\in \cB U(1)(M)$ have natural isomorphisms
\begin{equation}\label{actsglob1}(E\otimes P)\otimes Q\cong E\otimes (P\otimes Q)\cong P\otimes (E\otimes Q)\ .\end{equation}

The trivialised vector  bundle $M\times \R^{n}$ has a preferred trivial $Spin^{c}$-structure
$Q(n):=M\times Spin^{c}(n)$. We can use this to produce a canonical equivalence of groupoids
$$Spin^{c}(V)\cong Spin^{c}(V\oplus (M\times \R^{n}))\ ,  \quad P\mapsto P\otimes Q(n)\ .$$
On the level of $Spin^{c}$-structures we speak of stabilisations.

Let us now consider a pair $(M,f)$ of a compact   $n$-dimensional Riemannian manifold   and a map $f:M\to B$ which admits a refinement to a stable normal $B$-structure.
  Then we can assume that $f$ has a factorisation up to homotopy over $BSpin^{c}(k)$ as in the diagram
\begin{equation}\label{fakdia}\xymatrix{&&\xi_{k}^{Spin^{c}}\ar[r]\ar[d]&\xi_{k}^{SO}\ar[r]\ar[d]&\xi_{k}\ar[d]\\&&BSpin^{c}(k)\ar[r]\ar[d]&BSO(k)\ar[d]\ar[r]&BO(k)\ar[d]\\M\ar@{.>}[urr]^{\tilde f}  
\ar[r]^{f}&B\ar[r]&BSpin^{c}\ar[r]&BSO\ar[r]&BO}\ .\end{equation} The map $\tilde f$ classifies a
$Spin^{c}(k)$-principal bundle $\tilde f^{*}Q_{k}\to M$, where $Q_{k}\to BSpin^{c}(k)$ denotes the universal $Spin^{c}(k)$-bundle. Note that we have $\tilde f^{*}Q_{k}\in Spin^{c}( \tilde f^{*}\xi^{Spin^{c}}_{k})$. 

We let $\hat f:M\to BO(k)$ be induced by $\tilde f$ so that $\hat f^{*}\xi_k\cong \tilde f^{*}\xi_{k}^{Spin^{c}} $. With these identifications
the trivialisation (\ref{eq1360}) induces a bifunctor (\ref{actsv})
$$Spin^{c}(TM)\times Spin^{c}(\tilde f^{*}\xi^{Spin^{c}}_{k})\to Spin^{c}(M\times \R^{n+k})\ .$$

Since $\cB U(1)(M)$ acts simply transitively on isomorphisms classes
we conclude using (\ref{actsglob1}) that there is a unique isomorphism class of $Spin^{c}$-structures
$P\in Spin^{c}(TM)$ such that
\begin{equation}\label{eq31367}
P\otimes \tilde f^{*}Q_{k}\cong Q(n+k)\ .\end{equation} One can further check that this
isomorphism class only depends on the normal $B$-structure represented by $f$ and not on its representative. This is the tangential $Spin^{c}$-structure determined by the normal $Spin^{c}$-structure.

For constructions which involve glueing or in the notion of a $Spin^{c}$-map we need a rigidified notion of a tangential $Spin^{c}$-structure.  
\begin{ddd}\label{tang}
Assume that we have fixed a  normal $B$-structure in terms of the factorisation $\tilde f$ and the isomorphism  (\ref{eq1360}). Then we define a tangential  $Spin^{c}$-structure as a pair of a 
$Spin^{c}$-structure $P\in Spin^{c}(TM)$ together with a choice of an isomorphism in (\ref{eq31367}).
\end{ddd}
There are many tangential $Spin^{c}$-structures  associated to a  normal $Spin^{c}$-structure, but the main point is that two of them are isomorphic by a canonical isomorphism.

Let $h:M\to W$ be a smooth map and assume that we are given oriented euclidean vector bundles
$V_{M}\to M$ and $V_{W}\to W$ together with an isomorphism
\begin{equation}\label{eq31367okl}
V_{M}\oplus (M\times \R^{k})\cong h^{*}V_{W}\oplus (M\times \R^{l})\ .\end{equation}
Assume further that we are given $Spin^{c}$-structures
$P_{M}\in Spin^{c}(V_{M})$ and $P_{W}\in Spin^{c}(V_{W})$.
\begin{ddd}
A   a refinement of $h$ to a $Spin^{c}$-map is  a choice of an isomorphism
\begin{equation}\label{eq31367okll} h^{*}P_{W}\otimes Q(l)\cong P_{M}\otimes Q(k)\end{equation} in
$Spin^{c}(V_{M}\oplus (M\times \R^{k}))$ (this uses (\ref{eq31367okl})). 
The equivalence class of a $Spin^{c}$-map under stabilisation will be called a stable $Spin^{c}$-map.
 \end{ddd} 
Being a $Spin^{c}$-map is an additional datum, not just a property of the map. 
Observe that we can compose stable $Spin^{c}$-maps in a natural way.

We now assume that $(W,F)$ is a zero bordism of $(M,f)$. 
We choose a representative of the normal $B$-structure on $W$ involving the factorisation
$\tilde F:W\to BSpin^{c}(k)$. On $M$ we take the induced factorisation $\tilde f:=\tilde F_{|M}$.
Then we have a natural   decomposition of oriented euclidean vector bundles
\begin{equation}\label{eq1534}TW_{|M}\cong TM\oplus (M\times \R)\ ,\end{equation} where we trivialise the 
normal bundle by the outgoing unit normal vector field. Assume now that we have chosen tangential $Spin^{c}$-structures $P(TM)$ and $P(TW)$  on $M$ and $W$, respectively (see Definition \ref{tang}).  
In this situation we get a  natural refinement of the inclusion
$M\to W$ to a $Spin^{c}$-map. This refinement is distinguished by the condition that the following diagram in
$Spin^{c}(M\times \R^{n+1+k})$ 
$$\xymatrix{P(TM)\otimes Q(1)\otimes \tilde f^{*}Q_{k}\ar[r]^{\cong}\ar[d]^{\cong}&P(TW)_{|M}\otimes  \tilde f^{*}Q_{k}\ar[d]^{\cong}\\
Q(n+1+k)\ar@{=}[r]&Q(n+1+k)}$$
commutes. Here the upper corners are interpreted in 
$Spin^{c}(M\times \R^{n+1+k})$ using the   normal $B$-structures on $M$ or $W$, respectively. The vertical morphisms are given by the tangential   $Spin^{c}$-structures. Finally, the upper horizontal isomorphism uses  (\ref{eq1534}) and fixes the refinement of the inclusion $M\to W$ to a $Spin^{c}$-map.

\subsection{The definition of $\eta^{an}$}\label{sec81}

We consider a cycle $(M,f)$ for a  class $x=[M,f]\in \pi_{n}(MB)$ and assume in addition that $x$ is  torsion. Then there exists a non-zero integer $l\in \nat$ such that $l x=0$. We can thus find a zero bordism $(W,F)$ of the disjoint union of $l$ copies of $(M,f)$ which we denote by $l(M,f)$ .

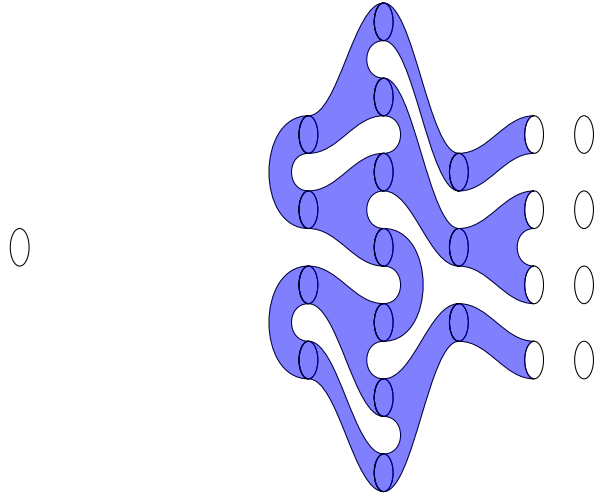
 \begin{center}\begin{align*}
 \begin{graphcalc}
 \gvec{1}
 \gsone
 \end{graphcalc}\hspace{3cm}
 \begin{graphcalc}
    \gvec{2}
   \ggamma
   \gtensor
   \ggamma
   \gcomp
      \gdelta    
      \gtensor
   \gdelta\gtensor\gdelta\gtensor \gid
   \gcomp
   \gid\gtensor \gmu\gtensor \gbeta\gtensor \gmu\gcomp
   \gid\gtensor\gdelta\gtensor\gid
    \end{graphcalc}\hspace{0.4cm}
     \begin{graphcalc}
 \gvec{4}
 \gsone\gtensor
 \gsone\gtensor
\gsone\gtensor
\gsone
 \end{graphcalc}
    \end{align*}
\end{center}
  \begin{center}A picture of $M$ and the zero bordism $W$ of $4M$\end{center}

 We will define
$\eta^{an}(x)\in Q_{n}(B)$ in terms  of a collection of indices  of associated $\Z/l\Z$-index problems \cite{MR1144425}. In order to formulate these index problems and to express the indices in terms of geometric and spectral invariants we must choose appropriate geometric structures.

We choose  
 a tangential $Spin^{c}$-structure   $(P,\kappa)\in Spin^{c}(TM)$ related to  the stable normal $Spin^{c}$-structure on $f$.
  A connection $\tilde \nabla^{TM}$ on $P$ induces via $\kappa$ a connection on $TM$. We say that $\tilde \nabla^{TM}$ is a $Spin^{c}$-extension of the Levi-Civita connection on $M$ if it induces the Levi-Civita connection $\nabla^{TM,LC}$ on $TM$.

The group $Spin^{c}(n)$ has a distinguished unitary  representation called the spinor representation $\Delta^{n}$. For even $n$ its dimension is $2^{n/2}$, and it has a decomposition $\Delta^{n}\cong \Delta^{n,+}\oplus \Delta^{n,-}$. It is related with the odd-dimensional case by  $\Delta^{n,+}_{|Spin^{c}(n-1)}\cong \Delta^{n-1}$.

The bundle $S(TM):=P\times_{Spin^{c}(n)}\Delta^{n}\to M$ is called the Spinor bundle of $M$.
If we have chosen a $Spin^{c}$-extension $\tilde \nabla^{TM}$ of the Levi-Civita connection on $M$, then the spinor bundle carries the structure of a Dirac bundle. We thus obtain the  $Spin^{c}$-Dirac operator $\Dirac_{M}$ which acts on sections of $S(TM)$. Standard references for these constructions are \cite[Ch. 3]{MR2273508},  \cite[App. D]{MR1031992}.

 If we are given a  class $\phi\in K^{0}(B )$, then  we can choose a $\Z/2\Z$-graded
vector bundle $V\to M$ whose $K$-theory class satisfies $[V]=f^{*}\phi\in K^{0}(M )$.
We choose a hermitean metric $h^{V}$ and a metric connection $\nabla^{V}$ which preserve the grading.  The triple $\bV:=(V,h^{V},\nabla^{V})$ will then be called a geometric vector bundle.
We let $\Dirac_{M}\otimes \bV$ be the Dirac operator twisted by $\bV$. It acts on sections of $S(TM)\otimes V$.

We now assume that $n=\dim(M)$ is odd.  The  $\eta$-invariant \cite{MR0397797}
of the twisted Dirac operator 
 $$\eta(\Dirac_{M}\otimes \bV)\in \R$$  is defined  as the value at $s=0$ of the meromorphic continuation of the $\eta$-function function $$\eta(\Dirac_{M}\otimes \bV)(s):=\Tr_{s} \: |\Dirac_{M}\otimes \bV|^{-s}\sign(\Dirac_{M}\otimes \bV)\ ,$$
 where $\Tr_{s}$ is the super trace with respect to the grading of $V$. 
Note that the trace exists if $\Ree(s)>n$, and that the meromorphic continuation of the $\eta$-function is regular at $s=0$ by the results of  \cite{MR0397797}.
 The $\eta$-invariant depends on the geometry of $M$ and $V$ in a possibly discontinuous way with jumps when eigenvalues of $\Dirac_{M}\otimes \bV$ cross zero. In order to get a quantity which depends continuously on the geometry one usually considers the reduced $\eta$-invariant for which in the present paper we will use the symbol $\xi$: 
\begin{equation}\label{eq132}\xi(\Dirac_{M}\otimes \bV):=[\frac{\eta(\Dirac_{M}\otimes \bV)+\dim\ker(\Dirac_{M}\otimes \bV)}{2}]\in \R/\Z .\end{equation}

In an appropriate model of $\pi_{n}(MB\Q/\Z)$ the zero bordism $(W,F)$  of $l(M,f)$   geometrically represents the lift of $x$ to  a class $$\hat x=[W,F]\in \pi_{n+1}(MB\Q/\Z)\ ,$$ using the notation of the diagram (\ref{eq3}).
We refer to Lemma  \ref{sec25} for more details.
We choose a  $Spin^{c}$-connection $\tilde \nabla^{TW}$ extending the 
 Levi-Civita connection  on $TW$   which extends the connection on the boundary of $W$  induced by $\tilde \nabla^{TM}$. 
The class $F^{*}\phi\in  K^{0}(W)$ extends the class $F_{|\partial W}^{*}\phi\in K^{0}(\partial W)$ which restricts to $f^{*}\phi\in K^{0}(M)$ on the copies of $M$ in the boundary of $W$. Hence  we can assume, after  adding some trivial bundles to the even and odd parts of $V$,  that the bundle on $\partial W$ induced by $V$  has an extension $U$ to $W$. We choose a hermitean metric $h^{U}$ and a metric connection $\nabla^{U}$ on $U$ which extend the  already given data on the boundary. In this way we get a geometric bundle $\bU:=(U,h^{U},\nabla^{U})$.

 We can now form the Atiyah-Patodi-Singer
boundary value problem for $\Dirac_{W}\otimes \bU$.
The analytic details of that boundary value problem are not important for our present purpose so that we refer to \cite{MR0397797} for a precise description. We only have to know that it produces a Fredholm operator $(\Dirac_{W}\otimes \bU)_{APS}$ which has a well-defined index
$$\ind((\Dirac_{W}\otimes \bU)_{APS})\in \Z\ ,$$ and that the following index formula
proved in \cite{MR0397797} holds true: 
\begin{equation}\label{eq4}\ind((\Dirac_{W}\otimes \bU)_{APS})=\int_{W} p_{n+1}(\Td(\tilde \nabla^{TW})\wedge \ch(\nabla^{U}))-l \frac{\eta(\Dirac_{M}\otimes \bV)+\dim\ker(\Dirac_{M}\otimes \bV)}{2}\ .\end{equation}
In this formula the  closed form $\Td(\tilde \nabla^{TW})\in \Omega^{0}(W)[b,b^{-1}]$ is the Chern-Weyl representative determined by the universal class $\Td\in HP\Q^{0}(BSpin^{c})$ and the connection $\tilde \nabla^{TM}$.
Similarly, the form $\ch(\nabla^U)\in \Omega^{0}(W)[b,b^{-1}]$ is the Chern-Weyl representative
determined by the class $\ch\in HP\Q^{0}(BU)$ and the connection $\nabla^{U}$. 
Note that we use powers of $b$ in order  to shift the higher form-degree components to total degree zero.

We consider the element
\begin{equation}\label{eq870v}e:=[\frac{\ind(\Dirac_{W}\otimes \bU)_{APS}}{l}]\in \Q/\Z\ .\end{equation}
Equivalently, by the index theorem  (\ref{eq4}) and (\ref{eq132}) we can write
\begin{equation}\label{eq870}e=[\frac{1}{l}\int_{W} p_{n+1}(\Td(\tilde \nabla^{TW})\wedge \ch(\nabla^{U}))]-\xi(\Dirac_{M}\otimes \bV)\end{equation} 
if we interpret this equality in $\R/\Z$.  
The quantity $e$ can be interpreted as a $\Z/l\Z$-index in the sense of \cite{MR1144425}.
In the following proposition we state how the number $e$ depends on the data.
\begin{prop}
\begin{enumerate}
\item
The value of $e$ does not depend on the choices of the geometric structures on $M$ and $W$.
\item The value of  $e$ only  depends on the $K$-theory class $\phi$. This dependence is additive and determines an element  
$\tilde e\in \Hom^{cont}(K^{0}(B),\Q/\Z)$.
\item The class $[\tilde e]\in Q_{n}(MB)$ of this homomorphism (using the presentation \eqref{spacerep1}) does not depend on $l$ or the choice of the zero bordism of $(W,F)$.
\item The element  $[\tilde e]\in Q_{n}(MB)$ described in 3. only depends on the bordism class $x$. This dependence is additive so that we obtain a well-defined homomorphism
$$\eta^{an}:\pi_{n}(MB)_{tors}\to Q_{n}(MB)\ .$$ 
\end{enumerate}
\end{prop}
\proof
On the one hand, we have $e\in \frac{1}{l}\Z/\Z\subseteq \R/\Z$. On the other hand, we know that
the right-hand side of (\ref{eq870}) depends continuously on the geometric data.
This shows that $e$ does not depend on the geometric structures at all since two choices  of  geometric structures can be connected by a family. This proves 1.

The element $e$ depends additively on the bundle $U$ and therefore only depends on the class $[U]\in K^{0}(M)$. Its construction  therefore  induces a  homomorphism
$\tilde e\in \Hom(K^{0}(B ),  \Q/\Z)$. 
Since it factors over the restriction along the map $f:M\to B$ and $M$ is compact this homomorphism is continuous. We use the identification 
$$\Q/\Z\cong \pi_{0}(K\Q/\Z)\xrightarrow{b^{-\frac{n+1}{2}}} \pi_{n+1}(K\Q/\Z)$$
in order to interpret $\tilde e$ as an element of 
$  \Hom^{cont}(K^{0}(B ),  \pi_{n+1}(K\Q/\Z))$.
This shows 2.

 Assume that we have a second zero bordism $(W^{\prime},F^{\prime} )$ of $l^{\prime} (M,f)$ yielding a homomorphism $\tilde e^{\prime}\in  \Hom^{cont}(K^{0}(B ),  \Q/\Z)$.
Then by glueing along boundary components we can form the closed Riemannian $n+1$-dimensional   $B$-manifold $\tilde W:=l^{\prime} W\cup_{ll^{\prime}M} lW^{\prime}$ which comes with a map $\tilde F:\tilde W\to B$. The latter has a natural refinement to a stable normal $B$-structure which restricts to the given stable normal $B$-structures  on $W$ and $W^{\prime}$.

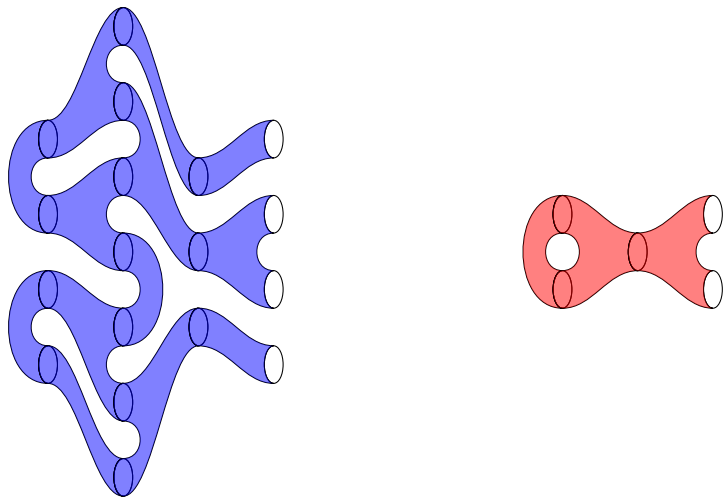
 \begin{center}\begin{align*}\begin{graphcalc}\def\gfillcolour{blue}
    \gvec{2}
   \ggamma
   \gtensor
   \ggamma
   \gcomp
      \gdelta    
      \gtensor
   \gdelta\gtensor\gdelta\gtensor \gid
   \gcomp
   \gid\gtensor \gmu\gtensor \gbeta\gtensor \gmu\gcomp
   \gid\gtensor\gdelta\gtensor\gid
    \end{graphcalc}    \hspace{3cm}
   \begin{graphcalc}\def\gfillcolour{red}
    \gvec{1}
\ggamma
\gcomp
\gmu\gcomp
\gdelta  
    \end{graphcalc}  \end{align*}
 \end{center}
 \begin{center}\begin{align*}\begin{graphcalc}\def\gfillcolour{blue}
    \gvec{4}
   \ggamma
   \gtensor
   \ggamma
    \gcomp
      \gdelta    
      \gtensor
   \gdelta\gtensor\gdelta\gtensor \gid
   \gcomp
   \gid\gtensor \gmu\gtensor \gbeta\gtensor \gmu\gcomp
   \gid\gtensor\gdelta\gtensor\gid
   \gcomp
  \def\gfillcolour{red} 
   \gmu\gtensor\gmu\gcomp
   \gdelta\gtensor \gdelta\gcomp
   \gbeta\gtensor \gbeta
    \end{graphcalc}     \end{align*}
 \end{center}
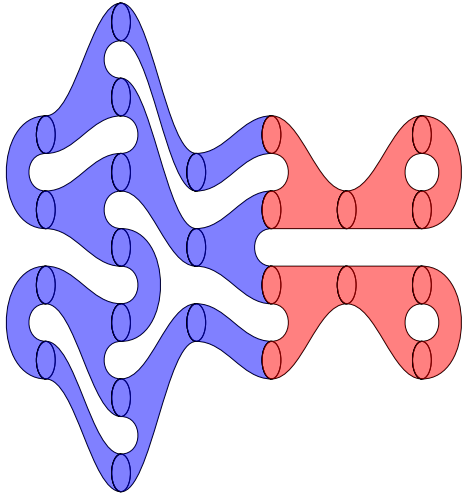
\begin{center} Pictures of $W$, $W^{\prime}$, and $\frac{1}{2}\tilde W$ with $l=4$ and $l^{\prime}=2$ \end{center}

Note that the tangential  $Spin^{c}$-structures
$(P,\kappa)$ and $(P^{\prime},\kappa^{\prime})$ come with isomorphisms
of the type (\ref{eq31367}). Compatibility with these fixes the morphism which we have to use order to glue $P$ with $P^{\prime}$. In this way we get a tangential  $Spin^{c}$-structure on $\tilde W$.
The triple $(\tilde W,\tilde F)$  is thus  cycle for a class $$y:=[\tilde W,\tilde F]\in \pi_{n+1}(MB)\ .$$
Then for $\phi\in K^{0}(B)$ we get from the right-hand side of (\ref{eq870}) that
$$\tilde e(\phi)-\tilde e^{\prime}(\phi)=[\frac{1}{ll^{\prime}} \langle \Td(T\tilde W)\cup \tilde F^{*}\ch(\phi),[\tilde W]\rangle]\ .$$
Since $\tilde F^{*}\Td^{-1}=\Td(T\tilde W)$ this is exactly the formula \eqref{localthomformula} for the evaluation of  $\epsilon(\frac{1}{ll^{\prime}}y)\in \pi_{n+1}(K\wedge MB\Q)$ against
$\Thom^{K}(\phi)\in K^{0}(MB )$. Therefore the class $[\tilde e]\in Q_{n}(MX)$ is independent of the choice of $l$ and the zero bordism $(W,F)$. This finishes the verification of 3.

We observe that the map
which associates to $(M,f)$ the class 
$[\tilde e]\in Q_{n}(MB)$ is additive under disjoint unions. Moreover, if
$(M,f)$ itself is zero bordant, i.e. we can find $(W,F)$ as above with $l=1$, then
$[\tilde e]=0$. It follows that the construction above uniquely descends to a homomorphism
\begin{equation}\label{etaandef111}\eta^{an}:\pi_{n}(MB)_{tors}\to Q_{n}(MB)\ .\end{equation}
\hB

Let us collect the essentials of this construction in the following definition.
\begin{ddd}\label{etaandef}
We define $\eta^{an}:=0$ for even $n$.
For odd $n$ we define the homomorphism 
$$\eta^{an}:\pi_{n}(MB)_{tors}\to Q_{n}(MB)$$ by the following prescription:
If $x\in \pi_{n}(MB)_{tors}$ is represented by $(M,f)$, then we choose a zero bordism
$(W,F)$ of $l(M,f)$ for a suitable non-zero $l\in \nat$. 
If $\phi\in K^{0}(B )$, then  we choose a geometric bundle $\bU\to W$ whose underlying $K$-theory class is $F^{*}\phi \in K^{0}(W)$ and whose restrictions to
 the $l$ copies of $M$ in the boundary are pairwise isomorphic.
 
We further choose a $Spin^{c}$-geometry for $W$ such that the restrictions to the $l$ copies of $M$ in the boundary of $W$ are again pairwise isomorphic. Then $\eta^{an}(x)\in Q_{n}(MB)$ is represented by the homomorphism
\begin{equation}\label{eq133}K^{0}(B )\ni \phi\mapsto [\frac{1}{l} \ind((\Dirac_{W}\otimes \bU)_{APS})]\in \Q/\Z\cong \pi_{n+1}(K\Q/\Z)\ .\end{equation}
\end{ddd}

\subsection{The secondary index theorem}\label{theproof}

In the Definitions \ref{etatopdefv} and  \ref{etaandef}
we have described  homomorphisms
$$\eta^{top}:\pi_{n}(MB )_{tors}\to Q_{n}(MB)\ ,\quad \eta^{an}:\pi_{n}(MB )_{tors}\to Q_{n}(MB)\ .$$ 
Both constructions follow a common idea. Given a torsion element
$x\in \pi_{n}(MB)_{tors}$ in a first step a lift $\hat x\in \pi_{n}(MB\Q/\Z)$, respectively a geometric representative of such a lift, is chosen. The homotopy theoretic invariant $\eta^{top}(x)$ is represented by the homomorphism
$K^{0}(MB)\to \Q/\Z$ induced by this lift via the homotopy theoretic pairing between $K$-homology and cohomology. The analytic variant $\eta^{an}(x)$ is represented by a  homomorphisms, which this time is obtained  from
a suitable family of Atiyah-Patodi-Singer index problems on the geometric representative of the lift $\hat x$. Because of these coincidences it is very natural to expect that the following theorem holds true.
\begin{theorem}\label{indthm}
$$\eta^{an}=\eta^{top}\ .$$
\end{theorem}
\proof 
An obvious option is to apply the $\Z/l\Z$-index theorem \cite{MR1144425} directly to $\eta^{an}$ in order to express it in homotopy theoretic terms. 
In this paper we decided to go a different path. It is interesting since it  explains in greater detail in which sense the homotopy theoretic construction of $\eta^{top}$ and the geometric or analytic constructions
involved in $\eta^{an}$  correspond to each other.  Our bride between analysis and topology will be
the identification of homotopy theoretic $K$-homology with the analytic picture \cite{MR679698} and 
the ordinary Atiyah-Singer index theorem for elliptic operators \cite{MR0236950}, respectively its local form described in \cite[Ch. IV]{MR2273508}. 

Some ideas of our proof of Theorem 
\ref{indthm}, in particular  the usage of Moore spaces, are taken from  \cite{MR1144425} and the proof of the $\R/\Z$-index theorem \cite[Thm 5.3]{MR0397799}.

We start with a description of Moore spaces for cyclic groups $\Z/l\Z$, $l\in \Z$. 
Let $S^{1}\to S^{1}$ be the $l$-fold covering of the pointed circle.
Its mapping cylinder $Z_{l}$ and mapping cone $C_{l}$ fit into the cofibre sequence of pointed spaces
\begin{equation}\label{eq11}S^{1}\to Z_{l}\to C_{l}\stackrel{\partial}{\to} \Sigma S^{1}\to \dots\ .\end{equation}
Note that the shifted suspension spectrum $\Sigma^{\infty-1}C_{l}$ is then a model for the Moore spectrum
$\tM \Z/l\Z$  discussed in \ref{sec141}. Further note that the inclusion of the cylinder basis $S^{1}\to Z_{l}$ is a homotopy equivalence. Hence  we have equivalences
$\Sigma^{\infty-1}Z_{l}\cong  \Sigma^{\infty-1}S^{1}\cong \tM\Z$.
Applying the functor $\Sigma^{\infty-1}$ to the sequence (\ref{eq11}) and using these identifications we get
the fibre sequence
\begin{equation}\label{eq12}\tM\Z\stackrel{l}{\to} \tM\Z \to \tM \Z/l\Z\stackrel{\partial}{\to}\Sigma \tM\Z\end{equation}
of Moore spectra.
We use the Moore spectra $\tM\Z/l\Z$ and the sequence (\ref{eq12}) as approximations for
$\tM\Q/\Z$  by spectra with finite skeleta and (\ref{eq631}) in the sense that
$$\tM \Q/\Z\cong \hocolim_{l} \:\tM \Z/l\Z\ .$$
The connecting maps are fixed by their compatibility with the inclusions \begin{equation}\label{eq4333}\Z/l\Z\to \Q/\Z\ ,\quad [n]\mapsto [\frac{n}{l}]\ .\end{equation}

Smashing the sequence (\ref{eq12}) with $MB$ and taking homotopy groups we get a long exact sequence of abelian groups
\begin{equation}\label{eq14}\dots\to\pi_{n+1}(MB\Z/l\Z)\stackrel{\partial}{\to} \pi_{n}(MB)\stackrel{l}{\to}\pi_{n}(MB)\to\pi_{n}(MB\Z/l\Z)\to\dots\ .\end{equation}

 The spectrum $MB\wedge X$ is related with Thom spectra by 
a fibre sequence
$$MB\to M(B\times X) \stackrel{\iota}{\to} MB\wedge X\to \Sigma MB\ ,$$
where we use the structure map
$B\times X\stackrel{\pr}{\to} B\to BGL(\R)$ in order to define the Thom spectrum
$M(B\times X)$ and the map
$MB\to M(B\times X)$ is induced by the base point of $X$. We further have
an equivalence   of spectra
\begin{equation}\label{eq13}\Sigma MB\Z/l\Z\cong  MB\wedge C_{l}\ .\end{equation}
 
We 
 use 
the notation $(F,G)$ in order to write maps from $M$ to $B\times C_{l}$.
In the following we construct a cycle $(\tilde W,(\tilde F,\tilde G))$ for a class in
$\pi_{n+2}(M(B\times C_{l}))$
such that
$$\iota [\tilde W,(\tilde F,\tilde G)]=\hat x\in \pi_{n+1}(MB\Z/l\Z)$$
under the identification \eqref{eq13}, where
$\tilde F:\tilde W\to B$ is the underlying map of a $B$-structure   and $\tilde G :\tilde W\to C_{l}$.
We will obtain   $(\tilde W,(\tilde F,\tilde G))$ by closing up the boundary of the triple $(W,(F,G))$
found in \ref{sec81}. 
 
 \begin{center}
 \begin{graphcalc}
    \gvec{2}
   \ggamma
   \gtensor
   \ggamma
   \gcomp
      \gdelta    
      \gtensor
   \gdelta\gtensor\gdelta\gtensor \gid
   \gcomp
   \gid\gtensor \gmu\gtensor \gbeta\gtensor \gmu\gcomp
   \gid\gtensor\gdelta\gtensor\gid
    \end{graphcalc}
\end{center}
  \begin{center}A picture of $S^{1}\times W$\end{center}
 The details are as follows.
We consider a two-sphere $S^{2}_{l}$ with $l$ holes.
  \begin{center}\begin{graphcalc}
    \gvec{2}
    \def\gfillcolour{red}
    \gmu
    \gtensor
    \gmu
    \gcomp
    \gmu
    \gcomp
    \gepsilon
  \end{graphcalc}
\end{center} 
 \begin{center}A picture of $S^{2}_{4}\times M$\end{center}
More precisely
we let $S^{2}_{l}\subset S^{2}$ be the compact submanifold with boundary
$\partial S^{2}_{l}\cong \bigsqcup_{i=1}^{l}S^{1}$ obtained by deleting the interiors of $l$ disjoint discs from $S^{2}$. We equip $S^{2}_{l}$ with a Riemannian metric with product structure such that all boundary components are isometric to the standard $S^{1}$.
The identification of the boundary with the $l$ copies of $S^{1}$ is fixed such that it preserves the natural orientations.
We now have an   identification
$$\partial (S^{1}\times W)\cong l(S^{1}\times M)\cong \partial (S^{2}_{l}\times M)\ .$$ We let
\begin{equation}\label{eq15cvc}\tilde W:=(S^{1}\times W)\cup_{l(S^{1}\times M)} (S^{2}_{l}\times M)\end{equation} 
be the manifold obtained by glueing along the boundary.
\begin{center}
\begin{align*}
\begin{graphcalc}
    \gvec{2}
   \ggamma
   \gtensor
   \ggamma
   \gcomp
      \gdelta    
      \gtensor
   \gdelta\gtensor\gdelta\gtensor \gid
   \gcomp
   \gid\gtensor \gmu\gtensor \gbeta\gtensor \gmu\gcomp
   \gid\gtensor\gdelta\gtensor\gid
        \gcomp
       \def\gfillcolour{red}
    \gmu
    \gtensor
    \gmu
    \gcomp
    \gmu
    \gcomp
    \gepsilon
 \end{graphcalc}
 \end{align*}
\end{center}
 \begin{center}A picture of $\tilde W$\end{center}

We define $\tilde F:\tilde W\to B$ such that it restricts to $$S^{1}\times W\stackrel{\pr_{W}}{\to} W\stackrel{F}{\to} B\ ,\quad 
S^{2}_{l}\times M\stackrel{\pr_{M}}{\to} M\stackrel{f}{\to} B\ .$$
We must refine the map $\tilde F$ to a normal $B$-structure.
We start with the usual  normal framing $$TS^{2}\oplus S^{2}\times \R\cong S^{2}\times \R^{3}$$ of $S^{2}$. If we take $k=1$ and let $f$ and  $\hat f$  in \eqref{eq1360} be the constant maps,
then we can interpret this isomorphism as a normal $B$-structure, too.
By restriction we obtain a normal $B$-structure on $S^{2}_{l}$. Furthermore, the construction explained at the end of \ref{sec1400}
provides   normal $B$-structures 
$$TS^{1}\oplus S^{1}\times \R^{2}\cong S^{1}\times \R^{3}$$
on the $l$ copies of $S^{1}$ in the boundary of  $S^{2}_{l}$ which are isomorphic to each other.
We are given a normal $B$-structure
$$TW\oplus \hat F^{*}\xi_{k}\cong W\times \R^{n+1+k}\ .$$
We therefore get an
induced normal $B$-structure  on the product $S^{1}\times W$ refining $\tilde F_{|S^{1}\times W}$:
$$T(S^{1}\times W)\oplus \hat F^{\prime,*}\xi_{k+2} \cong TS^{1}\oplus  TW\oplus S^{1}\times W\times \R^{2} \oplus (\hat F^{*}\circ \pr_{W})^{*}\xi_{k} \cong S^{1}\times W\times R^{n+k+4}\ ,$$
where $\hat F^{\prime}$ is the two-fold stabilisation of $\hat F\circ \pr_{W}$.
In a similar manner, using the induced normal $B$ structures on the copies of $M$ in the boundary of $W$ we get a normal $B$-structure 
$$T(S^{2}_{l}\times M)\oplus \hat f^{*}\xi_{k+1}\cong S^{2}_{l}\times M\times \R^{n+k+4}$$ on
the product $S^{2}_{l}\times M$ which refines $\tilde F_{|S^{2}_{l}\times M}$.
These isomorphisms coincide over  the locus of glueing $l(S^{1}\times M)$.
Hence we get a refinement of $\tilde F$ to a normal $B$-structure.


We now consider the map
\begin{equation}\label{eq15}S^{1}\times W\stackrel{\pr_{S^{1}}}{\to} S^{1}\stackrel{i}{\to} C_{l}\ ,\end{equation} 
where $i:S^{1}\to C_{l}$ is the identification of $S^{1}$ with the basis of the mapping cone.
Note that the map
$$\sqcup_{j=1}^{l}i :\partial S^{2}_{l}\cong \bigsqcup_{j=1}^{l }S^{1}\to C_{l}$$ can be extended
to a map 
\begin{equation}\label{eq51}g:S^{2}_{l}\to C_{l}\ .\end{equation} 
We can and will restrict the choice of  $g$ such that it is smooth on the preimage of a neighbourhood $U$ of  the cone basis  $\partial C_{l}\subset C_{l}$, and regular values of $g$ in the interior $U\setminus \partial C_{l}$ have exactly one preimage. 
The restriction of the map (\ref{eq15}) to $\partial (S^{1}\times W)\cong l(S^{1}\times M)$ thus  has an extension across the other part
$S^{2}_{l}\times M$ of $\tilde W$ given by 
$$S^{2}_{l}\times M\stackrel{\pr_{S^{2}_{l}}}{\to} S^{2}_{l}\stackrel{g}{\to} C_{l}\ ,$$
 Altogether we obtain the map
$\tilde  G :\tilde W\to C_{l}$. The cycle $(\tilde W,(\tilde F,\tilde G))$ represents a class in in 
$\pi_{n+2}(M(B\times C_{l}))$ and we consider its image
$$\hat x:= \iota[\tilde W,(\tilde F,\tilde G)]\in \pi_{n+2}(MB\wedge C_{l})\stackrel{(\ref{eq13})}{\cong} \pi_{n+1}(MB\Z/l\Z)\ .$$

Let $\partial:\pi_{n+1}(MB\Z/l\Z)\to \pi_{n}(MB)$
be  the boundary as in (\ref{eq14}). 
\begin{lem}\label{sec25}
We have $\partial \hat x=x$
\end{lem}
\proof  The  boundary operator $\partial$ in the Lemma is induced by 
the map denoted by the same symbol in (\ref{eq11})
$$\partial:C_{l}\stackrel{p}{\to} \Sigma S^{1}\cong S^{2}\ ,$$
where $p$ is the projection which contracts the cone basis to a point.
Therefore
$\partial \hat x\in \pi_{n+2}(MB\wedge S^{2})$ is represented by 
$(\tilde W,(\tilde F,p \circ \tilde G))$. We must show that it corresponds to
$x$ under the suspension isomorphism
$$\pi_{n}(MB)\cong \pi_{n+2}(MB\wedge S^{2})\ .$$
To this end we invert the suspension isomorphism in the geometric picture.
This inverse is of course given by taking the inverse image of  a regular point in $S^{2}$ of the corresponding component $p\circ \tilde G$ of the structure map. If we take the inverse image of
a point in the neighbourhood $U\setminus \partial C_{l}$ mentioned above we exactly recover the representative $(M,f)$ of $x$. \hB

The construction of $\eta^{top}$ involves the $K$-homology of a based space $Y$ defined homotopy theoretically as $\pi_{*}(K\wedge Y)$. It is equivalent to the analytic picture introduced in \cite{MR679698}. The analytic $K$-homology is subsumed in the more general bivariant $KK$-theory (see \cite{MR918241} and the text book \cite{MR1656031}) which allows to treat $K$-homology and cohomology on equal footing.
Of particular importance for our purpose is that the product in $KK$-theory provides a description of the $\cap$-product
between $K$-homology and $K$-theory which easily compares with the operation of twisting Dirac operators.

The unit of $K$-theory induces the map
\begin{equation}\label{eq20}\epsilon:\pi_{n+2}(MB\wedge C_{l})\to \pi_{n+2}(K\wedge MB\wedge C_{l})\ .\end{equation}
We use the Thom isomorphism for $MB$ in $K$-homology 
\begin{equation}\label{eq21}\Thom_{K}:\pi_{n+2}(K\wedge MB\wedge C_{l})\stackrel{\sim}{\to} \pi_{n+2}(K\wedge B_{+}\wedge C_{l})\ .\end{equation}
Finally we use $KK$-theory in order represent this  $K$-homology  of a pointed space analytically.
For the moment we assume that $X$ and $B$ are compact. This is no real restriction since we are calculating with a finite number of cycles at a time and their structure maps can only hit compact parts
of the spaces $B$ and $X$. For a compact based space $Y$ we let $C(Y)$ denote the $C^{*}$-algebra
of continuous $\C$-valued functions which vanish on the base point.
Then by the equivalence between homotopy theoretic and analytic $K$-homology \cite{MR679698} we have an isomorphism
 \begin{equation}\label{eq22}\pi_{n+2}(K\wedge B_{+}\wedge C_{l})\cong KK_{n+2}(C(B_{+}\wedge C_{l}),\C)\ . \end{equation}

The $Spin^{c}$-extension of the Levi-Civita connection on $W$ together with the standard $Spin^{c}$-geometry of $S^{1}$ induce a corresponding  product $Spin^{c}$-extension of the Levi-Civita connection on $S^{1}\times W$. The  $Spin^{c}$-geometry on $S^{1}$  also induces such a geometry on the boundary  $\partial S^{2}_{l}\cong lS^{1}$ which we extend to $S^{2}_{l}$, again with a product structure.
We get a corresponding product   $Spin^{c}$-extension of the Levi-Civita connection on $S^{2}_{l}\times M$. These geometric structures glue nicely and give a   $Spin^{c}$-extension of the Levi-Civita connection on  $\tilde W$. We let $\Dirac_{\tilde W}$ denote the corresponding Dirac operator. It acts on the complex spinor bundle $S(T\tilde W)$.
The Hilbert space $L^{2}(\tilde W,S(\tilde W))$ of square integrable sections of this bundle carries an action $\rho$ of the $C^{*}$-algebra
$C(\tilde W_{+})$ of continuous functions on $\tilde W$ by multiplication. The triple
$$(\Dirac_{\tilde W}):=(L^{2}(\tilde W,S(\tilde W)), \Dirac_{\tilde W},\rho)$$ is an unbounded Kasparov module for the pair of $C^{*}$-algebras $(C(\tilde W_{+}),\C)$ and represents a class
$$[\Dirac_{\tilde W}]\in KK_{n+2}(C(\tilde W_{+}),\C)\ .$$
The map $(\tilde F,\tilde G)$ induces a homomorphism of $C^{*}$-algebras
$$(\tilde F,\tilde G)^{*}:C(B_{+}\wedge C_{l})\to C(\tilde W_{+})$$ which in turn induces the push-forward in analytic $K$-homology in the statement of the following Lemma.

\begin{lem}\label{lem9991}
The image of the class $\hat x\in \pi_{n+2}(MB\wedge C_{l})$
under the composition of the unit (\ref{eq20}), Thom isomorphism, (\ref{eq21}) and the identification (\ref{eq22}) is given by 
 $$(\tilde F,\tilde G)_{*}[\Dirac_{\tilde W}]\in KK_{n+2}(C(B_{+}\wedge C_{l}),\C)\ .$$ 
 \end{lem}
 \proof
The image of $\hat x$ under the unit and Thom isomorphism is given 
$$\Thom_{K}(\epsilon(\hat x))=\beta(\Delta(\hat x)) \in \pi_{n+2}(K\wedge B_{+}\wedge C_{l})\ ,$$
where
$\Delta:MB\to MB\wedge B_{+}$ is the Thom diagonal and
$\beta:MB\to K$  
 is the $K$-orientation of $MB$ given by \eqref{betadef}.
We have $\Delta(\hat x)=\iota[\tilde W,(\tilde F,(\tilde F,\tilde G))]\in \pi_{n+2}(MB\wedge B_{+}\wedge C_{l})$.
Formally we can view this as the push-forward of the $B$-bordism fundamental class of $\tilde W$
along the map $(\tilde F,\tilde G)$. Its image under the $K$-orientation $\beta$ is then the push-forward of the $K$-theory fundamental class of $\tilde W$ associated to the $Spin^{c}$-structure along this map. In the analytic picture of $K$-homology  the  $K$-theory fundamental class of $\tilde W$ is represented by the $Spin^{c}$-Dirac operator. Hence  it is equal to $[\Dirac_{\tilde W}]$.
We thus get
$$ \beta(\iota[\tilde W,(\tilde F,(\tilde F,\tilde G))])= (\tilde F,\tilde G)_{*}[\Dirac_{\tilde W}]\ .$$
\hB

We let $\phi\in K^{0}(B)$. The pairing on  the left-hand side in the following calculation in $\pi_{n+2}(K\wedge C_{l})\cong \Z/l\Z$ is reminiscent to the evaluation occurring in the definition of $\eta^{top}$:
\begin{eqnarray*} \langle \Thom^{K}(\phi),\epsilon(\hat x)\rangle &=&\langle  \phi, \Thom_{K}(\epsilon(\hat x))\rangle \\&\stackrel{Lemma\: \ref{lem9991})}{=} &\langle \phi,(\tilde F,\tilde G)_{*}[\Dirac_{\tilde W}]\rangle \\ &=&\tilde G_{*}  ([\Dirac_{\tilde W}]\cap \tilde F^{*}\phi)\ .\end{eqnarray*} 

We choose a geometric bundle $\tilde \bV$ whose underlying $K$-theory class is equal to $\tilde F^{*}\phi$. The restriction of the map $\tilde F$  to the part $S^{1}\times W\subset \tilde W$ factors over the projection to $W$ and $F:W\to B$. Hence we can assume that the restriction of $\tilde \bV$ to $S^{1}\times W\subset \tilde W$
is isomorphic to the pull-back of the bundle $\bU$ on $W$, if we allow some stabilisation of $\tilde \bV$ and $\bU$.

In the $KK$-picture the $\cap$-product
$$[\Dirac_{\tilde W}]\cap \tilde F^{*}\phi\in KK_{n+2}(C(\tilde W_{+}),\C)$$ is realised by the unbounded Kasparov module
$(L^{2}(\tilde W,  S(\tilde W)\otimes \tilde \bV),\Dirac_{\tilde W}\otimes \tilde \bV,\rho)$
 associated to the twisted Dirac operator
$\Dirac_{\tilde W}\otimes \tilde \bV$, where $\rho$ again denotes the action of $C(\tilde W_{+})$ on  $L^{2}(\tilde W,  S(\tilde W)\otimes \tilde \bV)$
by multiplication. Hence
we have
$$[\Dirac_{\tilde W}\otimes \tilde \bV]=[\Dirac_{\tilde W}]\cap (\tilde F,\tilde G)^{*}\phi\ .$$
We conclude that
$\eta^{top}(x)\in Q_{n}(MB)$ is represented by the map
\begin{equation}\label{eq62}K^{0}(B)\ni \phi\mapsto \tilde G_{*}[\Dirac_{\tilde W}\otimes \tilde \bV]\in KK_{n+2}(C(C_{l}),\C)\cong \pi_{n+1}(K\Z/l\Z)\subset \pi_{n+1}(K\Q/\Z)\ ,\end{equation}
where the last inclusion is induced by (\ref{eq4333}).

Next we want to calculate the element in $\Z/l\Z$ given by $\tilde G_{*}([\Dirac_{\tilde W}\otimes \tilde \bV])$. Since the usual index theorem \cite{MR0236950} calculates integral indices we have to construct and calculate an integral representative of this $\Z/l\Z$-valued index. 
 The inclusion of the cone base
$i:S^{1}\to C_{l}$ induces a surjective map
$$\Z\cong \pi_{n+2}(K\wedge S^{1})\to \pi_{n+2}(K\wedge C_{l})\cong \Z/l\Z\ .$$
We try to construct a lift of
$\tilde G_{*}([\Dirac_{\tilde W}\otimes \tilde \bV])$ to $\pi_{n+2}(K\wedge S^{1})$
by providing a factorisation  $\gamma$ as in the diagram
$$\xymatrix{&S^{1}\ar[d]^{i}\\\tilde W\ar@{..>}[ur]^{\gamma}\ar[r]^{\tilde G}&C_{l}}\ .$$
For our given representative such a factorisation does not exist in general. The idea is to modify 
the representative  without changing its $\Z/l\Z$-valued index such that this lift exists for the modified cycle.

Note that $M$ is a closed odd-dimensional manifold. The Dirac operator $\Dirac_{M}\otimes \bV$ is 
selfadjoint. We can find a selfadjoint smoothing operator $Q$ on $L^{2} (M,S(M)\otimes V)$
such that $\Dirac_{M}\otimes \bV+Q$ is invertible.
In \cite{MR2191484} such a perturbation was called a taming.
As described in this reference a taming can be lifted to the product
$S^{2}_{l}\times M$ and also to a collar neighbourhood  $l(S^{1}\times (-\epsilon,0]\times M) \cong Z\subset S^{1}\times W$ of  $\partial (S^{1}\times  W)$. This lift is a selfadjoint operator $\bar Q$ on $L^{2}(Z\cup_{l(S^{1}\times M)} S^{2}_{l}\times M,S(\tilde W)\otimes \tilde V)$ which is an integral operator along $M$ and local in the remaining directions. Let $\chi:\tilde W\to [0,1]$ be a cut-off function which is supported on
$Z\cup_{l(S^{1}\times M)} S^{2}_{l}\times M$,  is  equal to one in a neighbourhood of the subset $S^{2}_{l}\times M$, and only depends on the normal variable near
$\partial (S^{1}\times W)$. We define the extension $\tilde Q:=\chi \bar Q \chi$ of $\bar Q$ to all of 
$\tilde W$. Note that $\tilde Q$ commutes with the image of $\tilde G^{*}(C(C_{l}))$.
Adding $\tilde Q$ to $\Dirac_{W}\otimes \tilde \bV$ gives a relatively compact perturbation.
Therefore we have
$$\tilde G_{*}[\Dirac_{\tilde W}\otimes \tilde \bV]= \tilde G_{*}[\Dirac_{\tilde W}\otimes \tilde \bV+\tilde Q]$$ in $KK_{n+2}(C(C_{l}),\C)$.
On the part $S^{2}_{l}\times M\subset \tilde W$ the perturbed operator $\Dirac_{W}\otimes \tilde \bV+\tilde Q$ is invertible along the fibres of the projection to $S^{2}_{l}$.

On $S^{1}\times (0,\infty)$ we fix a warped product Riemannian metric of the form $g:=dt^{2}+f(t) g^{TS^{1}}$, where $f_{|[0,1]}\equiv 1$ and $f(t)=t$ for $t\in(2,\infty)$. The precise form is not important, but we need that $\lim_{t\to \infty} f(t)=\infty$. We further choose a $Spin^{c}$-extension of this geometry.
We define the manifold
$$\tilde{\tilde W}:=S^{1}\times W\cup_{l(S^{1}\times M)}l(S^{1}\times  [0,\infty)\times M)\ .$$
Its $Spin^{c}$-geometry is given as the product of the geometries of $W$ and $S^{1}$ on the left-hand side, and by the product of the geometries on $lM$ and $S^{1}\times (0,\infty)$ on the right-hand side.
In a similar manner we define the geometric bundle $\tilde{\tilde \bV}$ on $\tilde{\tilde W}$ by a cylindrical extension of $\tilde \bV_{|S^{1}\times W}$. 
We define an operator
$\tilde{\tilde Q}$ similarly to $\tilde Q$ by lifting $Q$ to the cylinder $S^{1}\times [0,\infty)\times M$ and cutting off in the interior of $S^{1}\times W$.
Finally we let $\tilde {\tilde G}:\tilde{\tilde W}\to C_{l}$
be given by $G$ on $S^{1}\times W$ and the radially constant extension 
of $(\tilde G)_{|\partial (S^{1}\times W)}$ to the cylinder $l([0,\infty)\times S^{1}\times M)$.

To every complete   Riemannian manifold $(N,g)$ we associate the  commutative $C^{*}$-algebra 
$C_{g}(N)$ defined as the closure  in the $\sup$-norm of the algebra $C_{g}^{\infty}(N)$ of all bounded smooth functions $f$ on $N$ such that $|df|\in C_{0}(N)$. Note that if $N$ is compact, then
$C_{g}(N)=C(N)$.

The operator $\Dirac_{\tilde{\tilde W}}\otimes \tilde {\tilde \bV}+\tilde {\tilde Q}$ is invertible along the fibre
$M$ of the projection from the cylindrical end of $\tilde{\tilde W}$ to  $S^{1}\times  [0,\infty)$.
Therefore it is invertible at infinity in the sense of  \cite[Assumption 1]{MR1348799}. The arguments given \cite[Section 1]{MR1348799} show that
$(L^{2}(\tilde{\tilde W},S(T\tilde{\tilde W})\otimes \tilde{\tilde V}),\Dirac_{\tilde{\tilde W}}\otimes \tilde {\tilde \bV}+\tilde {\tilde Q}, \tilde{\tilde \rho})$ is an unbounded Kasparov module over  $C_{g}(\tilde{\tilde W})$. We write $[\Dirac_{\tilde{\tilde W}}\otimes \tilde {\tilde \bV}+\tilde {\tilde Q}]\in KK_{n+2}(C_{g}(\tilde{\tilde W}),\C)$ for its class.
 Because of the choice of the warped product metric we have the homomorphism of $C^{*}$-algebras
$\tilde{\tilde G}^{*}:C(C_{l})\to C_{g}(\tilde{\tilde W})$ so that the class
$\tilde{\tilde G}_{*}[\Dirac_{\tilde{\tilde W}}\otimes \tilde {\tilde \bV}+\tilde {\tilde Q}]\in KK_{n+2}(C(C_{l}),\C)$ is well-defined.

The operators  $\Dirac_{W}\otimes \tilde \bV+\tilde Q$ and $\Dirac_{\tilde{\tilde W}}\otimes \tilde {\tilde \bV}+\tilde {\tilde Q}$ coincide on $S^{1}\times W$ and are invertible along the fibres $M$ outside of this submanifold of $\tilde W$ and $\tilde{\tilde W}$. 

We let $N_{1}$ be the double of $l(M\times S^{1}\times [0,\infty))$ which carries an obvious geometric bundle $\bV_{1}$, taming $Q_{1}$, and admits a map $G_{1}:N_{1}\to C_{l}$. The associated $K$-theory class $[N_{1}]:=[\Dirac_{N_{1}}\otimes \bV_{1}+Q_{1}]\in   KK_{n+2}(C_{g}(N_{1}),\C)$  vanishes since
$\Dirac_{N_{1}}\otimes \bV_{1}+Q_{1}$ is invertible. Similarly we define
$N_{2}$, a map $G_{2}:N_{2}\to C_{l}$,  and a trivial class $[\Dirac_{N}\otimes \bV_{2}+Q_{2}]\in KK_{n+2}(C_{g}(N_{2}),\C)$ by attaching the (reflected) warped product $lM\times S^{1}\times [0,\infty)$
to $M\times S_{l}^{2}$.

In this situation, which is schematically pictured below, we can apply a relative index theorem (the proof of \cite[Theorem 1.14]{MR1348799} extends since all constructions there are compatible with the action of the $C^{*}$-algebra $C(C_{l})$) in order to get the equality $!$ in $KK_{n+2}(C(C_{l}),\C)$:
$$\tilde G_{*}[\Dirac_{\tilde W}\otimes \tilde \bV]=
\tilde G_{*}[\Dirac_{\tilde W}\otimes \tilde \bV]+G_{1,*}[N_{1}]\stackrel{!}{=} \tilde{\tilde  G}_{*}[\Dirac_{\tilde{\tilde W}}\otimes \tilde {\tilde \bV}+\tilde {\tilde Q}]+G_{2,*}[N_{2}]=\tilde{\tilde  G}_{*}[\Dirac_{\tilde{\tilde W}}\otimes \tilde {\tilde \bV}+\tilde {\tilde Q}]\ .
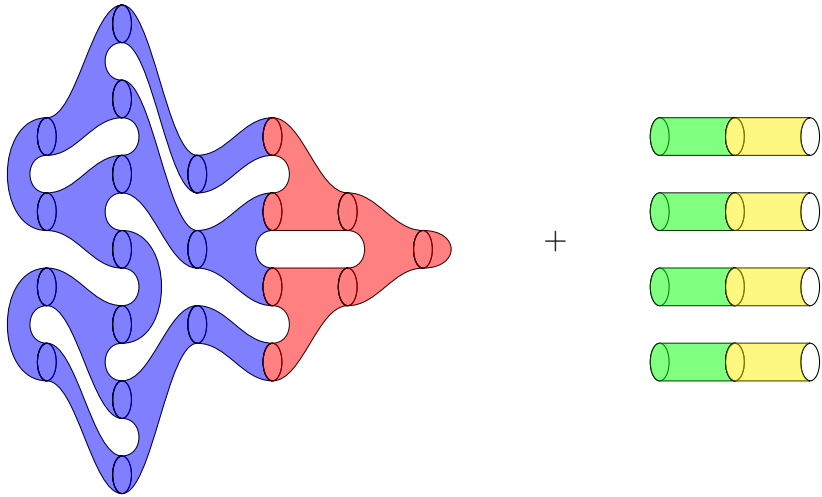
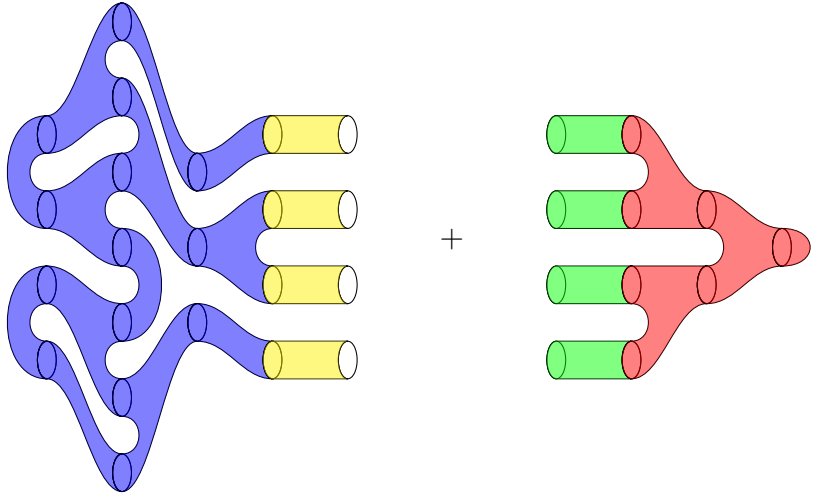
$$

\begin{center}
\begin{align*}
 \begin{graphcalc}
    \gvec{2}
   \ggamma
   \gtensor
   \ggamma
   \gcomp
      \gdelta    
      \gtensor
   \gdelta\gtensor\gdelta\gtensor \gid
   \gcomp
   \gid\gtensor \gmu\gtensor \gbeta\gtensor \gmu\gcomp
   \gid\gtensor\gdelta\gtensor\gid
        \gcomp
        \def\gfillcolour{red}
    \gmu
    \gtensor
    \gmu
    \gcomp
    \gmu
    \gcomp
    \gepsilon
 \end{graphcalc}  \hspace{1cm}+ \hspace{1cm}
  \begin{graphcalc}
    \gvec{4}
    \def\gfillcolour{green}
    \gid\gtensor
     \gid\gtensor
      \gid\gtensor
       \gid
       \gcomp
       \def\gfillcolour{yellow}
         \gid\gtensor
     \gid\gtensor
      \gid\gtensor
       \gid
 \end{graphcalc}
 \end{align*}
\end{center}

\begin{center}
\begin{align*}
 \begin{graphcalc}
    \gvec{2}
   \ggamma
   \gtensor
   \ggamma
   \gcomp
      \gdelta    
      \gtensor
   \gdelta\gtensor\gdelta\gtensor \gid
   \gcomp
   \gid\gtensor \gmu\gtensor \gbeta\gtensor \gmu\gcomp
   \gid\gtensor\gdelta\gtensor\gid
        \gcomp \def\gfillcolour{yellow}
         \gid\gtensor
     \gid\gtensor
      \gid\gtensor\gid
  \end{graphcalc}  \hspace{1cm}+ \hspace{1cm}
  \begin{graphcalc}
    \gvec{4}
    \def\gfillcolour{green}
    \gid\gtensor
     \gid\gtensor
      \gid\gtensor
       \gid
       \gcomp
        \def\gfillcolour{red}
    \gmu
    \gtensor
    \gmu
    \gcomp
    \gmu
    \gcomp
    \gepsilon
           \end{graphcalc}
 \end{align*}
\end{center}
 \begin{center}{\it A picture of the relative index theorem.
 The operator is invertible on the parts which are not blue.
 The index of the operator associated to  the upper picture is the index of its left part $\tilde W$. 
 The index is preserved under cut-and paste as indicated.
 The index of the operator associated to the lower picture is again the index of the left part $\tilde{\tilde W}$.
}  \end{center}

\newcommand{\bL}{\mathbf{L}}
Note the factorisation
$\tilde{\tilde G}:\tilde{\tilde W}\stackrel{\pr_{S^{1}}}{\to}S^{1}\stackrel{i}{\to} C_{l}$,
where the last map is the embedding of the cone basis.
Therefore
\begin{equation}\label{ine517271}\pr_{S^{1}*}[\Dirac_{\tilde{\tilde W}}\otimes \tilde {\tilde \bV}+\tilde {\tilde Q}]\in KK_{n+2}(C(S^{1}),\C)\cong \pi_{n+2}(K\wedge S^{1})\cong \Z\end{equation}
represents the desired integral lift. We reduce the calculation of this integer to a calculation of a Fredholm index by suspending once more.
Let $\bL$  be a geometric line bundle on $S^{1}\times S^{1}$ such that
  $c_{1}(L)\in H^{2}(S^{1}\times S^{1};\Z)$ is a generator and the sign is fixed such that
  \begin{equation}\label{hebfjefefiof}\ind(\Dirac_{S^{1}\times S^{1}}\otimes \bL)=1\ .\end{equation} We define $\tilde{\tilde{\tilde W}}:=S^{1}\times \tilde{\tilde W}$ with the product $Spin^{c}$-geometry. Then the desired integer \eqref{ine517271} is the index
of the Fredholm operator 
$\Dirac_{\tilde{\tilde{\tilde W}}}\otimes \pr_{\tilde{\tilde W}} \tilde{\tilde \bV}\otimes \pr_{S^{1}\times S^{1}}^{*} \bL+ \tilde{\tilde {\tilde Q}}$
where $\tilde{\tilde {\tilde Q}}$ is the taming induced by $\tilde {\tilde Q}$.
In order to see this note that the identification 
$KK_{1}(C(S^{1}),\C)\stackrel{\cong}{\to} \Z$ is given by the iterated Kasparov product
$x\mapsto [\bL]\otimes_{C(S^{1}\times S^{1})} ([\Dirac_{S^{1}}]\otimes_{\C} x)\in KK_{2}(\C,\C)\cong \Z$, 
where $[\Dirac_{S^{1}}]\in KK_{1}(C(S^{1}),\C)$ is the class of the standard $Spin^{c}$-Dirac operator on $S^{1}$ and
$[\bL]\in KK_{0}(\C,C(S^{1}\times S^{1}))$
is the class represented by the line bundle $\bL$.
 We now have $$\pr_{S^{1}\times S^{1}*}[\Dirac_{\tilde{\tilde{\tilde W}}}\otimes \pr_{\tilde{\tilde W}} \tilde{\tilde \bV} + \tilde{\tilde {\tilde Q}}]=[\Dirac_{S^{1}}]\otimes_{\C}   \pr_{S^{1}*}[\Dirac_{\tilde{\tilde W}}\otimes \tilde {\tilde \bV}+\tilde {\tilde Q}]$$
 in $KK_{2}(C(S^{1}\times S^{1}),\C)$ and therefore, using the relation of twisting Dirac operators and Kasparov products,
$$[
\Dirac_{\tilde{\tilde{\tilde W}}}\otimes \pr_{\tilde{\tilde W}} \tilde{\tilde \bV}\otimes \pr_{S^{1}\times S^{1}}^{*} \bL+ \tilde{\tilde {\tilde Q}}]=[\bL]\otimes_{C(S^{1}\times S^{1})} \pr_{S^{1}\times S^{1}*}[\Dirac_{\tilde{\tilde{\tilde W}}}\otimes \pr_{\tilde{\tilde W}} \tilde{\tilde \bV} + \tilde{\tilde {\tilde Q}}]$$
in $KK_{2}(\C,\C)$.

We deform the warped product metric on the end of $\tilde{\tilde {\tilde{W}}}$ to a product metric. This produces a continuous family of Fredholm operators and therefore does not change the index. After this deformation we see that
 $\tilde{\tilde {\tilde{W}}}\cong S^{1}\times S^{1}\times \hat W$ geometrically, where
$\hat W:=W\cup_{\partial W} l([0,\infty)\times M)$  and  $\hat W$  carries the  cylindrical extension of the geometry of $W$. Similarly we let $\hat \bU$ be the geometric bundle on $\hat W$ obtained by the cylindrical extension of $\bU$. Then the resulting   operator 
represents the product
$$[\Dirac_{\hat W}\otimes \hat \bU+\hat Q]\otimes_{\C}[\Dirac_{S^{1}\times S^{1}}\otimes \bL]\ ,$$
where  $\hat Q$  is the  taming   $\hat W$ uniquely determined by the property that its lift  is  the taming $\tilde{\tilde{\tilde Q}}$. Because of \eqref{hebfjefefiof}
 the  integer \eqref{ine517271}
  is equal to 
 $\ind(\Dirac_{\hat W}\otimes \hat \bU+\hat Q)$.

We now calculate this index. The index theory for these kinds of perturbations of Dirac operators
has been developed in  \cite{MR2191484}. In the language of this reference
 the operator $Q$ defines a
 taming $(M\otimes \bV)_{t}$  of the geometric manifold $M\otimes \bV$ and
a boundary taming $(W\otimes \bU)_{bt}$ of the geometric manifold $W\otimes \bU$. 
The following equality holds true by  definition of the right-hand side:
$\ind((\Dirac_{\hat W}\otimes \hat \bU+\hat Q)= \ind((W\otimes \bU)_{bt})$. The index theorem  
\cite[Thm. 4.18]{MR2191484} gives
  $$\ind((W\otimes \bU)_{bt})=\int_{W} p_{n+1}(\Td(\tilde \nabla^{TW})\wedge \ch(\nabla^{U}))-l\eta((M\otimes \bV)_{t}) \ .$$
 In $\R/\Z$ we have
 $[\eta((M\otimes \bU)_{t})]=\xi(D_{M}\otimes \bV)$.
 Hence by comparison with (\ref{eq4}) we get the equality in $\Q/\Z$
$$[\frac{1}{l}\ind((W\otimes \bU)_{bt})]=[\frac{1}{l}\ind((\Dirac_{W}\otimes \bU)_{APS})]\ .$$
In view of the construction of $\eta^{an}$, in particular of (\ref{eq133}),  we see that the map (\ref{eq62})
also represents $\eta^{an}(x)$.
This finishes the proof of Theorem \ref{indthm}. 
\hB
\section{An intrinsic formula}\label{sec444}

\subsection{Motivation}

In a typical situation for the theory of the present paper one is given a geometric representative $(M,f)$ for a torsion class $x=[M,f]\in \pi_{n}(MB)$ and wants to calculate the universal $\eta$-invariants $\eta^{top}(x)=\eta^{an}(x)\in Q_{n}(MB)$. The expressions for the universal $\eta$-invariant  that we have at our disposal at the moment  share the disadvantage that one has to find a lift $\hat x\in \pi_{n+1}(MB\Q/\Z)$ or a geometric zero bordism $(W,F)$ of $l$ copies of $(M,f)$ explicitly.
It is at this point where differential and spectral geometry helps. 
In the present section we develop a generalisation of Chern-Weyl theory
 which is designed to finally obtain formulas for the universal $\eta$-invariant which are intrinsic in the cycle $(M,f)$.

The main new object is the notion of a geometrisation of $(M,f,\tilde \nabla)$ which is defined in Definition \ref{def903}. It  involves differential $K$-theory which is reviewed 
in Subsection \ref{subsec41}.  
In Subsection \ref{sec18555} we  show the existence of geometrisations and study their functorial properties.  
In Subsection \ref{sec801} we introduce a special class of geometrisations which we call good.
In contrast to general geometrisations they have the property that they extend over zero bordisms.
The main result 
is the intrinsic formula for the universal $\eta$-invariant formulated in Theorem \ref{them2}.

\subsection{Review of differential $K$-theory}\label{subsec41}

The definition of a geometrisation utilises differential $K$-theory. We refer to 
\cite{MR1769477}, \cite{MR2192936}, \cite{MR2664467} for constructions and further information. In the following we review the basic structures  which by  \cite{MR2608479} uniquely characterise differential $K$-theory.
%
Differential $K$-theory is a five-tuple
$$(\hat K,I,R,a,\int)$$
of the following objects.
The first entry is a contravariant functor
$$\hat K:{\rm smooth\: manifolds}\longrightarrow {\rm \Z/2\Z-graded \:commutative\:rings}\ .$$
The remaining entries are natural transformations between functors. 
The domains and ranges of the first three are given by 
$$I:\hat K\to K\ ,\quad R:\hat K\to \Omega P_{cl}\ .$$
$$a:\Omega P/\im(d)[1]\to \hat K$$
 Here the evaluation of $\Omega P$ at $M$ is the graded vector space  $\Omega P(M):=\Omega(M)[b,b^{-1}]$  of two-periodic smooth real differential forms on $M$ which carries a differential $d$. By 
$ \Omega P_{cl}(M)\subseteq \Omega P(M)$ we denote its  subspace of closed forms.
The transformations $R$ and $I$ preserve the ring structures while $a$ is just additive.
These transformations are compatible in the sense that
for every manifold $M$ the following differential cohomology  diagram commutes
$$\xymatrix{&\Omega P^{*-1}(M)/\im(d)\ar[dr]^{a}\ar[rr]^{d}&&\Omega P_{cl}^{*}(M)\ar[dr]^{Rham}&\\HP\R^{*-1}(M)\ar[ur]\ar[dr]&&\hat K^{*}(M)\ar[ur]_{R}\ar[dr]^{I}&&HP\R^{*}(M)\\&K\R/\Z^{*-1}(M)\ar[ur]\ar[rr]^{Bockstein}&&K^{*}(M_{+})\ar[ur]^{\ch}&}\ .$$
Here we define the spectrum $HP\R$ representing periodic real cohomology similarly as $HP\Q$ in (\ref{eq8560}).
Furthermore, for $\alpha\in \Omega P^{*}(M)/\im(d)$ and $x\in \hat K^{*}(M)$ we have the identity
\begin{equation}\label{compaprod}a(\alpha)\cup x=a(\alpha\wedge R(x))\ .\end{equation}
The flat part of differential $K$-theory  is defined as the kernel of the curvature transformation $R$. It is canonically isomorphic to $\R/\Z$-$K$-theory (with a shift):\begin{equation}\label{eq700hhh}\hat K^{*}_{flat}(M):=\ker(R:\hat K^{*}(M)\to \Omega P^{*}_{cl}(M))\cong K\R/\Z^{*-1}(M ) \ .\end{equation}
The sequence
\begin{equation}\label{eq700}K^{*-1}(M )\stackrel{\ch}{\to} \Omega P^{*-1}(M)/\im(d)\stackrel{a}{\to} \hat K^{*}(M)\stackrel{I}{\to} K^{*}(M )\to 0\end{equation}
is exact. The integration is a natural (in $M$) transformation
$$\int:\hat K^{*}(S^{1}\times M)\to \hat K^{*-1}(M)$$
whose existence and compatibility with the other structures fixes the odd part of the  differential extension uniquely up to unique isomorphism as discussed in \cite{MR2608479}.
Since we do not need the integration in the present paper we will not write out the long list of these compatibilities explicitly.

Differential $K$-theory is not homotopy invariant. The deviation from homotopy invariance is quantified by the homotopy formula. If $\hat x\in \hat K^{*}([0,1]\times M)$, then it states that
\begin{equation}\label{homotopyformula}i_{1}^{*}\hat x-i_{0}^{*}\hat x=a(\int_{[0,1]\times M/M} R(\hat x))\ . \end{equation}

Let $\bV=(V,h^{V},\nabla^{V})$ be a geometric bundle on  manifold $M$, where $h^{V}$ is a hermitean metric which is  preserved by the connection $\nabla^{V}$. Then we have a natural class
\begin{equation}\label{eq129}[\bV]\in \hat K^{0}(M)\ .\end{equation}
This class is in fact tautological in the model \cite{MR2664467} in view of \cite[2.1.4]{MR2664467}. 
It satisfies
\begin{equation}\label{ccase}I([\bV])=[V]\in K^{0}(M)\ ,\quad R([\bV])=\ch(\nabla^{V})\in \Omega P_{cl}^{0}(M)\ ,\end{equation}
where $$\ch(\nabla^{V}):=\Tr(\exp(-\frac{bR^{\nabla}}{2\pi i  }))$$ is the normalized Chern character form.

\begin{rem}\label{remark222}{\rm
If we replace $HP\R$ and $\Omega P$ by $HP\C$ and $\Omega P\otimes \C$, then we get a complex version $\hat K_{\C}$ of differential $K$-theory with similar properties. A complex vector bundle with connection $\bV=(V,\nabla^{V})$  on a manifold $M$, where $\nabla^{V}$ is not necessarily hermitean, gives a class $[\bV]\in \hat K^{0}_{\C}(M)$ such that analogues of the equalities \eqref{ccase} still hold true. For the flat part we get the equivalence
$$\hat K^{0}_{\C,flat}(M)\cong K\C/\Z^{-1}(M)\ .$$
}\end{rem}

\subsection{Geometrisations}\label{sec18555}

Let $M$ be a compact manifold equipped with a map $f:M\to B$.  At the moment we do not require any connection of $f$ with the tangent bundle. Nevertheless we must imitate this situation.
 We can assume that
$f$ has a factorisation over $\tilde f:M\to BSpin^{c}(k)$  as in \eqref{fakdia} which classifies a $Spin^{c}(k)$-bundle $\tilde f^{*}Q_{k}\in Spin^{c}(\tilde  f^{*}\xi^{Spin^{c}}_{k})$ on $M$. The role of the tangent bundle is taken by the choice of a complementary $Spin^{c}$-bundle. In detail, we choose an  $l$-dimensional oriented euclidean  vector  bundle $\eta\to M$ for some $l\ge 0$
  together with an orientation preserving isomorphism  of euclidean vector bundles.
  \begin{equation}\label{er1}\eta\oplus \tilde f^{*}\xi^{Spin^{c}}_{k}\cong M\times \R^{l+k}\ .\end{equation}
Then we choose a $Spin^{c}$-structure $P\in Spin^{c}(\eta)$ together with an isomorphism
$$P\otimes \tilde f_{k}^{*}Q_{k}\cong Q(l+k)\ ,$$ where we use
the isomorphism (\ref{er1}) in order view the left-  and right-hand  sides in the same groupoid $Spin^{c}(M\times \R^{l+k})$ (see Subsection  \ref{sec1677} for details).

  We choose a connection $\tilde \nabla$ on $P$
  and get an induced Todd form $\Td(\tilde \nabla)\in \Omega P_{cl}^{0}(M)$ which represents the class $f^{*}\Td^{-1}\in HP\Q^{0}(M)$.

We now consider a continuous homomorpism 
$$\cG:K^{0}(B)\to \hat K^{0}(M)\ ,$$
where the domain has the  topology described in Remark \ref{proffinite1} 
and the target is discrete.
\textcolor{black}{\begin{ddd} A cohomological character for $\cG$ is a continuous  homomorphism 
$c_{\cG}:HP\Q^{0}(B)\to\Omega P^{0}_{cl}(M) $ such that the following diagram commutes: $$\xymatrix{K^{0}(B)\ar[r]^{\cG}\ar[d]^{\Td^{-1}\cup \ch(\dots)}&\hat K^{0}(M)\ar[d]^{\Td(\tilde \nabla)\wedge R(\dots)}\\HP\Q^{0}(B)\ar[r]^{c_{\cG}}&\Omega P^{0}_{cl}(M)}\ .$$\end{ddd} 
\begin{lem} Given $\cG$ there exists a cohomological character $c_{\cG}$.
If $B$ is compact, then it is unique.
\end{lem} 
\proof
Since
$\cG$ is continuous there exists a map $r:A\to B$ from a finite $CW$-complex $A$ such that
$\cG$ factors over the quotient
$K^{0}(B)/\ker(r_{K^{0}}^{*})$. For a space $X$ we write
  \begin{equation}\label{mirette20}HP\Q^{0}(X)_{0}:= \overline{\im(K^{0}(X)\otimes \Q\stackrel{\Td^{-1}\cup \ch(\dots)}{\longrightarrow}  HP\Q^{0}(X) )} \ .\end{equation} 
We get the diagram
$$\xymatrix{K^{0}(B) \ar[r]\ar@/^1cm/[rr]^{\cG}\ar[d]^{\Td^{-1}\cup \ch(\dots)}  &K^{0}(B) /\ker(r_{K^{0}}^{*})\ar[r]^{\bar \cG}\ar[d]^{\phi}  &\hat K^{0}(M)\ar[d]^{\Td(\tilde \nabla)\wedge R(\dots)}\\
HP\Q^{0}(B)_{0}\ar[d]\ar[r]^-{p_{0}}  & HP\Q^{0}(B)_{0}/\ker(r_{HP\Q^{0}}^{*})\cap HP\Q^{0}(B)_{0} \ar[d]  & \Omega P^{0}_{cl}(M)\\HP\Q^{0}(B)\ar@{.>}@/_2cm/[rru]_{c_{\cG}}\ar[r]^{p}&HP\Q^{0}(B)/\ker(r_{HP\Q^{0}}^{*})\ar@/_{1cm}/@{-->}[u]^{\pr}&}\ .$$   
The map $\phi$  induces an isomorphism after tensoring with $\Q$. In order to see this we consider the diagram
$$\xymatrix{0\ar[r]&K^{0}(B)\otimes \Q/\ker(r_{K^{0}\otimes \Q}^{*})\ar[d]^{\phi\otimes \Q}\ar[r]&K^{0}(A)\otimes \Q\ar[d]_{\cong}^{(r_{HP\Q^{0}}^{*}\Td^{-1}\wedge \ch(\dots))\otimes \Q}\\
0\ar[r]&HP\Q^{0}(B)_{0}/\ker(r_{HP\Q^{0}}^{*})\cap HP\Q^{0}(B)_{0}\ar[r]&HP\Q^{0}(A)}$$
with exact horizontal lines. We immediately see that $\phi\otimes \Q$ is injective. On the other hand, by definition \eqref{mirette20}
the image of $K^{0}(B)\otimes \Q\xrightarrow{(\Td^{-1}\cup \ch(\dots))\otimes \Q}HP\Q^{0}(B)$ is dense in $HP\Q^{0}(B)_{0}$. Note that the quotient
$HP\Q^{0}(B)/\ker(r_{HP\Q^{0}}^{*})$ carries the discrete topology. Therefore
 the composition 
$$K^{0}(B)\otimes \Q\xrightarrow{(\Td^{-1}\cup \ch(\dots))\otimes \Q}HP\Q^{0}(B)_{0}\to HP\Q^{0}(B)_{0}/\ker(r_{HP\Q^{0}}^{*})\cap HP\Q^{0}(B)_{0}$$
is surjective. 
This implies that $\phi\otimes \Q$ is surjective, too. }

\textcolor{black}{
After choosing some   linear projection $$\pr:HP\Q^{0}(B)/\ker(r_{HP\Q^{0}}^{*}) \to HP\Q^{0}(B)_{0}/\ker(r_{HP\Q^{0}}^{*})\cap HP\Q^{0}(B)_{0}$$ we can  define a cohomological character by
$$c_{\cG}(x):=\Td(\tilde \nabla)\wedge ((R\circ \bar \cG)\otimes \Q)\circ (\phi\otimes \Q)^{-1}(\pr(p(x)))\ , \quad x\in HP\Q^{0}(B)\ .$$ Since it factorizes over $p$ it is continuous.
 }

\textcolor{black}{
If $B$ is compact, then $HP\Q^{0}(B)_{0}=HP\Q^{0}(B)$ and this implies uniqueness of $c_{\cG}$.
\hB 
}

\textcolor{black}{
\begin{ex}[V\"olkl]{\rm If $B$ is not compact, then a cohomological character is not necessarily unique.
Consider $B:=K(\Z,4)$ and assume that $M$ is compact. By \cite[Thm. II]{MR0231369}  the reduced $K$-theory group $\tilde K^{0}(B)$ consists of phantom classes. Consequently, a  continuous map $\cG:\tilde K^{0}(B)\to K^{0}(M)$ must be trivial. Furthermore, we have $ HP\Q(B)_{0}=b^{0} \Q$.  Let $\cG:=0$. In this case any continuous  homomorphisms
$HP\Q^{0}(B) \to \Omega P^{0}_{cl}(M)$ vanishing on $b^{0} \Q$ can serve as a cohomological character $c_{\cG}$. Note that $HP\Q^{0}(B) \cong \Q[[q]]$ with $q:=b^{2}u$ for the canonical class $u\in H^{4}(K(\Z,4),\Q)$, so that there are many such homomorphisms.
}
\end{ex}}

  We say that the cohomological character $c_{\cG}$ preserves degree if it preserves the decompositions
$$HP\Q^{0}(B)\cong \prod_{k\in \Z} b^{-k} H^{2k}(B;\Q)\ ,\quad 
\Omega P^{0}_{cl}(M)\cong \prod_{k\in \Z} b^{-k}\Omega^{2k}_{cl}(M)\ .$$

\begin{ddd}\label{def903} A geometrisation of $(M,f,\tilde \nabla)$ is a continuous homomorphism 
$$\cG:K^{0}(B)\to \hat K^{0}(M)$$ such that the following diagram $$\xymatrix{ 
 &\hat K^{0}(M)\ar[d]^{I}  \\
K^{0}(B) \ar@{..>}[ur]^{\cG}\ar[r]^{f^{*}}&K^{0}(M) }$$ commutes, \textcolor{black}{and which admits a degree-preserving cohomological character.}
\end{ddd}
  
Note that the notion of a geometrisation only depends on the form $\Td(\tilde \nabla)$
representing the class $\Td^{-1}$. The additional rigidity of the construction related with the choice of $\eta$, $\tilde f$ and the isomorphism \eqref{er1}
is only used for the construction of the pull-back of geometrisations below.
  
\begin{ex}\label{example1}{\rm The notion of a geometrisation generalises the notion of a connection.
This is demonstrated in Lemma  \ref{lem77} for the case $B=BSpin$.
At this place we will discuss another example where we put $B:=BSpin^{c}\times B\Gamma$ for some \textcolor{black}{compact} Lie group $\Gamma$ and $B\to BSpin^{c}$ is the projection. We write maps to $BSpin^{c}\times B
\Gamma$ as pairs $(f,g)$. 
 
 Let us assume that we already have a geometrisation $\cG^{0}$ of
$(M,f,\tilde \nabla)$ \textcolor{black}{with a degree-preserving cohomological character $c_{\cG^{0}}$}. Its existence is guaranteed by Proposition \ref{prop904}. 
The map $g:M\to B\Gamma$ classifies a $\Gamma$-principal bundle
$R\to M$. We choose a connection  $\nabla^{R}$ on $R$.

\begin{lem}\label{cgl1} There exists a natural 
  geometrisation 
$\cG$ of $(M,(f,g),\tilde \nabla)$
associated to this data.
 \end{lem}
 \proof
The completion theorem \cite{MR0259946} gives an isomorphism $K^{0}(B\Gamma)\cong R(\Gamma)^{\hat{}}_{I_{\Gamma}}$ of topological groups, where  $I_{\Gamma}\subseteq R(\Gamma)$ is the dimension-ideal of the integral representation ring.
We consider a representation   $\sigma:\Gamma\to U(m_{\sigma})$ which represents an element $[\sigma]\in K^{0}(B\Gamma)$.
The
 associated complex vector bundle
 $V_{\sigma}:=R\times_{\Gamma,\sigma}\C^{m_{\sigma}}$ on $M$ then represents the element
 $[V_{\sigma}]=f^{*}[\sigma]\in K^{0}(M)$.
This bundle comes  with a hermitean metric
 $h^{V_{\sigma}}$ and a metric connection $\nabla^{V_{\sigma}}$ induced by   $\nabla^{R}$.
We therefore get a geometric bundle  $\bV_{\sigma}:=(V_{\sigma},\nabla^{V_{\sigma}},\nabla^{V_{\sigma}})$. It represents the class $[\bV_{\sigma}]\in \hat K^{0}(M)$ such that   $[V_{\sigma}]=I([\bV_{\sigma}])$ in $K^{0}(M)$, see (\ref{eq129}).
Let $\phi\in K^{0}(BSpin^{c})$. Then we get the element $\phi\times [\sigma]\in K^{0}(BSpin^{c}\times  B\Gamma)=K^{0}(B)$. We define
$$\cG(\phi\times [\sigma]):=\cG^{0}(\phi)\cup [\bV_{\sigma}]\ .$$
By linear extension this construction defines the map $\cG$   on a dense subgroup of
$K^{0}(B)$.

We now show that
the map $\cG$ extends by continuity to all of $K^{0}(B)$
and defines a geometrisation of $(M,(f,g),\tilde \nabla)$.
 Indeed, the map
$R(\Gamma)\to \hat K^{0}(M)$ induced by  $\sigma\mapsto [\bV_{\sigma}]$ is multiplicative and annihilates $I^{2n+1}_{\Gamma}$, where $n:=\dim(M)$.
Therefore, since $\cG^{0}$ is continuous,  the map $\cG$ is continuous as well.  
\textcolor{black}{We now use the fact that $HP\Q^{0}(B\Gamma)$ is topologically generated by the classes $\ch([\sigma])$ for $\sigma\in R(\Gamma)$.}
We let $c_{\Gamma}:HP\Q^{0}(B\Gamma)\to \Omega P_{cl}^{0}(M)$ be the unique continuous map such that
$\ch(\nabla^{V_{\sigma}})=c_{\Gamma}(\ch([\sigma]))$. Note that $c_{\Gamma}$ preserves degree. Since the cohomological character  $c_{\cG^{0}}$   preserves degree, \textcolor{black}{we can take 
$c_{\cG}:=c_{\cG^{0}}\wedge c_{\Gamma}$ of $\cG$ as a degree-preserving cohomological character for $\cG$.} 
\hB 
The geometrisation $\cG$ allows to recover the Chern character form of $\nabla^{V_{\sigma}}$ by
$$\ch(\nabla^{V_{\sigma}})=\Td(\tilde \nabla)^{-1}\wedge R(\cG(1\otimes [\sigma]))\ .$$
It also allows to partially recover transgressions as we will explain in the following. If $\nabla^{R\prime}$ is a second connection on $R$ and
$\cG^{\prime}$ is the associated geometrisation,
then
$$\cG^{\prime}(1\otimes [\sigma])-\cG_{M}^{\prime}(1\otimes [\sigma])=a(\Td(\tilde \nabla)\wedge \tilde \ch(\nabla^{V_{\sigma}\prime},\nabla^{V_{\sigma}}))\ .$$
Here $\tilde \ch(\nabla^{V_{\sigma}\prime},\nabla^{V_{\sigma}})\in \Omega P^{-1}(M)$ denotes the transgression form which satisfies
$$ d\tilde \ch(\nabla^{V_{\sigma}\prime},\nabla^{V_{\sigma}})=\ch(\nabla^{V_{\sigma}\prime})-\ch(\nabla^{V_{\sigma}})\ .$$[\textbf{End} of Example \ref{example1}] \hB}
\end{ex}

The following Proposition \ref{prop904} asserts
that geometrisations exist. Its proof uses the functoriality of  geometrisations in the space $B$.
Consider a  map $\phi$ over $BSpin^{c}$, i.e. a homotopy commutative diagram
$$\xymatrix{B^{\prime}\ar[rr]^{\phi}\ar[dr]&&B\ar[dl]\\&BSpin^{c}&}\ .$$
Given a geometrisation $\cG$ of $(M,\phi\circ f,\tilde \nabla)$ we get a geometrisation
\begin{equation}\label{eq601}\phi_{*}\cG:=\cG\circ \phi^{*}\end{equation} of $(M,f,\tilde \nabla)$.

Note that our standing assumption is that $M$ is compact. 
 \begin{prop}\label{prop904}
Given $(M,f,\tilde \nabla)$ there exists a geometrisation.
\end{prop}
\proof
Since $M$ is compact the map $f$ factors over a compact subspace of $B$.
In view of the functoriality of the geometrisation (\ref{eq601}) we can assume that $B$ is compact. Then $K^{0}(B)$ is a finitely generated abelian group.  We choose a   decomposition
$$K^{0}(B)\cong A_{tors}\oplus A_{free}$$
into a torsion and a free part. We write
$$A_{tors}:=\bigoplus_{y\in I} y\Z/\ord(y)\Z $$
for some set  of  generators $I\subset A_{tors}$.
For all $y\in I$, using the exactness at the right end of (\ref{eq700}), we choose  $\tilde y_{0}\in \hat K^{0}(M)$ such that $I(\tilde y_{0})=
f^{*}y$. Then
$\ord(y) \tilde y_{0}=a(\omega_{y})$ for some $\omega_{y}\in \Omega P^{-1}(M)/\im(d)$, again by  (\ref{eq700}).
We define
$$\tilde y:=\tilde y_{0}-a(\frac{1}{\ord(y)}\omega_{y})\ .$$
Then $\ord(y)\tilde y=0$ and we can define
$\cG_{|A_{tors}}:A_{tors}\to \hat K^{0}(M)$ such that $\cG(y)=\tilde y$ for all $y\in I$.
Since $\Td^{-1}\wedge \ch$ vanishes on $A_{tors}$ and $\cG_{|A_{tors}}$ maps to flat classes
it is clear that the cohomological character of this part of $\cG$  preserves degree.

We now come to the free part. We choose a basis $J\subset A_{free}$ and  classes
$\tilde z_{0}\in \hat K^{0}(M)$ such that $I(\tilde z_{0})=f^{*}z$ for all $z\in J$.
We further choose a basis 
$J^{\prime}\subset A_{free}\otimes \Q$ such that
$\{\Td^{-1}\wedge \ch(z^{\prime})\}_{z^{\prime}\in J^{\prime}}$ is a homogeneous basis with respect to the decomposition
 $$HP\Q^{0}(B)\cong \bigoplus_{m\in \Z} b^{m}H^{2m}(B;\Q)\ .$$ 
 We define the even integers  $n_{z^{\prime}}:=\deg(\Td^{-1}\wedge \ch(z^{\prime}))$ for all $z^{\prime}\in J^{\prime}$.
Then there exists an invertible rational $(J,J^{\prime})$-indexed matrix $A$ such that
$z=\sum_{z^{\prime}\in J^{\prime}}A_{z z^{\prime}} z^{\prime}$ for all $z\in J$. 
We now can choose forms $\alpha_{z^{\prime}}\in \Omega P^{-1}(M)/\im(d)$ for all $z^{\prime}\in J^{\prime}$ such that
$$\sum_{z\in J}A^{-1}_{z^{\prime}z} \;\: \Td(\tilde \nabla)\wedge R(\tilde z_{0})-  d\alpha_{z^{\prime}}\in b^{\frac{n_{z^{\prime}}}{2}}\Omega_{cl}^{n_{z^{\prime}}}(M)\subseteq \Omega P_{cl}^{0}(M)$$
for all $z^{\prime}\in J^{\prime}$.
We define
$$\cG_{|A_{free}}:A_{free}\to \hat K^{0}(M)$$ by linear extension
such that
$$\cG(z)=\tilde z_{0}-a(\Td(\tilde \nabla)^{-1}\wedge \sum_{z^{\prime}\in J^{\prime}} A_{zz^{\prime}}\alpha_{z^{\prime}})\ .$$
By construction its \textcolor{black}{(uniquely determined)} cohomological character preserves degree.

\hB

Geometrisations can be pulled back along stable $Spin^{c}$-maps over $X$.  In detail the construction goes as follows.
Let $(M^{\prime},f^{\prime} )$ be a compact manifold with a map  $f^{\prime}:M^{\prime}\to B$ . We consider a smooth map $h:M^{\prime}\to M$  such that $f\circ h$ is homotopic to $f^{\prime}$. This implies that we can choose a stable isomorphism of complementary bundles
\begin{equation} \label{iso345q}\eta^{\prime}\oplus (M^{\prime}\times \R^{s})\cong h^{*}\eta\oplus (M^{\prime}\times \R^{t})\ .\end{equation} We refine $h$ to a $Spin^{c}$-map  
  by choosing an isomorphism
\begin{equation} \label{iso345}  P^{\prime}\otimes Q(s)\cong h^{*}P\otimes Q(t)\ .\end{equation}

Assume now that we have connections $\tilde \nabla$ on $P$ and $\tilde \nabla^{\prime}$ on $P^{\prime}$. They induce connections on the stabilisations   $P\otimes Q(t)$ and  $ P^{\prime}\otimes Q(s)$.   We thus can define the transgression
 $$\tilde  \Td(h^{*}\tilde \nabla,\tilde \nabla^{\prime})\in \Omega P^{-1}(M^{\prime})/\im(d)$$
 where we use the isomorphism (\ref{iso345}) in order to compare the stabilisation of  $h^{*}\tilde \nabla$ with that of $\tilde \nabla^{\prime}$ on the same bundle. 
 The transgression satisfies
 $$d\tilde  \Td(h^{*}\tilde \nabla,\tilde \nabla^{\prime}) =h^{*}\Td(\tilde \nabla)-\Td(\tilde \nabla^{\prime})\ .$$
  
 Let  $\cG$ be a geometrisation of 
$(M,f,\tilde \nabla)$.
\begin{lem}\label{lem906} If $h:M^{\prime}\to M$ is a  $Spin^{c}$-map between compact manifolds, then
there exists a  construction of a pull-back $\cG^{\prime}:=h^{*}\cG$  of $(M^{\prime},f^{\prime}, \tilde \nabla^{\prime})$. This pull-back only depends on the joint stable homotopy class of isomorphisms \eqref{iso345q} and \eqref{iso345} and is functorial under compositions.
\end{lem}
\proof
By our assumptions the equivalence class
\begin{equation}\label{eqq4}\beta:=\tilde \Td(h^{*}\tilde \nabla,\tilde \nabla^{\prime})\wedge\Td(\tilde \nabla^{\prime})^{-1}\in \Omega P^{-1}(M^{\prime})/\im(d)\end{equation}
of forms is well-defined.
It  satisfies
\begin{equation}\label{betafix1}d\beta=h^{*}\Td(\tilde \nabla)\wedge\Td(\tilde \nabla^{\prime})^{-1}-1\ .\end{equation}
We define the pull-back $\cG^{\prime}:=h^{*}\cG$ by 
\begin{equation}\label{eq1300}\cG^{\prime}(y):=h^{*}\cG(y)+a(\beta\wedge h^{*}R(\cG(y)))\ ,\quad  y\in K^{0}(B )\ ,\end{equation}
where $a$ and $R$ belong to the structure maps of differential $K$-theory.
We have by construction
$$\Td(\tilde \nabla^{\prime})\wedge R(\cG^{\prime}(y))= h^{*}(\Td(\tilde \nabla)\wedge R(\cG(y)))$$ \textcolor{black}{
and hence can take
$c_{\cG^{\prime}}:=h^{*}c_{\cG}$ as a cohomological character for $\cG^{\prime}$.} 
Since the cohomological character $c_{\cG}$ preserves degree,  so does  the cohomological character of $\cG^{\prime}$.

We show that the pull-back is functorial.
We consider a second triple  $(M^{\prime\prime},f^{\prime\prime}, \tilde \nabla^{\prime\prime})$ with a $Spin^{c}$-map
 $h^{\prime}:M^{\prime\prime}\to M^{\prime}$   and the associated transgression form $\beta^{\prime}$.
Then we have for the iterated pull-back
\begin{eqnarray*}\cG^{\prime\prime}(y)&=& h^{\prime *}(h^{*}(\cG(y)))+h^{\prime *} (a(\beta\wedge h^{*}R(\cG(y))))+a(\beta^{\prime}\wedge h^{\prime *} R(\cG^{\prime}(y)))\\
&=&(h\circ h^{\prime *})^{*}(\cG(y))\\&&\hspace{-2cm}+a(h^{\prime *} \beta\wedge h^{\prime *}(h^{*}(R(\cG(y)))+\beta^{\prime}\wedge h^{\prime *}( h^{*}R(\cG(y)))+\beta^{\prime}\wedge h^{\prime *}d\beta\wedge h^{\prime *}(h^{*}(R(\cG(y)))))\\\end{eqnarray*}
Let $\tilde \beta$ be the transgression form for the composition $h\circ h^{\prime}$ of $Spin^{c}$-maps over $X$. Then we must show that
$$\tilde \beta-(h^{\prime *}\beta +\beta^{\prime}+h^{\prime *}\beta^{\prime}\wedge d\beta)\in \im(d)$$
This follows from $$d(h^{\prime *}\beta +\beta^{\prime}+h^{\prime *}\beta^{\prime}\wedge d\beta)=h^{\prime* }h^{*}\Td(\tilde \nabla)^{-1}\wedge\Td(\tilde \nabla^{\prime\prime})-1=d\tilde \beta$$
and the fact that all these forms are defined by transgressions and the contractibility of the space of connections.
The assertion about homotopy invariance easily follows from the homotopy formula \eqref{homotopyformula} for differential $K$-theory.
\hB

Note that the form $\beta$  is determined up the closed forms by \eqref{betafix1}.
The refinement of the map $h$ to a $Spin^{c}$-map is necessary in order to rigidify the choice of $\beta$
up to exact forms by \eqref{eqq4} by constructing it via transgression.

\begin{rem}{\rm
The identity of $M$ refines to a $Spin^{c}$-map   in a natural way by choosing the identity in (\ref{iso345}).
The pull-back of geometrisations for the identity of $M$ can be used to transfer a geometrisation defined for one choice of the connection $\tilde \nabla$ to a second choice.  This allows to define a notion of geometrisation which is independent of the choice of the connection. This could play a role of one wants to classify geometrisations. We will not pursue that goal in the present paper.}
\end{rem}

\subsection{Good geometrisations}\label{sec801}

Assume that $(W,F)$ is a zero-bordism of the $n$-dimensional cycle $(M,f)$. 
We fix tangential  $Spin^{c}$-structures $P(TW)\in Spin^{c}(TW)$ and $P(TM)\in Spin^{c}(TM)$ associated to  the normal $Spin^{c}$-structure on $W$ and $M$, see Definition \ref{tang}.  As explained in Subsection \ref{sec1677} there is a natural   isomorphism  of $Spin^{c}$-structures
 \begin{equation}\label{eq1573}P(TM)\otimes Q(1)\cong P(TW)_{|M}\end{equation}
 which turns the inclusion $$i :M\to W$$ into a  $Spin^{c}$-map.
 
We choose a $Spin^{c}$-extension of the Levi-Civita connection
$\tilde \nabla^{TW}$ on $W$ with product structure and a  $Spin^{c}$-extension of the Levi-Civita connection $\tilde \nabla^{TM}$ on $M$ such that the isomorphism (\ref{eq1573}) preserves the connections. In this situation  the form (\ref{eqq4}) is trivial.
Assume now that we have a geometrisation of $(W,F,\tilde \nabla^{TW})$.
Then we can define the restriction
$\cG_{\partial W}:=(\cG_{W})_{|\partial W}$ as in Lemma \ref{lem906}. It is given by
\begin{equation}\label{eqq6}\cG_{\partial W}(\phi)=\cG_{W}(\phi)_{|\partial W}, \quad \phi\in K^{0}(B)\ .\end{equation}
 In general we do not expect that a given geometrisation $\cG_{M}$ of $(M,f,\tilde \nabla^{TM})$
 can be obtained by restricting  a geometrisation $\cG_{W}$ of $(W,F,\tilde \nabla^{TW})$. In this respect geometrisations are more rigid than connections.

\begin{ex}\label{nonexr}{\rm

Here is a very simple example of a geometrisation which does not extend. We consider the case $B=*$ and let $(S^{3},f)$  be a cycle for $\pi_{3}(S)$.   We choose a normal framing of $S^{3}$ such that  it extends over $D^{4}$ so that the  framed bordism class $[S^{3},f]$ is trivial. Furthermore we equip $S^{3}$ with its standard Riemannian metric.
 
 We let $\cG_{0}$ be the good geometrisation of $(S^{3},f,\tilde \nabla^{TS^{3}})$
defined  in Subsection \ref{subsec51}. 
We have
$K^{0}(B)\cong  K^{0}(*)\cong \Z$ so that a geometrisation is fixed by the image of $1$ in 
$\hat K^{0}(S^{3})$.  
Let $\omega\in \Omega^{3}(S^{3})$ be some form. Then we can define a new geometrisation
$\cG_{\omega}$ of $(S^{3},f,\tilde \nabla^{TS^{3}})$ by $$\cG_{\omega}(1):=\cG_{0}(1)+a(\omega)\ .$$ It is easy to check,  using the fact that $\cG_{0}$ does extend by Lemma \ref{lem907},  that $\cG_{\omega}$  extends to  $D^{4}$ if and only if  $\int_{S^{3}}\omega\in \Z$.
[\textbf{End} of example \ref{nonexr}] \hB}\end{ex}

 In order to deal with the problem of non-extendability of geometrisations appropriately we introduce the notion of a $k$-good geometrisation. If $\cG_{M}$ is $k$-good with $k\ge \dim(M)+1$, then it will extend to zero bordisms.

For $k\in \nat$ we define the notion of $k$-good geometrisations constructively.  
We consider a map    $f_{u}:M_{u}\to B$ from a smooth compact manifold $M_{u}$, a lift $\tilde f_{u}$ together with a choice of a complementary bundle $\eta_{u}$, a connection $\tilde \nabla_{u}$ on $P_{u}\in Spin^{c}(\eta_{u})$, and
a geometrisation
$\cG_{u}$ of
$(M_{u},f_{u},\tilde \nabla_{u})$. If the map  $f:M\to B$  has a   factorisation up to homotopy through a smooth map $h:M\to M_{u}$, then we can refine $h$ to    
a stable $Spin^{c}$-map since, after stabilisation, there exists an isomorphism between $h^{*}\eta_{u}
$ and $TM$.
\begin{ddd}\label{goodgeom}
The geometrisation $\cG_{M}$ is called $k$-good if $\cG_{M}=h^{*}\cG_{u}$
for some choices as above such that $f_{u}$ is $k$-connected.
We say that $\cG_{M}$ is good, if it is $\dim(M)+1$-good. \end{ddd}

\begin{rem}\label{kdsdkfsdf}{\rm
We consider a sequence $(M_{u,i},f_{u,i},\tilde \nabla^{u,i})$, $i\in \nat$, of data as
above together with $Spin^{c}$-maps $h_{u,i}:M_{u,i}\to M_{u,i+1}$ for all $i\in \nat$ such that
$f_{u,i+1}\circ h_{u,i}\sim f_{u,i}$ and $f_{u,i}$ is $i$-connected. A family of geometrisations $\cG_{u,i}$ of $(M_{u,i},f_{u,i},\tilde \nabla^{u,i})$ for all $i\in \nat$  such that $h_{u,i}^{*}\cG_{u,i+1}=\cG_{u,i}$ will be called a universal geometrisation. Universal geometrisations will be constructed and classified in a forthcoming paper by M.V\"olkl.

Let us fix a universal geometrisation. A geometrisation of $(M,f,\tilde \nabla^{TM})$ will be called very good (relative to the chosen universal geometrisation), if it is isomorphic to $
h^{*}\cG_{u,i}$ for some $i\in \nat$ and $Spin^{c}$-map $h:M\to M_{u,i}$ such that $f_{u,i}\circ h\sim f$. Note that such a geometrisation is $\ell$-good for every $\ell\in \nat$.}
\end{rem}

 \begin{lem}\label{lem1101}
If $B$ has the  homotopy type of a $CW$-complex with finite skeleta, then for every triple $(M,f,\tilde \nabla^{TM})$ and every $k\in \nat$ there exists a $k$-good geometrisation.
\end{lem}
\proof
By the assumption, for every $k\in \nat$ we can find a compact manifold $M_{u}$  and a  map
$f_{u}:M_{u}\to B$  such that
$f_{u}:M_{u}\to B$ is $k$-connected. 
We choose  complementary data $\tilde f_{u}$, $\eta_{u}$, $P_{u}$ and  $\tilde \nabla^{u}$  as above. 
By Proposition \ref{prop904} there exists a geometrisation $\cG_{u}$ of the triple
$(M_{u},f_{u},\tilde \nabla^{u})$. Given $(M,f)$ with  $\dim(M)\le k-1$ there exists a factorisation up to homotopy
\begin{equation}\label{eq734}\xymatrix{&M_{u}\ar[d]^{f_{u}}\\M\ar[r]^{f}\ar@{..>}[ur]^h&B}\end{equation}
and a refinement of $h$ to a $Spin^{c}$-map.
Then $\cG_{M}:=h^{*}\cG_{u}$ is a $k$-good geometrisation. \hB

\begin{lem}\label{lem907}
Let $\cG_{M}$ be  a good geometrisation of $(M,f,\tilde \nabla^{TM})$. If $(W,F,G)$ is a zero bordism of $(M,f)$  with connection $\tilde \nabla^{TW}$, then there exists a geometrisation
$\cG_{W}$ of $(W,F,\tilde \nabla^{TW})$ which restricts to $\cG_{M}$.
\end{lem}
\proof
Since $f_{u}:M_{u}\to B$ is an $n+1$-equivalence and $\dim(W)=n+1$
we can   extend the factorisation (\ref{eq734})  to a factorisation
\begin{equation}\label{eq734vvvv}\xymatrix{&M_{u}\ar[d]^{f_{u}}\\W\ar[r]^{F}\ar@{..>}[ur]^H&B}\end{equation}
There exists a refinement of $H$ to a stable $Spin^{c}$-map such that $H\circ i=h$ in the sense $Spin^{c}$-maps.
Then we can define the pull-back  $\cG_{W}:=H^{*}\cG_{u}$ and get $(\cG_{W})_{|M}=\cG_{M}$
 by the   functoriality  of the pull-back. \hB 

\begin{rem}\label{chadcldc}{\rm Note that Lemma \ref{lem1101}  does not imply the existence of very good geometrisations. Furthermore observe the following potentially bad behaviour with respect to disjoint unions. Let $(M_{i},f_{i},\tilde\nabla^{TM_{i}})$, $i=0,1$ be two geometric cycles with  $k$-good geometrisations $\cG_{i}$. Then we can form the disjoint union
$(M,f,\tilde \nabla^{TM}):= (M_{0},f_{0},\tilde\nabla^{TM_{0}})\sqcup (M_{1},f_{1},\tilde\nabla^{TM_{1}})$ which carries a geometrisation $\cG$ naturally induced by $\cG_{i}$. It is not clear that
 this geometrisation is   $k$-good. For this it would good to know that the geomtrizations $\cG_{i}$ are pulled back from the {\em same} geometrisation $\cG_{u}$ by different maps. This is  the motivation behind the notion of a very good geometrisation.
 
}\end{rem}

\subsection{An intrinsic formula for $\eta^{an}$}\label{intr1z}

 The main goal of the present Subsection is to give an intrinsic formula for $\eta^{an}(x)$ which only involves structures on the cycle  $(M,f)$ for $x\in \pi_{n}(MB)_{tors}$. 

The geometric and analytic terms in the formula (\ref{eq870}) for $\eta^{an}(x)$
separately have values in $\R/\Z$; only their sum belongs to $\Q/\Z$. In order to deal with these terms separately it is useful to use a real version $Q_{n}^{\R}(E)$ of the group $Q_{n}(E)$.
We start with introducing this group. We further show that there is no  loss of information when going over the  real version. 
We let (compare with \eqref{sec31u})
\begin{equation}\label{eq1002}U^{\R}\subseteq \Hom^{cont}(K^{0}(E),\pi_{n+1}(K\R/\Z))\end{equation} be the subgroup
given by evaluations against elements in
$\pi_{n+1}(E\R)$ and define
\begin{equation}\label{eq1003}Q^{\R}_{n}(E):=\frac{\Hom^{cont}(K^{0}(E),\pi_{n+1}(K\R/\Z))}{U^{\R}}\ .\end{equation}
The inclusion $\pi_{n+1}(K\Q/\Z)\to \pi_{n+1}(K\R/\Z)$ induces a map
$$i_{\R}:Q_{n}(E)\to Q_{n}^{\R}(E)\ .$$ 
\begin{lem}
The map $i_{\R}:Q_{n}(E)\to Q_{n}^{\R}(E)$ is injective.
\end{lem}
\proof
Let $\kappa\in Q_{n}(E)$ be represented by
$\hat \kappa\in \Hom^{cont}(K^{0}(E),\pi_{n+1}(K\Q/\Z))$.
Since $\kappa$ is continuous it factors over a finitely generated quotient of 
$K^{0}(E)$. Hence there exists $N\in \nat$ such that $N\hat\kappa$ vanishes.

Assume now that $i_{\R}(\kappa)=0$. Then there exists
$w\in \pi_{n+1}(E\R)$ such that $\hat \kappa(\phi)=[\langle w,\phi\rangle]\in \pi_{n+1}(K\R/\Z)$ for all $\phi\in K^{0}(E)$. 
Since $\pi_{n+1}(E\R )\cong \pi_{n+1}(E)\otimes \R$  (see (\ref{torseq})) there exists a finite subset
$I\subset \pi_{n+1}(E)$ and a map $\lambda:I\to \R$ such that
$w=\sum_{v\in I} \lambda(v)v$. We have
$\hat \kappa(\phi)=\sum_{v\in I}[\lambda(v)\langle \phi,v\rangle]$, where here $\langle \phi,v\rangle\in \pi_{n+1}(K)$.
For $v\in I$ we define
$\hat v\in \Hom^{cont}(K^{0}(E),\pi_{n+1}(K))$ by
$\hat v(\phi):=\langle \phi,v\rangle$. The set $\{\hat v|v\in I\}$ generates a free
abelian subgroup $A\subseteq\Hom^{cont}(K^{0}(E),\pi_{n+1}(K))$.
We can choose a minimal subset $J\subseteq I$ which generates a subgroup of $A$ of full rank.
Then there exists a suitable 
map $\mu:J\to \R$   such that
$\hat \kappa(\phi)=\sum_{v\in J}[\mu(v)\hat v(\phi)]$  for all $\phi\in K^{0}(E)$.

The image of $K^{0}(E)\to \Hom(A,\pi_{n+1}(K))$ has full rank.
Hence for every $v\in J$ there exists $\phi_{v}\in K^{0}(E)$
such that $\hat v(\phi_{v})\not=0$ and $\hat v^{\prime}(\phi_{v})=0$ for all $J\ni v^{\prime}\not= v$.
It follows that $\hat \kappa(\phi_{v})=[\mu(v) \hat v(\phi_{v})]$.
Since
$0=N\hat \kappa(\phi_{v})=[N\mu(v) \hat v(\phi_{v})]$ it follows that
$\mu(v)\in \Q$. We set $w_{\Q}:=\sum_{v\in J} \mu(v)v\in \pi_{n+1}(E)\otimes \Q\cong \pi_{n+1}(E\Q)$.  Then we have
$\hat \kappa(\phi)=[\langle \phi, w_{\Q}\rangle]$ for all $\phi \in K^{0}(E)$.
This shows that $\hat \kappa\in U$ and $\kappa=0$.
\hB

Let $x\in \pi_{n}(MB)_{tors}$ be an $l$-torsion element and $(M,f)$ be a cycle for $x$.
We choose a  $Spin^{c}$-extension $\tilde \nabla^{TM}$ of the Levi-Civita
connection on $M$. We further assume that we have a good geometrisation $\cG_{M}$ of $(M,f,\tilde \nabla^{TM})$ (see Definition \ref{goodgeom}) \textcolor{black}{with a choice of a degree-preserving cohomological character $c_{\cG}$}. If $B$ has finite skeleta, then its existence is guaranteed by Lemma \ref{lem1101}.

For every $\phi$ in $K^{0}(B)$ we choose a $\Z/2\Z$-graded vector bundle $V_{\phi}\to M$  such that   $[V_{\phi}]=f^{*}\phi$ in $K^{0}(M)$. We furthermore choose a
hermitean metric $h^{V_{\phi}}$ and a metric connection $\nabla^{V_{\phi}}$ so that we get the geometric bundle
$\bV_{\phi}=(V_{\phi},h^{V_{\phi}},\nabla^{V_{\phi}})$. It represents a differential $K$-theory class $[\bV_{\phi}]\in \hat K^{0}(M)$ such that $I([\bV_{\phi}])=[V_{\phi}]=f^{*}\phi=I(\cG_{M}(\phi))$.
By the exactness of (\ref{eq700}) 
 we get a uniquely determined element  $$\gamma_{\phi}\in \Omega P^{-1}(M)/\im(\ch)$$  such that
\begin{equation}\label{eq1450}\cG_{M}(\phi)-[\bV_{\phi}]=a(\gamma_{\phi})\ .\end{equation}
\begin{ddd}\label{cgkl1}
We will refer to $\gamma_{\phi}$ as the correction form associated to $\phi$.
\end{ddd}
\begin{theorem}\label{them2}
The element $$i_{\R}(\eta^{an}([M,f]))\in Q^{\R}_{n}(MB)$$ is represented by the homomorphism
\begin{equation}\label{eq1000} K^{0}(MB)\stackrel{\Thom^{K}}{\cong }K^{0}(B)\ni \phi \mapsto [-\int_{M} \Td(\tilde \nabla^{TM})\wedge \gamma_{\phi}] -\xi(\Dirac_{M}\otimes \bV_{\phi}) \in \R/\Z\ .\end{equation}
\end{theorem}
\proof
The integral in formula \eqref{eq1000} belongs to  the group $\R[b,b^{-1}]^{-n-1}$ which will be identified with $\R$ using the generator
$b^{-\frac{n+1}{2}}$.
First note that, despite of the fact that $\gamma_{\phi}$ is only defined up the image of $\ch:K^{-1}(M)\to HP\Q^{-1}(M)$,  the class
$$[\int_{M} \Td(\tilde \nabla^{TM})\wedge \gamma_{\phi}]\in \R/\Z$$ is well-defined. Indeed, we have  
$\langle [M],\Td(TM)\cup \ch(\psi)\rangle\in \Z$ for all $\psi\in K^{-1}(M)$ by the odd version of  Atiyah-Singer index theorem.
We use  (\ref{eq870}) in order to express the right-hand side of
(\ref{eq133}) as
$$[\frac{1}{l}\int_{W}\Td(\tilde \nabla^{TW})\wedge \ch(\nabla^{U_{\phi}})]-\xi(\Dirac_{M}\otimes\bV_{\phi})\ .$$
The whole idea  is now to turn the integral over $W$ into an integral over $M$.
To this end we assume by Lemma \ref{lem907} that the good geometrisation $\cG_{M}$ has an extension $\cG_{W}$ to $W$. 
The $K$-theory class $f^{*}\phi$ extends across $W$ as $F^{*}\phi$.
We can thus assume, after adding a bundle of the form $\bW\oplus \bW^{op}$,  that the bundle $\bV_{\phi}$ has an extension $\bU_{\phi}$ as a geometric bundle to $W$. Note that this sort of stabilisation does not
effect the correction form $\gamma_{\phi}$ and the reduced $\eta$-invariant
$\xi(\Dirac_{M}\otimes \bV_{\phi})$. From now on we assume that $\bV_{\phi}$ extends.
We let $\gamma^{W}_{\phi}\in \Omega P^{-1}(W)/\im(\ch)$ be the correction form defined by 
$$\cG_{W}(\phi)-[\bU_{\phi}]=a(\gamma_{\phi}^{W})\ .$$
By (\ref{eqq6}) we conclude that $(\gamma^{W}_{\phi})_{|\partial W}$ coincides up to the image of $\ch$ with $\gamma_{\phi}$ on all copies of $M$. We now use Stokes' theorem    in order to rewrite
\begin{eqnarray*}
[\frac{1}{l}\int_{W}\Td(\tilde \nabla^{TW})\wedge \ch(\nabla^{U_{\phi}})] &=& [
\frac{1}{l}\int_{W} \Td(\tilde \nabla^{TW})\wedge R(\cG_{W}(\phi)) - \frac{1}{l}\int_{W} \Td(\tilde \nabla^{TW})\wedge d\gamma^{W}_{\phi}]          \\
&=&[\frac{1}{l}\int_{W} \Td(\tilde \nabla^{TW})\wedge R(\cG_{W}(\phi))]-[\int_{M} \Td(\tilde \nabla^{TM})\wedge \gamma_{\phi}] \end{eqnarray*}
  We want to show that  homomorphism
 $$\kappa:K^{0}(B)\ni \phi\mapsto [\frac{1}{l}\int_{W} \Td(\tilde \nabla^{TW})\wedge R(\cG_{W}(\phi))]\in \R/\Z$$ belongs to $U^{\R}$. The integrand of the integral over $W$ \textcolor{black}{can be expressed in terms of the cohomological character $c_{\cG_{W}}$} of $\cG_{W}$. Therefore $\kappa$ has 
    a factorisation as
 $$K^{0}(B)\xrightarrow{ \Td^{-1}\cup\ch} HP\Q^{0}(B)\xrightarrow{c_{\cG_{W}}} \Omega P^{0}_{cl}(W)\xrightarrow{\frac{1}{l}\int_{W}}\R\xrightarrow{[\dots]_{\R/\Z}} \R/\Z\ .$$
Since the cohomological character $c_{\cG_{W}}$ preserves degree and $\frac{1}{l}\int_{W}$ factorises over the degree $n+1$-part the homomorphism $\kappa$ actually factorizes over 
$$K^{0}(B)\xrightarrow{ p_{n+1}(\Td^{-1}\cup\ch)} H\Q^{n+1}(B)\xrightarrow{(c_{\cG_{W}})_{|H\Q^{n+1}(B)}} \Omega_{cl}^{n+1}(W)\xrightarrow{\frac{1}{l}\int_{W}}\R\xrightarrow{[\dots]_{\R/\Z}} \R/\Z\ .$$
Since $c_{\cG_{W}}$ is continuous,  the composition
$$H\Q^{n+1}(B)\xrightarrow{(c_{\cG_{W}})_{|H\Q^{n+1}(B)}} \Omega^{n+1}_{cl}(W)\xrightarrow{\frac{1}{l}\int_{W}}\R$$
is continuous and therefore given by the pairing against an element of $H\R_{n+1}(B)$. 
We now use  the $\R$-version of   \eqref{localthomformula} in order to conclude that $\kappa\in U^{\R}$.
 Therefore
 $i_{\R}(\eta^{an}([M,f]))$ is also represented by 
 the map \eqref{eq1000}. \hB

 \begin{rem}{\rm
Let us mention the following aspect of the intrinsic formula (\ref{eq1000}) which is not yet completely
understood at the moment. For the intrinsic formula to make sense we do not need the 
zero bordism $(W,F)$ of $l$-copies of $(M,f)$. Therefore formula \eqref{eq1000} provides an element
$\eta^{intrinsic}_{\cG}(M,f,\tilde \nabla^{TM})\in Q^{\R}_{n}(MB)$. 
One can check that $\eta^{intrinsic}_{\cG}(M,f,\tilde \nabla^{TM})=0$
if the data $(M,f,\tilde \nabla^{TM},\cG_{M})$ extends to a zero bordism.
Note that in general two good geometrisations of $(M,f,\tilde \nabla^{TM})$ cannot be connected over the cylinder but  this is possible if both are very good with respect to the same universal geometrisation. More general, if we work with very good geometrisation associated to a fixed universal geometrisation (see Remark \ref{kdsdkfsdf}), then the problem with disjoint unions
mentioned in Remark \ref{chadcldc} disappears, too.
Therefore, if we fix a universal geometrization, then  we would get a homomorphism
$$\eta_{\cG}^{intrinsic}:\pi_{n}(MB)\to Q^{\R}_{n}(MB)$$
which restricts to $i_{\R}\circ \eta^{top}=i_{\R}\circ \eta^{an}$ on
$ \pi_{n}(MB)_{tors}$. In general we do not know the topological contents of this extension of the universal $\eta$-invariant.
It might be related to the effect observed  in Remark \ref{hiurfhrifhfx}.
 }\end{rem}

\section{Examples}\label{sect5}

\subsection{Adams' $e$-invariant}\label{subsec51}

We consider the example $B=*$. The associated Thom spectrum is the sphere spectrum $S\cong MB$.
By Serre's theorem \cite{MR0059548} the homotopy groups $\pi_{n}(S)$ are finite for $n\ge 1$ (we refer to  \cite{MR860042} for more details about their structure). 
Therefore the universal $\eta$-invariant is defined on all of $\pi_{n}(S)$.



We have an identification
$K^{0}(S)\cong \Z$.
Furthermore, since $\pi_{n+1}(S\Q)\cong 0$ for $n\ge 0$   the group $U$ defined in \eqref{sec31u} is trivial. From now one let $n\in \nat$ be odd. After identifying 
$\pi_{n+1}(K\Q/\Z)\cong \Q/\Z$   we obtain  the identification
$$Q_{n}(S)\cong \Hom(\Z,\Q/\Z)\cong \Q/\Z$$ given by evaluation against $1\in \Z$.
For every odd $n\in \nat$ the universal $\eta$-invariant is thus interpreted as a homomorphism
\begin{equation}\label{eq1200} \eta:\pi_{n}(S)\to \Q/\Z\ .\end{equation}
The universal $\eta$-invariant essentially coincides with Adams' $e$-invariant  \cite{MR0198470}
$$e^{Adams}:\pi_{n}(S)\to \Q/\Z$$ which was introduced in order to detect  
the image of the $J$-homomorphism
\begin{equation}\label{eq30001}J:KO_{n+1}\to \pi^{S}_{n}\ .\end{equation} 
We will see below that we have the  relation
\begin{equation}\label{gleichheit}\eta=e^{Adams}_{\C}:=\left\{\begin{array}{cc}e^{Adams}&\frac{n-1}{2} \:even\\
2e^{Adams}&\frac{n-1}{2}\: odd\end{array}\right. \ . \end{equation}

We consider the Adams filtration
  of  $\pi_{*}(S)$ associated to the $K$-theory based Adams spectral sequence (see e.g. \cite{MR1324104} or Section \ref{secneu8}). It follows from Proposition \ref{prop1}, $4.$, that  the universal $\eta$-invariant  (\ref{eq1200}) induces an injection
$$\eta^{top}:\Gr^{1}\pi_{n}(S) \hookrightarrow \Q/\Z$$ (recall that $n$ is assumed to be odd).
The relation of the Adams $e$-invariant to spectral geometry  has first been observed in \cite{MR0397799}. The spectral geometric calculation of the Adams $e$-invariant interprets $\pi_{n}(S)$ as a framed bordism group.
 It has the favourable
property that it provides an intrinsic formula for $e^{Adams}$, a fact which  has been successfully exploited e.g. in \cite{MR762355}, \cite{MR687857}. 

The goal of the following discussion is to derive, using Theorems \ref{them2} and \ref{indthm},   the   intrinsic formula \eqref{hfjhgerjfgjerferffref} for $\eta^{top}$ and compare it with the formula \cite[Thm 4.14]{MR0397798} for
$e^{Adams}_{\C}$.
This  is finally our argument for \eqref{gleichheit}.

We let $(M,f)$ be a cycle for $\pi_{n}(S)$ as in Subsection \ref{sec1400}, where the constant map $f:M\to *$  is refined to   a stable framing  $TM\oplus (M\times \R^{k})\cong M\times \R^{n+k}$  of the tangent bundle. 
 A tangential   $Spin^{c}$-structure is now given by
a trivialisation \begin{equation}\label{eq20100}P(TM)\otimes Q(k)\cong M\times Spin^{c}(n+k)\ .\end{equation} We will  in fact assume that $P(TM)$ comes from a $Spin$-structure. In this case the Levi-Civita connection induces a canonical $Spin^{c}$-connection
$\tilde \nabla^{TM}$.

In order to apply  Theorem \ref{them2}
we must first choose a good geometrisation of $(M,f,\tilde\nabla^{TM})$. We will use the notation of 
 Subsection \ref{sec801}.  We can  choose the  manifold $M_{u}$ to be a point. The datum $(M_{u},f_{u},\tilde \nabla)$ (where $f_{u}:M_{u}\to *$ is constant and $\tilde \nabla$ is the trivial connection on the trivial bundle $P_{u}$) has  a unique geometrisation $\cG_{u}$. Let $h:M\to M_{u}$ be the constant map.  The $Spin^{c}$-bundle $h^{*}P_{u}$ is trivial. Hence the given trivialisation (\ref{eq20100}) 
 refines  $h$ to a $Spin^{c}$-map.  Using this refinement
we define
the geometrisation $\cG:=h^{*}\cG_{u}$ which turns out to be $\ell$-good for every $\ell\in \nat$.

 For $1\in K^{0}(*)\cong \Z$ we have $\cG_{u}(1)=1\in \hat K^{0}(M_{u})$. We now use Lemma \ref{lem906} in order to calculate
$\cG(1)\in \hat K^{0}(M)$. 
  By equation (\ref{eq1300}) we have
$$\cG(1)=1+a\left(\frac{\tilde \Td(\tilde \nabla^{triv},\tilde \nabla^{TM})}{\Td(\tilde \nabla^{TM})}\right)\ ,$$
where $\nabla^{triv}$ is the connection on $P(TM)\otimes Q(k)$ induced by the trivialization (\ref{eq20100}).
Let $\bV_{1}$ be the trivial one-dimensional geometric bundle on $M$. Then
$[\bV_{1}]=1\in \hat K^{0}(M)$ and in view of Equation (\ref{eq1450}) we must take the correction form 
$$\gamma_{1}:=\frac{\tilde \Td(\tilde \nabla^{triv},\tilde \nabla^{TM})}{\Td(\tilde \nabla^{TM})}\ .$$
We now specialize Theorem \ref{them2} to the present situation and obtain the intrinsic formula 
\begin{equation}\label{hfjhgerjfgjerferffref} i_{\R}(\eta^{an}([M,f]))(1)=[-\int_{M}\tilde \Td(\tilde \nabla^{triv},\tilde \nabla^{TM})]-\xi(\Dirac_{M})\in \R/\Z
\end{equation} This formula directly compares with the formula for $e_{\C}^{Adams}([M,f])$  derived by specializing \cite[Thm 4.14]{MR0397798}.

\subsection{$\rho$-invariants and the index theorem for flat bundles}\label{sec501}

We consider a closed odd-dimensional $Spin^{c}$-manifold $M$. The geometry on $M$ is given by a 
  Riemannian metric and the  $Spin^{c}$-extension of the Levi-Civita connection.
 If $\bV=(V,\nabla^{V},h^{V})$ is a  flat hermitean vector bundle of dimension $k$ on $M$,
%
%
 then the difference of reduced $\eta$-invariants
 \begin{equation}\label{eq3000}\rho(\Dirac_{M},\bV):=\xi(\Dirac_{M}\otimes \bV)-k\xi(\Dirac_{M})\end{equation}
 is invariant under variations of the geometry of $M$.
 The $\rho$-invariant  is thus a differential topological invariant of the $Spin^{c}$-manifold $M$ with a flat hermitean bundle $\bV$.
 
  The $\rho$-invariant is a classical  example of a topological invariant derived from the $\eta$-invariant  which has been studied a lot. 
  For example,  it has been    used    successfully    to detect elements in $Spin^{c}$-bordism groups of classifying spaces of finite cyclic groups \cite{MR870805}, \cite{MR883375}. We refer to these reference for examples of explicit calculations of $\rho$-invariants.

The precise homotopy theoretic description of $\rho$-invariants is given by  
 the index theorem for flat bundles \cite[Thm. 5.3]{MR0397799}.
The goal of the following discussion is to explain the relation of  the relation $\eta^{an}=\eta^{top}$ shown in Theorem \ref{indthm} with the index theorem for flat bundles.
Roughly speaking,   this goes as follows. The index theorem for flat bundles is about the pairing of the $K$-homology class represented by the $Spin^{c}$-Dirac operator with the torsion $K$-cohomology classes obtained from the flat bundle, while our index theorem considers the pairing of a torsion $K$-homology class with $K$-theory classes. Clearly the case of intersection is when both classes are torsion.


 We first translate the data of the $Spin^{c}$-manifold $M$ of odd dimension $n$ with a flat hermitean bundle $\bV$ into the bordism picture.  Let $U(k)^{\delta}$ denote the unitary group equipped with  the discrete topology.
 Its classifying space $BU(k)^{\delta}$ is universal for flat hermitean vector bundles
 of dimension $k$. We consider the bordism group based on the Thom spectrum of
   $$B:=BSpin^{c}\times BU(k)^{\delta}\stackrel{\pr}{\to} BSpin^{c}\ . $$
 
 We let $(M,f)$ be a cycle for $\pi_{n}(MSpin^{c})$ and consider in addition a map
  $g:M\to BU(k)^{\delta}$  which classifies a  flat bundle $\bV$. In this way we get a class
  $$[M,(f,g)]\in \pi_{n}(MB)\ .$$ We assume that this class is torsion in order to apply 
  the universal $\eta$-invariant.

We consider the $K$-theory class $\lambda_{k}\in K^{0}(BU(k)^{\delta})$ of the  universal $\C^{k}$-bundle on $BU(k)^{\delta}$ and    the projection
\begin{equation}\label{eq12000}q:BSpin^{c} \times  BU(k)^{\delta}\to BU(k)^{\delta}\ .\end{equation} 
Since $\ch(\lambda_{k}-k)=0$ the evaluation against the difference $q^{*}\lambda_{k}-k$ provides a well-defined homomorphism
$$\ev_{q^{*}\lambda_{k}-k} :Q_{n}(MB)\to \Q/\Z\ .$$



\begin{lem} 
We assume that $[M,(f,g)]\in \pi_{n}(MB)$ is a torsion class. Then we have the equality
  \begin{equation}\label{rihsqw}\ev_{q^{*}\lambda_{k}-k}(\eta^{an}([M,(f,g)])) =\rho(\Dirac_{M},\bV)\ .\end{equation}
   \end{lem}
   \proof
    We are going to use the notation introduced in Subsection \ref{sec81}. As an intermediate step we choose, for a suitable non-vanishing integer $l$,  a zero bordism $(W,(F,G))$ of the union of $l$  copies of the cycle $(M,(f,g))$ with $Spin^{c}$-geometry. The geometric bundle $\bU$ is then the flat hermitean bundle classified by $F$, and we have by Definition \ref{etaandef}
$$\ev_{q^{*}\lambda_{k}-k}(\eta^{an}([M,(f,g)]))=[\frac{1}{l}\ind(\Dirac_{W}\otimes \bU)]-[\frac{k}{l}\ind(\Dirac_{W})]\ .$$ If we use (\ref{eq870}) instead, then we express this evaluation
in terms of an integral over local data on $W$ and the reduced $\eta$-invariants. Because of $\ch(\nabla^{U})=k$
the local contributions cancel out and we remain with
$$\ev_{q^{*}\lambda_{k}-k}(\eta^{an}([M,(f,g)]))=\rho(\Dirac_{M},\bV)\ .$$
\hB

We now calculate the topological version of the universal $\eta$-invariant explicitly.
We again assume that $x=[M,(f,g)]\in \pi_{n}(MB)$ is a torsion element and let $\tilde x\in \pi_{n+1}(K\Q/\Z\wedge MB)$  as in \eqref{eq3}. By definition of $\eta^{top}$ we get
the equality  \begin{equation}\label{yx45}\ev_{q^{*}\lambda_{k}-k}(\eta^{top}([M,(f,g)]))=\langle \Thom^{K}(q^{*}\lambda_{k}-k),\tilde x \rangle \ .\end{equation}
The right-hand side of \eqref{rihsqw} is the analytic side of the index theorem for flat bundles in  \cite[Thm. 5.3]{MR0397799}.
The topological side of the index theorem for flat bundles  \cite[Thm. 5.3]{MR0397799} is not given as the pairing of a $K\Q/\Z$-homology class with a $K$-theory class, but rather by a pairing between a $K$-homology class and  a $K\R/\Z$-cohomology class. In the following we rewrite  the right-hand side of (\ref{yx45}) in this way.   

There exists a class $\Lambda_{k} \in K\R/\Z^{-1}(BU(k)^{\delta})$ with $\partial \Lambda_{k}=\lambda_{k}-k$, where  $\partial$ is the  Bockstein operator
$\partial:  K\R/\Z^{-1}(BU(k)^{\delta})\to K^{0}(BU(k)^{\delta})$, and which is uniquely characterised  by the following property:
If $h:N\to BU(k)^{\delta}$ is a map from a smooth manifold $N$, then
$h^{*}(\Lambda_{k})=[\bU]-k$, where
$\bU$ is the  unitary flat bundle classified by $h$. In order to interpret this equality we employ the identification $
    K\R/\Z^{-1}(N)\cong \hat K_{flat}^{0}(N) $ given by \eqref{eq700hhh}. 

%

Recall that $q$ denotes the projection  (\ref{eq12000}). We use Lemma
  \ref{lem900} for the first equality in the chain  $$\ev_{q^{*}\lambda_{k}-k}(\eta^{top}(x))=\langle\Thom^{K\R/\Z}(q^{*}\Lambda_{k}),\epsilon (x)\rangle=\langle q^{*}\Lambda_{k}, \Thom_{K}(\epsilon(x))\rangle=\langle [\bV]-k, [M]_{K}\rangle\ ,$$
  where $[M]_{K}$ denotes the $K$-theory fundamental class of the $Spin^{c}$-manifold $M$. 
The right-hand side of this equality  is the  topological side  of the index theorem for flat bundles of 
\cite[Thm. 5.3]{MR0397799}. 
The following Corollary now immediately  follows
from the equality $\eta^{top}=\eta^{an}$.
 \begin{kor}\label{cor45}
Let $n\in \nat$ be odd and  $M$ be a closed $n$-dimensional $Spin^{c}$-manifold with a flat hermitean $k$-dimensional vector bundle
 $\bV$. 
 We assume in addition that the corresponding class
 $[M,(f,g)]\in \pi_{n}(M(BSpin^{c}\times BU(k)^{\delta})$ is torsion.
 Then we have the following equality in $\R/\Z$:
\begin{equation}\label{eq4000}\rho(\Dirac_{M},\bV)=\langle[\bV]-k,[M]_{K}\rangle \end{equation}
\end{kor}
 In this way Theorem \ref{them2} implies a special case of
\cite[Thm. 5.3]{MR0397799}. Let us again remark, that by \cite[Thm. 5.3]{MR0397799}
the equality (\ref{eq4000}) holds true without the additional assumption that
 $[M,(f,g)]$ is a torsion class.

\subsection{Algebraic $K$-theory}\label{alghhh}

In this Subsection we use the universal $\eta$-invariant in order detect algebraic $K$-theory classes of $\C$. We will observe that  for odd $n\in \nat$ the well-known homomorphism
$$\varepsilon:K_{n}(\C)\to \Q/\Z$$ (see \eqref{ephvjhee}) 
 can be obtained from an appropriate evaluation of the universal $\eta$-invariant. Our main result is
Theorem \ref{them1000v} which provides a  formula for $\varepsilon$ in terms of  geometric cycles for $K$-theory classes. 
We will explain how the results of \cite{MR0482758} and  \cite{westburyjones} can be interpreted as constructions with the universal $\eta$-invariant.

For $n\in \nat$ the algebraic $K$-theory groups $K_{n}(\C)$ of the field $\C$ are defined as the homotopy groups of the connective  algebraic $K$-theory spectrum $K(\C)$. This spectrum is connected with classifying spaces  through Quillen's $+$-construction (see \cite[Ch. 3]{MR505692} for a detailed description) 
  \begin{equation}\label{eq1350}p:BGL(\C^{\delta})\to BGL(\C^{\delta})^{+} \end{equation} by the
 equivalence \begin{equation}\label{theequivb}\Omega^{\infty}K(\C)\cong \Z\times BGL(\C^{\delta})^{+}\ .\end{equation}

There exists a class
$$\Lambda_{0}\in K\R/\Z^{-1}(BGL(\C^{\delta})^{+})$$ which is uniquely characterised by the following property: If $k\in \nat$ and $g_{k}:N\to BGL(k,\C^{\delta})$ is a map from a smooth manifold, then
$$\iota_{*}(p\circ i_{k}\circ g_{k})^{*}\Lambda_{0}=[\bV]-k-a(\Ree(\tilde \ch(\nabla^{V,u},\nabla^{V})))\in K\C/\Z^{-1}(N)\ , $$
where on the right-hand side we use the identification of $K\C/\Z^{-1}(N)$ with the flat part of the  complex version of differential $K$-theory mentioned in Remark \ref{remark222}, $i_{k}:BGL(k,\C^{\delta})\to BGL(\C^{\delta})$ is the canonical map,  $\bV=(V,\nabla^{V})$ is the flat complex vector bundle of dimension $k$ classified by $g_{k}$, $\nabla^{V,u}$ is some choice of unitarisable connection on $V$, and $\iota:K\R/\Z^{-1}\to K\C/\Z^{-1}$ is the canonical map.  
We further define
$$\Theta_{0}:=\partial \Lambda_{0}\in K^{0}(BGL(\C^{\delta})^{+})\ , $$
where    $\partial:K\R/\Z^{-1}\to K^{0}$ is the Bockstein operator.


We will consider the Thom spectrum $MB$ associated to the projection $$B:=BSpin^{c}\times BGL(\C^{\delta})^{+}\to BSpin^{c}\ .$$  Furthermore we let
$q:B\to BGL(\C^{\delta})^{+}$ be the projection. 
Let $n\in \nat$ be odd.
Since $\ch(\Theta_{0})=0$, by  Lemma 
 \ref{lem1} 
the second map in the following definition of $\varepsilon_{0}$ is well-defined:
 $$\varepsilon_{0}:\pi_{n}(MB)_{tors}\stackrel{\eta^{top}}{\to} Q_{n}(MB)\stackrel{\ev_{q^{*}\Theta_{0}}}{\to} \Q/\Z\ .$$

We consider a cycle $(M,f)$ for $\pi_{n}(MSpin^{c})$ and a flat complex vector bundle $\bV$ on $M$ of dimension $k$. It is classified by a map $g_{k}:M\to BGL(k,\C^{\delta})$, and we consider the induced map $$g_{0}:=p\circ i_{k} \circ g_{k}:M\to BGL(\C^{\delta})^{+} .$$  Then
$(M,(f,g_{0}))$ is a cycle for $\pi_{n}(MB)$.

\begin{theorem}\label{them1000}   We assume that 
$[M,(f,g_{0})]\in \pi_{n}(MB)_{tors}$. If $\bV$ carries a flat hermitean metric, then
\begin{equation}\label{eq991v}
\varepsilon_{0}([M,(f,g_{0})])=\rho(\Dirac_{M},\bV) \ .\end{equation}
In general we have  
\begin{equation}\label{eq992v}
\varepsilon_{0}([M,(f,g_{0})])=\langle[M]_{K}, g_{0}^{*}\Lambda_{0}\rangle\ . \end{equation}
 \end{theorem}
\proof 
We first prove the general case (\ref{eq992v}) by the following chain of equalities
\begin{eqnarray*}\varepsilon_{0}([M,(f,g_{0})])&\stackrel{Def}{=}&\ev_{q^{*}\Theta_{0}}(\eta^{top}([M,(f,g_{0})]))\\&\stackrel{Lemma \:\ref{lem900}}{=}&\langle \Thom^{K\R/\Z}(q^{*}\Lambda_{0}), \epsilon([M,f,g_{0}])\rangle\\&=&\langle g_{0}^{*}\Lambda_{0}, [M]_{K}\rangle\ .\end{eqnarray*}

%
In the unitary case we observe that $g_{0}^{*}\Lambda_{0}=[\bV]-k$.  The equality (\ref{eq991v}) now follows from 
(\ref{eq4000}) and
the chain of equalities
$$\varepsilon_{0}([M,(f,g_{0})])\stackrel{(\ref{eq992v})}{=}\langle g_{0}^{*} \Lambda_{0}, [M]_{K}\rangle=\langle [\bV]-k, [M]_{K}\rangle \stackrel{(\ref{eq4000})}{=} \rho(\Dirac_{M},\bV)\ .$$
Note that in the present paper we have shown (\ref{eq4000}) under the assumption  that 
$ [M,(f,g_{k}^{u})]$ is a torsion class in $\pi_{n}(M(BSpin^{c}_{+}\times BU(k)^{\delta}))$,
where $g^{u}_{k}:M\to BU(k)^{\delta}$  classifies the hermitean flat bundle $\bV$.
 Since this might not be the case in general we have to appeal to the proof of this formula  (\ref{eq4000}) without such an assumption
given in \cite[Thm. 5.3]{MR0397799}.
\hB 

We have a canonical map of spectra
$$\Theta:K(\C)\to K\ .$$ Assume again that $n\in \nat$ is odd. Since
$\ch\circ \Theta=\dim$ is concentrated in degree zero, by Lemma 
 \ref{lem1}  the evaluation in the second map of the following definition is well-defined
 \begin{equation}\label{ephvjhee}\varepsilon:K_{n}(\C)_{tors}\stackrel{\eta^{top}}{\to} Q_{n}(K(\C))\stackrel{\ev_{q^{*}\Theta}}{\to}\Q/\Z\ .\end{equation}
We know from  Suslin \cite[Thm 4.9]{MR772065}  that the map $\varepsilon$ induces an isomorphism
 \begin{equation}\label{eq2100} K_{n}(\C)_{tors}\cong \left\{\begin{array}{cc}\Q/\Z&n\:odd\\0&n\:even\end{array}\right.
\end{equation} 
in the case  that $n$ is odd, and that the kernel of $\varepsilon$
is a uniquely divisible group.

Our goal is to provide a formula for $\varepsilon$ using geometric cycles for algebraic $K$-theory classes of $\C$.  Let $(M,f)$ be a cycle for $\pi_{n}(S)$, where $f$ stands for the constant map $M\to *$  refined by a stable normal framing. Furthermore, let $\bV$ be a flat complex vector bundle of dimension $k$ classified by a map $g_{k}:M\to BGL(k,\C^{\delta})$. Using the equivalence \eqref{theequivb} we define
$$g:=(\{0\}\times \id)\circ p\circ i_{k}\circ g_{k}:M\to \Omega^{\infty} K(\C)\ .$$   
In this way we get a class
$[M,(f,g)]\in \pi_{n}(\Sigma^{\infty}_{+}\Omega^{\infty}K(\C))$. Employing the canonical map
$u:\Sigma^{\infty}_{+}  \Omega^{\infty}K(\C) \to K(\C)$  we can form the class  $$u_{*}[M,(f,g)]\in K_{n}(\C)\ .$$

\begin{rem}{\rm
Every element $x\in K_{n}(\C)_{tors}
$ can be represented geometrically in this way.  More precisely
 there exists $(M,(f,g))$ as above such that
$x=u_{*}[M,(f,g)]$ and $[M,(f,g)]\in \pi_{n}(\Sigma^{\infty}_{+}  \Omega^{\infty}K(\C))_{tors}$.
In order to see this first represent $x$ by a map
$\gamma:S^{n}\to \Omega^{\infty} K(\C)$. Using the standard normal framing  of $S^{n}$ this map represents 
a torsion class $[S^{n},(f,\gamma)]\in \pi_{n}(\Sigma^{\infty}_{+} \Omega^{\infty} K(\C))_{tors}$. Since $\Sigma_{+}^{\infty}(p)$
is an equivalence by the universal property of the $+$-construction there exists a class
$z\in \pi_{n}(\Sigma^{\infty}_{+} BGL(\C^{\delta}))_{tors}$ such that
$p_{*}(z)=[S^{n},(f,\gamma)]$. Finally, the class $z$ can be represented in the form
$z=[M,(f,i_{k}\circ g_{k})]$ for some $k\in \nat$. Using this data we get
$u_{*}[M,(f,g)]=x$. 
}
\end{rem}

\begin{theorem}\label{them1000v} 
We assme that $[M,(f,g)]\in \pi_{n}(\Sigma^{\infty}_{+}  \Omega^{\infty}K(\C))_{tors}$.
 If $\bV$ carries a flat hermitean metric, then
\begin{equation}\label{eq991}
\varepsilon(u_{*}[M,(f,g)])=\rho(\Dirac_{M},\bV)+k e^{Adams}_{\C}([M,f]) \ .
\end{equation}
In general we have  
\begin{equation}\label{eq992}
\varepsilon(u_{*}[M,(f,g)])=\langle g_{0}^{*}\Lambda_{0} , [M]_{K}\rangle+ke^{Adams}_{\C}([M,f])\ . \end{equation}
 \end{theorem}
\proof
We again start with the general case (\ref{eq992}).
We calculate
\begin{eqnarray}\varepsilon(u_{*}[M,(f,g)])&=&\ev_{\Theta}(\eta^{top}(u_{*}[M,(f,g)]))\label{calghgiue}
 \\&=&\ev_{u^{*}\Theta}(\eta^{top}([M,(f,g)]))\nonumber\\&=&\ev_{ i_{k}^{*}u^{*}\Theta } ( \eta^{top}([M,(f,g_{0})]) ) \nonumber\\& =&\ev_{ \Theta_{0}+k}(\eta^{top}([M,(f,g_{0})]))\nonumber \ .\end{eqnarray}  We can now use the unstable case (\ref{eq992v}) and (\ref{gleichheit}) in order to deduce
 $$\varepsilon(u_{*}([M,(f,g)]))=\langle g^{*}_{0}\Lambda_{0}, [M]_{K}\rangle+ke^{Adams}_{\C}([M,f])\ .$$

The unitary case (\ref{eq991}) now follows from (\ref{eq992}) and (\ref{eq991v}). \hB

 \begin{rem}{\rm
We now show how one can deduce a special case of  \cite[Thm A]{westburyjones} from
 (\ref{eq992}). In \cite{westburyjones} an algebraic $K$-theory class is constructed from
 a homology sphere $M$ of dimension $n$ and a representation
 $\alpha:\pi_{1}(M)\to GL(k,\C^{\delta})$. 
 One  gets an induced map 
 $$\tilde g:M\stackrel{\alpha}{\to} BGL(k,\C^{\delta})\stackrel{i_{k}}{\to} BGL(\C^{\delta})$$ to which Quillen's $+$-construction is applied. The fundamental group of a homology sphere is perfect so that  the $+$-construction
 $M^{+}$ of $M$ is homotopy equivalent to a simply-connected homology sphere, hence to $S^{n}$.
 Thus we get a map
 $g^{+}: S^{n}\simeq M^{+} \stackrel{\tilde g^{+}}{\to}  BGL(\C^{\delta})^{+}$ which represents a class
 $[g^{+}]\in K_{n}(\C)$. The homology sphere $M$  admits a stable normal framing (see e.g. \cite{MR0148075} or \cite[Lem. 1]{MR0287549}) which refines the constant map $f:M\to *$.
 We further use  the composition
  $g:M\to M^{+}\stackrel{\tilde g^{+}}{\to}  BGL(\C^{\delta})^{+}$ in order to define  the class
 $[M,(f,g)]\in \pi_{n}(\Sigma^{\infty}_{+}\Omega^{\infty}K(\C)_{+})$ such that $u_{*}[M,(f,g)]=[g^{+}]$.

We consider the map of fibre sequence
$$\xymatrix{\Sigma^{-1}K\ar[r]\ar[d]&\Sigma^{-1}HP\C\ar[d]\\K^{rel}(\C)\ar[d]\ar@{.>}[r]&\Sigma^{-1}K\C/\Z\ar[d]\\
K(\C)\ar[r]^{\Theta}\ar[d]^{\Theta}&K\ar[d]^{\ch}\\
K\ar[r]^{\ch}&HP\C}\ .$$
 (the left column defines the relative $K$-theory spectrum $K^{rel}(\C)$). For odd $n\in \nat$ the dotted horizontal arrow
induces a map
\begin{equation}\label{eq9001}e:K_{n}(\C)\cong \frac{\pi_{n}(K^{rel}(\C))}{\im(\pi_{n}(\Sigma^{-1}K)\to \pi_{n}(K^{rel}(\C)))}\to \pi_{n}(\Sigma^{-1} K\C/\Z)\cong \C/\Z\ .\end{equation}

Let $\bV^{flat}$ denote the flat vector bundle determined by the representation $\alpha$.
The statement of  \cite[Thm A]{westburyjones} is the equality
\begin{equation}\label{eq9002}e([g^{+}])=\rho_{\C}(\Dirac_{M},\bV^{flat})\end{equation}
in $\C/\Z$. The subscript $\C$ on the right-hand side indicates  a complex version of the $\rho$-invariant defined for flat vector bundles without requiring a flat hermitean metric.
We have the equality
 $$\Ree(\rho_{\C}(\Dirac_{M},\bV^{flat}))=\langle g^{*}\Lambda_{0},[M]_{K} \rangle\ .$$
  
 Note that $(\C/\Z)_{tors}=\Q/\Z\subset \R/\Z\subset \C/\Z$.
If $[g^{+}]\in K_{n}(\C)$ is a torsion class, then using Lemma \ref{lem900}  we get the identification
$$e([g^{+}])=\varepsilon([g^{+}])-k e_{\C}^{Adams}([M,f])\in \R/\Z$$ by an argument which is similar to the calculation   \eqref{calghgiue}. 
Hence if $[g^{+}]$ is a torsion class, then the real part of \eqref{eq9002} is equivalent to \eqref{eq992}.
 
  }\end{rem}

\begin{rem}{\rm  
 Let us comment on the fact that Adams' $e$-invariant appears on the right-hand sides of (\ref{eq991}) and (\ref{eq992}).
Note that $K(\R)$ is a ring spectrum with unit $\epsilon_{K(\R)}:S\to K(\R)$. 
The unit  induces a homomorphism $\pi_{*}(S)\to K_{n}(\R)$. Since the image $\im(J)$ of the $J$-homomorphism
(\ref{eq30001}) is a well-known summand of $\pi_{*}(S)$ it was an interesting question to determine
its image under the unit $\epsilon_{K(\R)}$. Let us consider the case
$$\frac{\Z}{\frac{B_{m}}{4m}\Z}\cong \im(J)_{4m-1}\subseteq \pi_{4m-1}(S)$$
for $m\in \nat$.
In \cite{MR0482758} it has been shown that this piece goes  injectively to algebraic $K$-theory. 
This was deduced from the following  two facts:
\begin{enumerate}
\item $\im(J)_{4m-1}$ is detected by (the real version of) Adams' $e$-invariant
$e:\pi_{4m-1}(S)\to \Q/\Z$.
\item The $e$-invariant has a factorisation over the analog 
$$K_{n}(\R)(\dots)\to \pi_{n}(\Sigma^{-1}KO\C/\Z)$$ 
of the homomorphism (\ref{eq9001}).
\end{enumerate}
The complex case of 
this factorisation is easily seen from Theorem 
 \ref{them1000v}. 
 In fact the elements in $\im(J)_{4m-1}$ can be represented by cycles of the form
 $u_{*}[S^{4m-1},(f,g)]$, where $f$ carries  a normal framing obtained from the standard
 framing by twisting with  an element of $\pi_{4m-1}(O)$, and $g$ is obtained from the classifying map of the trivial one-dimensional bundle $\bV$. Then we get from Theorem \ref{them1000v}, \eqref{eq991},
$$\varepsilon(u_{*}[S^{4m-1},(f,g)])=e_{\C}^{Adams}([S^{4m-1},f])\ .$$
  
 }\end{rem}

\subsection{$String$-bordism}\label{subsec55}

In this Subsection we describe the connection of the present paper with constructions 
 in \cite{2009arXiv0912.4875B}. In this reference we introduced an invariant $b^{an}$ of elements of the bordism group $\pi_{4m-1}(MString)$   using a formula which shares a lot of similarities with the intrinsic formula (\ref{eq1000}) for $\eta^{an}$. We will recall the definition of the complex version $b^{an}_{\C}$ during the proof of Theorem \ref{thm7}.

 The spectrum of topological modular forms $\tmf$ has been constructed by Miller, Goerss and Hopkins, and in an alternative way by Lurie, see the survey \cite{MR2597740}. It is related to 
$K$-theory and $String$-bordism by a factorisation of the Witten genus
\begin{equation}\label{gleluzfe}\xymatrix{MString\ar@/^2cm/@{.>}[rrr]^{\sigma^{\C}_{Witten}}\ar@/^1cm/@{.>}[rr]^{\sigma_{Witten}}\ar[r]^{\sigma_{AHR}}&\tmf\ar@/_{1cm}/@{.>}[rr]_{W_{\C}}\ar[r]^{W}\ar[r]&KO[[q]]\ar[r]&K[[q]]}\ ,\end{equation}
where $\sigma_{AHR}$ is the String-orientation of $\tmf$ constructed by  Ando, Hopkins and Rezk.
 One of the interesting features of the restriction of $b^{an}$ to the kernel of the map to $Spin$-bordism is that it has a factorisation over $\sigma_{AHS}$.  Since $b^{an}$ is calculable in interesting cases it can be used to detect the $\tmf$-class represented by a closed $String$-manifold.

Our goal here is  Theorem \ref{thm7}  which gives the precise relation between $b^{an}$ and the universal $\eta$-invariant. In Corollary \ref{thm7v}
  we show one of the conjectures stated  in  \cite{2009arXiv0912.4875B} asserting that the factorisation of $b^{an}$ over topological modular forms holds true on the whole $String$-bordism group, i.e. we get rid of the restriction to the kernel to the $Spin$-bordism.    Formally our proof is complete in dimensions $8m-1$, while in dimensions $8m-5$ we lose some two-torsion since in the present paper we work with complex $K$-theory instead of real $K$-theory. We strongly believe that the relevant part of the theory has a real version which does prove the case in dimension $8m-5$ completely, too.

In  Theorem  \ref{prop77} we show how the Riemannian geometry on a $String$-manifold together
with a geometric string structure give rise to a geometrisation, and we derive the corresponding intrinsic formula for $b^{an}$.  We consider the proof of Theorem \ref{prop77} as a model for many other situations where a construction of a geometrisation is required.

We start with recalling the definition of $String$-bordism.  The homotopy type of the space $BString$ is defined as a stage in the Whitehead tower of
 $BO$:
$$ \xymatrix{BString\ar[d]^{p}&\\BSpin\ar[d]\ar[r]^{\frac{p_{1}}{2}}&K(\Z,4)\\BSO\ar[d]\ar[r]^{w_{2}}&K(\Z/2\Z,2)\\BO\ar[r]^{w_{1}}&K(\Z/2\Z,1)}\ .$$
The space $BString=BO\langle 8\rangle$ is just a low instance of a whole tower of higher connected coverings  $BO\langle n\rangle$
of the classifying space $BO$. Starting with $BString$ these higher spaces are no longer
associated to classical families of compact Lie groups. 


In Lemma \ref{lem77} we demonstrate that a connection on a $Spin$-principal bundle gives rise to a
geometrisation. While the connection on the principal bundle allows to define connections on all associated vector bundles, the geometrisation partially keeps this information in terms of the differential $K$-theory classes represented by these vector bundles with connections. 
The geometrisation associated to a geometric $String$ structure in this sense replaces the theory of
connections on the non-existing principal bundle with structure group $String$.
We think that the methods used in the case of $BString=BO\langle 8\rangle$ can easily be adapted to the higher stages $BO\langle n\rangle$.

We let $MString$ be the Thom spectrum associated to the map
$BString\to BSpin\stackrel{can}{\to} BSpin^{c}$.
The $String$-bordism spectrum $MString$ is rationally even (see \cite{MR2465746}, \cite{MR1297530}, \cite{MR1455523} for more calculations), so that for $m\in \nat$ and $n=4m-1$ the group
$\pi_{n}(MString)$ is torsion. Hence the universal $\eta$-invariant gives a map
$$\eta^{top}=\eta^{an}:\pi_{n}(MString)\to Q_{n}(MString)\ .$$
We will show that we can obtain $b^{an}$ from the universal $\eta$-invariant by defining an interesting homomorphism out of the universal target $Q_{n}(MString)$. It involves  evaluations against a collection of elements
$R_{k}(\lambda^{String})\in K^{0}(BString)$ for all $k\ge 0$. It is useful to organise this collection in a formal power series
$$R(\lambda^{String}):=\sum_{k\ge 0} q^{k} R_{k}(\lambda^{String})\in K[[q]]^{0}(BString)$$
which we will describe in the following.
By  $K[[q]]$ we denote the multiplicative cohomology theory (resp. the corresponding spectrum) which associates to a  space
$Y$ the  ring $$K[[q]]^{*}(Y):=K^{*}(Y)[[q]]$$ of formal power series with  coefficients in $K^*(Y)$.
The following constructions with real vector bundles are standard in the theory of the Witten genus (\ref{eqwitten}), compare e.g. with  \cite{MR1189136}, \cite{2009arXiv0912.4875B}.
 Given a real vector bundle $V\to Y$
we consider the element
$R(V)\in K[[q]]^{0}(Y)$
defined by
\begin{equation}\label{eq1105}R(V):=\sum_{k=0}^{\infty} R_{k}(V)q^{k}\ ,\end{equation}
where $R_{k}(V)$ is the $K$-theory class of the virtual bundle
given by the coefficient in front of $q^{k}$ in the expansion of 
$$\prod_{k\ge 1} (1-q^{k})^{2\dim(V)} \bigotimes_{k\ge 1} \Sym_{q^{k}}(V\otimes_{\R}\C)\ ,$$
where $$\Sym_{p}(V):=\bigoplus_{k\ge 0} p^{k}\Sym^{k}(V)\ .$$
The transformation $V\mapsto R(V)$ is exponential, i.e. it satisfies
$R(V\oplus W)=R(V)\cup R(W)$. Moreover, it  has values in the group of multiplicative units $K[[q]]^{0}(Y)^{\times}$ since the power series \eqref{eq1105} starts with $1$. In view of the universal property of $KO^{0}$ it therefore extends to a natural transformation
$$R:KO^{0}(Y)\to K[[q]]^{0}(Y)^{\times} \ .$$

The composition
$$BString\to BO\stackrel{x\mapsto (0,x)}{\longrightarrow} \Z\times BO\cong \Omega^{\infty} KO$$
classifies the universal class $\theta^{String}\in KO^{0}(BString)$.
We fix $n=4m-1$ and let 
\begin{equation}\label{eq1100}\lambda^{String}_{n}:=n-\theta^{String}\in KO^{0}(BString)\ .\end{equation}
If $(M,f)$ is a cycle for $\pi_{n}(MString)$, then we have the equality
\begin{equation}\label{eq2001v}[TM]+1 = f^{*}\lambda^{String}_{n+1}\end{equation}
in $KO^{0}(M)$.
We  have well-defined classes $R_{k}(\lambda^{String}_{n+1})\in K^0(BString)$ for all $k\ge 0$
and therefore $R(\lambda_{n+1}^{String})\in K[[q]]^{0}(BString)$. With this notation the Witten genus
$$\sigma_{Witten}^{\C}:\pi_{n+1}(MString)\to \pi_{n+1}(K[[q]])$$ satisfies 
\begin{equation}\label{eqwitten}\sigma_{Witten}^{\C}(x)=\langle \Thom^{K}(R(\lambda_{n+1}^{String})),\epsilon(x)\rangle \ .\end{equation}
We use the superscript $\C$ in order to indicate that we work with the image of the $KO[[q]]$-valued Witten genus in complex $K$-theory $K[[q]]$. 

We organize the evaluations of  $Q_{n}(MString)$ against the family of classes $R_{k}(\lambda^{String}_{n+1})$ into a formal power series and
define a homomorphism
\begin{equation}\label{eq1100v}W:\Hom(K^{0}(BString),\pi_{n+1}(K\Q/\Z))\to \Q/\Z[[q]]:=\prod_{k\ge 0} \Q/\Z\: q^{k}\end{equation}
by
$$W(\phi):=\sum_{k\ge 0} \ev_{R_{k}(\lambda^{String}_{n+1})}(\phi)q^{k}\ .$$
 
The homomorphism (\ref{eq1100v}) does not yet factorize over  $Q_{n}(MString)$ since it does not vanish on the subgroup $U^{\prime}$ appearing in  \eqref{spacerep1}. In order to get such a factorisation we must replace the target $\Q/\Z[[q]]$ of $W$ by  the quotient  by a subgroup which contains $W(U^{\prime})$. This subgroup will be defined using modular forms.
 We let $\cM^{R}_{2m}$ denote the space of modular forms for $SL(2,\Z)$ of weight $2m$
 whose $q$-expansion have coefficients in the subring $R\subseteq \C$ (see \cite{MR1189136} for an introduction).
 In particular, we let
 $$\cM^{\Q}_{2m}[[q]]\subseteq \Q[[q]]$$ be the finite-dimensional vector space of $q$-expansions of rational modular forms $\cM^{\Q}_{2m}$ of weight $2m$. Its image in $\Q/\Z[[q]]$ will be denoted by 
 $\overline{\cM^{\Q}_{2m}[[q]]}$.
 We define
 \begin{equation}\label{eq1101}T_{2m}:=\frac{\Q/\Z[[q]]}{\overline{\cM^{\Q}_{2m}[[q]]}}\ .\end{equation}
 Up to the replacement of $\Q$ by $\R$ this is exactly the group defined in \cite[Def. 1.1]{2009arXiv0912.4875B}.
 \begin{lem}\label{lem1092}
 The composition of (\ref{eq1100v}) with the projection to the quotient (\ref{eq1101})
 induces a well-defined map
 $$\bar W:Q_{4m-1}(MString)\to T_{2m}\ .$$
 \end{lem}
 \proof
 We must show that under this composition 
 the subgroup $U$ defined in (\ref{sec31u})
 is mapped to $\overline{\cM^{\Q}_{2m}[[q]]}$. By (\ref{eqwitten}) we have for 
 $y\in \pi_{n+1}(MString)$ that  
$$   \langle \Thom^{K}(R(\lambda_{n+1}^{String})), \epsilon(y)\rangle=\sigma^{\C}_{Witten}(y)\in \pi_{n+1}(K[[q]])\cong \Z[[q]]\ .$$ We now use the fact that
 the Witten genus has values in $\cM_{2m}^{\Z}[[q]]\subset \Z[[q]] $.
 More generally, for $y\in \pi_{n+1}(MString\Q)$ we get 
 $$\langle  \Thom^{K}(R(\lambda_{n+1}^{String})),  \epsilon(y)\rangle\in \cM_{2m}^{\Q}[[q]]\ .$$  
 This shows that $\bar W(U^{\prime})\subseteq \overline{\cM^{\Q}_{2m}[[q]]}$. \hB 
 
 In  \cite[Sec 3.3]{2009arXiv0912.4875B} we have constructed  homomorphisms
 $$b^{an}: \pi_{4m-1}(MString)\to T_{2m}\ ,\quad b^{top}:A_{4m-1}\to T_{2m}$$
 where $$A_{4m-1}:=\ker(\pi_{4m-1}(MString)\to \pi_{4m-1}(MSpin))\ .$$
 Since in the present paper we work we complex $K$-theory as opposed to real $K$-theory  in \cite[Sec 3.3]{2009arXiv0912.4875B}   we define
 \begin{equation}\label{eq2769}b^{*}_{\C}:=\left\{\begin{array}{cc}b^{*} &m\:even\\ 2b^{*}&m\:odd\end{array}\right.\ ,\quad *\in \{an,top,\tmf\}\ .\end{equation}
 \begin{theorem}\label{thm7}
 We have the equalities
$$\bar W\circ \eta^{top}_{|A_{4m-1}}=b^{top}_{\C}$$
 and
 $$\bar W\circ \eta^{an}=b^{an}_{\C}$$
\end{theorem} 
\proof
We extend the map $MString\to MSpin$ to a fibre sequence
$$\Sigma^{-1} MSpin\to \cA \to MString \to MSpin$$ which defines the spectrum
 $\cA$.
The smash product of the fibre sequence with the fibre sequence \eqref{eq631}
yields the following quadratic diagram
$$\xymatrix{\Sigma^{-2}MSpin\Q\ar[d]\ar[r]&\Sigma^{-1}\cA\Q\ar[d]\ar[r]&\Sigma^{-1}MString\Q\ar[d]\ar[r]&{}^{\hat w}\:\Sigma^{-1}MSpin\Q\ar[d]\\
\Sigma^{-2}MSpin\Q/\Z\ar[d]\ar[r]&\Sigma^{-1}\cA\Q/\Z\ar[d]\ar[r]&{}^{\hat x}\:\Sigma^{-1}MString\Q/\Z\ar[d]\ar[r]&{}^{\hat y}\:\Sigma^{-1}MSpin\Q/\Z\ar[d]\\
\Sigma^{-1}MSpin\ar[d]\ar[r]&{}^{\hat z}\:\cA\ar[d]\ar[r]&{}^{x}\:MString\ar[d]\ar[r]&{}^{0}\:MSpin\ar[d]\\{}^{-\hat w}\:
\Sigma^{-1}MSpin\Q\ar[r]&{}^{\tilde z}\cA\Q\ar[r]&{}^{0}\:MString\Q\ar[r]&MSpin\Q}$$
We start with $x\in A_{4m-1}\subseteq \pi_{4m-1}(MString)$. This element goes to zero if it is mapped
to the right or down.  The class
$\bar W(\eta^{top}(x))$
is represented by the power series
\begin{equation}\label{eq1102}\sum_{k\ge 0} \langle \Thom^{K}(R_{k}(\lambda^{String}_{4m})), \epsilon(\hat x) \rangle q^{k}\in  \Q/\Z[[q]] \ .\end{equation}
Note that we can define classes
$\theta^{Spin}$ and $\lambda_{n}^{Spin}:=n-\theta^{Spin}\in KO^{0}(BSpin)$
analogously to (\ref{eq1100}). Then we have equalities of the  evaluations
\begin{eqnarray}\label{eq1103} \langle \Thom^{K}(R_{k}(\lambda^{String}_{4m})), \epsilon(\hat x)\rangle&=& \langle\Thom^{K}(R_{k}(\lambda^{Spin}_{4m})), \epsilon(\hat y) \rangle\\
&=& [\langle\Thom^{K}(R_{k}(\lambda^{Spin}_{4m})),\epsilon(\hat w)\rangle]\in \Q/\Z\ ,\nonumber\end{eqnarray}
where the elements $\hat y$ and $\hat w$ are images and lifts as indicated in the above diagram, 
and where we use the compatibility of the $K$-theory Thom isomorphisms for $MString$ and $MSpin$.
In the construction of $b^{top}$ in \cite[Sec 4.1]{2009arXiv0912.4875B}  we go the other way. We 
first lift $x$ to $\hat z$ which maps to $\tilde z$ which is then again lifted to
$\Sigma^{-1}MSpin\Q$. It is a general fact of such a diagram chase in a product of fibre sequences that
modulo the obvious ambiguities 
this element in the lower left  corner is the negative of $\hat w$ from the upper right corner.
By the definition of $b^{top}$ in \cite[Sec 4.1]{2009arXiv0912.4875B} we see that 
$b_{\C}^{top}(x)$ is represented by
$$\sum_{k\ge 0} [\langle \Thom^{K}(R_{k}(\lambda^{Spin}_{4m})), \epsilon (\hat w)\rangle] q^{k}\in \Q/\Z[[q]]\ .$$
Combining this with (\ref{eq1102}) and (\ref{eq1103}) we see that
$$\bar W\circ \eta^{top}_{|A_{4m-1}}=b_{\C}^{top}\ .$$
This proves the first assertion of Theorem \ref{thm7}.

We now show the second.
Let $x=[M,f]\in \pi_{n}(MString)$ be an $l$-torsion element represented by the cycle $(M,f)$ and $(W,F)$ be a zero-bordism of the union of $l$ copies of $(M,f)$. 
The $Spin^{c}$-structures come from $Spin$-structures so that the Levi-Civita connections have  canonical $Spin^{c}$-extensions $\tilde \nabla^{TM}$ and $\tilde \nabla^{TW}$.
In view of Equation (\ref{eq2001v}) the class
  $f^{*}R(\lambda^{String}_{4m})\in K[[q]]^{0}(M)$ can be represented by a formal power series of $\Z/2\Z$-graded bundles $R(TM\oplus 1)$ associated to the tangent bundle. The Riemannian metric and the Levi-Civita connection turn $TM$ into a geometric bundle. The construction of $R(TM\oplus 1)$ therefore produces a formal power series of    geometric bundles $\bR(TM\oplus 1)$. 

The construction of $b^{an}$ involves the choice of a geometric $String$-structure $\alpha$ on $M$. 
This notion has been introduced in \cite{2009arXiv0906.0117W}. In particular, a geometric string structure  produces a
form $H_{\alpha}\in \Omega^{3}(M)$ with the property that 
\begin{equation}\label{eq2303}2dH_{\alpha}=p_{1}(\nabla^{TM,LC})\ .\end{equation}
In the following we use  characteristic forms associated to certain power series $$\tilde \Phi\ , \:\:\Phi\ , \:\:\Theta\:\in \Q[[q]][b,b^{-1}][[b^{2}p_{1},b^{4}p_{2},\dots]]\ .$$  We refer to   \cite[Sec 3.3]{2009arXiv0912.4875B}  or equation  (\ref{er10v}) for an explicit definition. In the notation of the latter we have
$$\Phi:=\Phi_{R(\lambda_{4m}^{Spin})}\ ,\quad  \tilde \Phi:=\tilde \Phi_{R(\lambda_{4m}^{Spin})}\ , \quad 
\Theta:=\Phi-p_{1}\tilde \Phi\ .$$
In the present paper we distribute the powers of $b$ such that
$\Phi$ and $\Theta$ have total degree zero, and $\tilde \Phi$ has total degree $-4$.
The notation $\tilde \Phi(\tilde \nabla^{TM})$ is as in (\ref{er10}).
We start with the representative of $b_{\C}^{an}(x)$ given in
 \cite[Def 4.1]{2009arXiv0912.4875B}  
\begin{equation}\label{eq1007}[2\int_{M} H_{\alpha}  \wedge \tilde \Phi(\tilde \nabla^{TM})]- \xi(\Dirac_{M}\otimes \bR(TM\oplus 1))]\in \R/\Z[[q]]\ ,\end{equation}
where here and below we ignore the power $b^{2m}$.
We use the APS index formula (\ref{eq4}) in order to express the reduced $\eta$-invariant appearing in  (\ref{eq1007}).
Using the equality 
$$\Phi(\tilde \nabla^{TW})=\Td(\tilde \nabla^{TM})\wedge \ch(\nabla^{\bR(TM\oplus 1)})$$
we get
$$(\ref{eq1007})= [2\int_{M} H_{\alpha}  \wedge \tilde \Phi(\tilde \nabla^{TM})-\frac{1}{l}\int_{W} \Phi(\tilde \nabla^{TW})]+[\frac{1}{l}\ind(\Dirac_{W}\otimes \bR(TW))_{APS}]\ .$$
We now use Stokes' theorem and the relation (\ref{eq2303})  in order to calculate
\begin{eqnarray*}\lefteqn{
2\int_{M} H_{\alpha}  \wedge \tilde \Phi(\tilde \nabla^{TM})-\frac{1}{l}\int_{W}\Phi(\tilde \nabla^{TW})}\hspace{1cm}&&\\&=&\frac{1}{l}\int_{W} \left( p_{1}(  \nabla^{TW})\wedge \tilde \Phi(\tilde \nabla^{TM})-\Phi(\tilde \nabla^{TW})\right)\\&=&\frac{1}{l}\int_{W} \Theta(\tilde \nabla^{TW})\in \cM_{2m}^{\R}[[q]]\ .
\end{eqnarray*}
For the last inclusion we use the crucial fact that
$$p_{4m}(\Theta)\in \cM_{2m}^{\Q}[[q]][b^{2}p_{1},b^{4}p_{2},\dots]\subset  \Q[[q]][[b^{2}p_{1},b^{4}p_{2},\dots]] \ ,$$
see \cite[Sec. 3.3]{2009arXiv0912.4875B}.
Therefore
$$ [\frac{1}{l}\ind(\Dirac_{W}\otimes \bR(TW))_{APS}]\in \R/\Z[[q]]$$ is a representative of
$b_{\C}^{an}(x)\in T_{2m}$, too. But in view of Definition \ref{etaandef}
and the construction of $\bar W$ this is also a representative of 
$\bar W(\eta^{an}(x))\in T_{2m}$. This shows
$$\bar W\circ \eta^{an}=b^{an}_{\C}\ .$$
\hB

As a consequence of the equality $\eta^{an}=\eta^{top}$ shown in Theorem \ref{indthm}
we get another proof of \cite[Thm. 2.2]{2009arXiv0912.4875B}. \begin{kor}
We have the equality
$$b^{an}_{\C|A_{4m-1}}=b_{\C}^{top}\ .$$
\end{kor}

We now recall  from  \cite[Sec.4.3]{2009arXiv0912.4875B} the construction of  the homomorphism
$$b^{\tmf}:\pi_{4m-1}(\tmf)\to T_{2m}$$
which is very similar to that of $\eta^{top}$.
Note that $\pi_{4m-1}(\tmf)$ is a torsion group (see \cite{MR1989190}, \cite{MR2508200} for more calculations of $\pi_{*}(\tmf))$. Therefore an element $y\in \pi_{4m-1}(\tmf)$ can be lifted to an element
$\hat y\in \pi_{4m}(\tmf \Q/\Z)$. Then
$$b^{\tmf}(y):=[W(\hat y)]\in T_{2m}\ ,$$
where $[W(\hat y)]$ denotes the class in $T_{2m}$ of the element
$W(\hat y)\in \pi_{4m}(KO[[q]]\Q/\Z)\cong \Q/\Z[[q]]$ and $W$ is as in \eqref{gleluzfe}.
The complex version  $b_{\C}^{\tmf}$ of $b^{\tmf}$ is defined similarly by
$$b_{\C}^{\tmf}(y):=[W_{\C}(\hat y)]\in T_{2m}\ ,$$
or alternatively, by (\ref{eq2769}).

\begin{prop}
We have the equality
$$b_{\C}^{\tmf}\circ \sigma_{AHR}=\bar W\circ \eta^{top}:\pi_{4m-1}(MString)\to T_{2m}\ .$$
\end{prop}
\proof
If $x\in \pi_{4m-1}(MString)$ and $\hat x\in \pi_{4m}(MString\Q/\Z)$ is a lift, then
$$\langle \Thom^{K}(R(\lambda^{String}_{4m})), \epsilon(\hat x)\rangle\in \Q/\Z[[q]]$$ represents
$\bar W\circ \eta^{top}(x)$. We have already seen in the proof of Lemma \ref{lem1092} that this expression is equal to the Witten genus (extended to $\Q/\Z$-theory)
$$\langle\Thom^{K}(R(\lambda^{String}_{4m})), \epsilon(\hat x)\rangle=\sigma^{\C}_{Witten}(\hat x)\ .$$
The Witten genus (see \eqref{gleluzfe}) can now be decomposed as
$$\sigma^{\C}_{Witten} (\hat x)=W_{\C}(\sigma_{AHR}(\hat x))\ .$$
We can take $\sigma_{AHR}(\hat x)\in \pi_{4m}(\tmf\Q/\Z)$ as the lift of $\sigma_{AHR}(x)\in \pi_{4m}(\tmf)$ so that 
$W_{\C}(\sigma_{AHR}(\hat x))$ represents  $b_{\C}^{\tmf}(\sigma_{AHR}(x))$.
Hence we can conclude that
$$b_{\C}^{\tmf}\circ \sigma_{AHR}(x)=\bar W\circ \eta^{top}(x)\ .$$
\hB

  Using   Theorem \ref{indthm}  stating  the equality $\eta^{an}=\eta^{top}$ and $\bar W\circ \eta^{an}=b_{\C}^{an}$ given by Theorem \ref{thm7}) we get
\begin{kor}\label{thm7v}
We have
$$b^{an}_{\C}= b_{\C}^{\tmf}\circ \sigma_{AHR}\ .$$
\end{kor}
This proves the complex version of the conjecture 3 in 
 \cite[Sec 1.5]{2009arXiv0912.4875B}. 
 In fact, for even $m$ there is no difference between the real and complex case, but
 in the case of odd $m$ the complex version implies the real version up to two-torsion which was known before.
 We believe  that a real version of the present theory would prove the conjecture 
 completely.

 \bigskip
 
 The formula for $b^{an}$ given in  \cite[Sec 3.3]{2009arXiv0912.4875B} and reproduced here as (\ref{eq1007}) is an intrinsic formula which uses the notion of a geometric $String$-structure \cite{2009arXiv0906.0117W}. In the following we show that a geometric $String$-structure gives rise to a good geometrisation $\cG^{String}$ of $(M,f,\tilde \nabla^{TM})$ such that the intrinsic formula  Theorem \ref{them2}  specialises to the one for $b^{an}$.
 Since $String$-structures refine $Spin$-structures we start with the construction of a geometrisation for a $Spin$-structure.
 
We consider a cycle $(M,f)$ for $\pi_{n}(MSpin)$. 
 We are going to use a version of Subsection \ref{sec1677} for $Spin$-structures. If $V\to M$ is a real euclidean oriented vector bundle, then  the $Spin$-gerbe  $\underline{Spin(V)}$ of $V$ associates to each open subset $A\subseteq M$ the groupoid $Spin(V_{|A})$ of $Spin$-structures on the restriction of $V$ to $A$. This gerbe has the band
$\Z/2\Z$, and its isomorphism class is  classified by the Dixmier-Douady class 
$$\DD(\underline{Spin(V)})=w_{2}(V)\in H^{2}(M;\Z/2\Z)\ ,$$
the second Stiefel-Whitney class of $V$.

We choose a tangential $Spin$-structure on $TM$  given by a geometric $Spin$-structure
$P\in Spin(TM)$ together with a trivialisation 
\begin{equation}\label{eq65789}P\otimes \tilde f^{*}Q^{Spin}_{k}\cong Q(n+k)\ ,\end{equation}
where $\tilde f:M\to BSpin(k)$ is some factorisation of $f$, and $Q^{Spin}_{k}\to BSpin(k)$   denotes
the universal $Spin$-bundle.  It naturally induces a tangential $Spin^{c}$-structure by extension of structure groups along $Spin(l)\to Spin^{c}(l)$ (see \cite[Example D.5]{MR1031992}).

The Levi-Civita connection gives rise to a connection $\nabla^{TM}$ on $P$ which in turn has a natural  $Spin^{c}$-extension $\tilde \nabla^{TM}$.
We furthermore choose a connection $\nabla^{k}:=\nabla^{\tilde f^{*}Q^{Spin}_{k}}$ on $\tilde f^{*}Q^{Spin}_{k}$ and let $\tilde \nabla^{k}$ be its $Spin^{c}$-extension.
\begin{prop}\label{lem77}
There exists a  geometrisation $\cG^{Spin}$ of
$(M,f,\tilde \nabla^{TM})$ which is $\ell$-good for every $\ell\in \nat$.
\end{prop}
\proof
The connections $\tilde \nabla^{TM}$ and $\tilde \nabla^{k}$ together induce a $Spin^{c}$-extension 
$\tilde \nabla^{TM}\otimes \tilde \nabla^{k}$
of the connection $\nabla^{TM}\otimes \nabla^{k}$ 
on $P\otimes \tilde f^{*}Q^{Spin}_{k}$ which can be compared with a trivial connection using the isomorphism
(\ref{eq65789}). Therefore the transgression form
$\tilde \Td(\tilde \nabla^{TM}\otimes\tilde  \nabla^{k},\nabla^{triv})\in \Omega P^{-1}(M)/\im(d)$ is defined and
satisfies
$$d\tilde \Td(\tilde \nabla^{TM}\otimes \tilde \nabla^{k},\nabla^{triv})=\Td(\tilde \nabla^{TM})\wedge \Td(\tilde \nabla^{k})-1\ .$$
We let
\begin{equation}\label{eq87100v}\mu:= \Td(\tilde \nabla^{k})^{-1}\wedge \tilde \Td\tilde (\nabla^{TM}\otimes \tilde \nabla^{k},\nabla^{triv})\ .\end{equation}
For any space or spectrum $Y$ we let
$\bar K^{*}(Y)$ denote the completion  the topological group $K^{*}(Y)$ equipped with the profinite topology (see  \cite[Def. 4.9] {MR1361899} or Remark \ref{proffinite1}).
We have $BSpin:=\hocolim_{l} \: BSpin(l)$ and therefore  $\bar K^{*}(BSpin)\cong \lim_{l}\: \bar K^{*}(BSpin(l))$.
The completion theorem \cite{MR0259946} gives
$$ K^{*}(BSpin(k))= \bar K^{*}(BSpin(k))\cong  R(Spin(k))^{\hat{}}_{I_{Spin(k)}}\ .$$
We therefore get the following description of the completion of the $K$-theory of $BSpin$:
$$\bar K^{*}(BSpin)\cong \lim_{l} \: K^{*}(BSpin(k+l))\cong \lim_{l}\: R(Spin(k+l))^{\hat{}}_{I_{Spin(k+l)}}\ .$$

We fix an integer $l\ge 0$ and form  the $l$-fold stabilization 
$\tilde f^{*}Q^{Spin}_{k}\otimes Q(l)$ of $\tilde f^{*}Q^{Spin}_{k}$. This stabilization is a $Spin(k+l)$-principal bundle with the connection  $\nabla^{k}\otimes \nabla^{Q(l)}$.
Given a  representation $\rho$ of $Spin(k+l)$
we define a geometric bundle $\bV_{\rho}$ as an associated bundle to $\tilde f^{*}Q^{Spin}_{k}\otimes Q(l)$.  We define
$$G_{l}(\rho):=[\bV_{\rho}]-a(\mu\wedge \ch(\nabla^{V_{\rho}}))\in \hat K^{0}(M)\ .$$
We have chosen the form $\mu$ in (\ref{eq87100v}) such that the following equality holds true in $\Omega P^{0}_{cl}(M)$:
\begin{equation}\label{eq87100}\Td(\tilde \nabla^{TM})\wedge R(G(\rho))=\Td(\tilde \nabla^{k})^{-1}\wedge \ch(\nabla^{V_{\rho}})\ .\end{equation}
The map $\rho\mapsto G_{l}(\rho)$
extends to a map
$G_{l}:R(Spin(k+l))\to \hat K^{0}(M)$ by linearity.
This extension  annihilates the $2n+1$'th power
$I_{Spin(k+l)}^{2n+1}\subseteq R(Spin(k+l))$ of the dimension ideal.
In order to see this note that  if $\rho\in I^{p}_{Spin(k+l)}$ and  $2p>n$, then we have $\ch(\nabla^{V_{\rho}})=0$.
For those $\rho$ we have
$G_{l}(\rho)=[\bV_{\rho}]$, and this class is flat.
If $p>n$, then we have $[V_{\rho}]=0$ so that
$G_{l}(\rho)=a(\omega)$ for some $\omega\in HP\R^{-1}(M)$.
The product of a flat class with a class of this form vanishes by \eqref{compaprod}.
Hence $G_{l}(\rho)=0$ if $p>2n$.
The map $G_{l}$ thus further extends be continuity 
to
a map
$$\cG^{Spin(k+l)}:K^{0}(BSpin(k+l))\to \hat K^{0}(M)\ .$$
 
One now checks that for $l\ge 1$ we have
$$G_{l}(\rho)=G_{l-1}(\rho_{|Spin(k+l-1})\ .$$
In this way the maps $\cG^{Spin(k+l)}$ for the various $l$ are compatible. We consider the continuous map
$$\cG^{Spin}:\bar K^{0}(BSpin)\to K^{0}(BSpin(k))\stackrel{\cG^{Spin(k)}}{\to} \hat K^{0}(M)\ .$$

We now show that $\cG^{Spin}$ and $\cG^{Spin(k+l)}$ are geometrisations.   To this end we must show that they admit degree-perserving cohomological characters.  By their compatibility it suffices to consider $\cG^{Spin(k)}$.
It follows from (\ref{eq87100}) and $\Td(\tilde \nabla^{TM})=\hA(\nabla^{TM})$ (since $\tilde \nabla^{TM}$ is induced from a $Spin$-connection $\nabla^{TM}$) that
$\Td(\tilde \nabla^{TM})\wedge R(G_{0}(\rho))$ is the Chern-Weyl representative of the class
$\hA^{-1}\cup \ch([\rho])\in HP\Q^{0}(BSpin(k))$ associated to the connection
$\nabla^{k}$, where $[\rho]\in K^{0}(BSpin(k))$ is the class represented by $\rho$. Note that
\begin{equation}\label{eqf1}H^{*}(BSpin(k);\Q)\cong \Q[p_{1},p_{2},\dots,p_{r_{k}}]\end{equation}
is the polynomial ring generated by the universal Pontrjagin classes. \textcolor{black}{For the cohomological    character $$c_{\cG^{Spin(k)}}: HP\Q^{0}(BSpin(k))\to \Omega P_{cl}^{0}(M)$$  of $\cG^{Spin(k)}$ we choose the map which sends}   the class  $b^{2i}p_{i}\in HP\Q^{0}(BSpin(k))$ to the corresponding characteristic form $b^{2i}p_{i}(\nabla^{k})\in \Omega P^{0}_{cl}(M)$.
 This map clearly preserves degrees.
 
 \bigskip

  We now show that $\cG^{Spin}$ is $\ell$-good for every $\ell\in \nat$.  By an inspection of the construction of $\cG^{Spin}$  we observe that the connection of the map $f:M\to B$ with the stable normal bundle of $M$ has not been used. This map can be arbitrary if we replace $TM$ by some complement $\eta\to M$ of $\hat f^{*}\xi_{k}$ as in  Subsection \ref{sec18555} and choose some connection $ \nabla^{\eta}$ of the associated complementary $Spin$-principal bundle $P^{Spin}\in Spin(\eta)$ in the place of $\nabla^{TM}$. We obtain the $Spin^{c}$-bundle $P$ with connection
 $\tilde \nabla^{\eta}$ which replaces $\tilde \nabla^{TM}$ by extension of the structure group.

Let $\ell\in \nat $ be such that $\ell\ge n+1$.
We choose an $\ell$-connected approximation $f_{u}:M_{u}\to BSpin$ such that there is a factorisation of $f$ over a closed embedding $h:M\to  M_{u}$. As in Subsection \ref{sec801} we obtain
a natural refinement of $h$ to a $Spin^{c}$-map. Since $h$ is a closed embedding we can choose the
connections $\tilde \nabla^{u}$ on $P_{u}$ and $\nabla^{k,u}$ on $\tilde f_{u}^{*}Q^{Spin}_{k}$ such that
$h^{*}\tilde \nabla^{u}=\tilde \nabla^{TM}$ stably  and $h^{*}\nabla^{k,u}=\nabla^{k}$.
  We now define
$\cG^{Spin}_{u}$ as above. Then by construction
$\cG^{Spin}=h^{*}\cG_{u}^{Spin}$ since the correction forms (\ref{eqq4})
vanishes.
Since $\ell$ can be taken arbitrary large we see that $\cG^{Spin}$ is $\ell$-good for every $\ell\in \nat$. 
\hB 

 Note that the geometrizations $\cG^{Spin(k+l)}$ constructed in the proof of Proposition \ref{lem77}
 are $\ell$-good for every $\ell\in \nat$, too.

 \bigskip
 
Let   $p:BString\to BSpin$ be the natural map. We  consider a cycle
$(M,f)$ for $\pi_{n}(MString)$.  Then $(M,p\circ f)$ naturally becomes a cycle for $\pi_{n}(MSpin)$.
 
\begin{theorem}\label{prop77}
A choice of a geometric string structure
 $\alpha$ on $(\tilde f^{*}Q^{Spin}_{k},\nabla^{k})$ naturally determines 
a geometrisation   $\cG^{String}$ of $(M,f,\tilde \nabla^{TM})$ which is $\ell$-good for every $\ell\in \nat$.
For $\phi\in \bar K^{0}(BSpin)$ it is given by
$$\cG^{String}(p^{*}\phi)=\cG^{Spin}(\phi)-a(\nu_{\phi})$$
with
$$\nu_{\phi}:=
\Td(\tilde \nabla^{TM})^{-1}\wedge 2H_{\alpha}\wedge \tilde \Phi_{\phi}(\nabla^{m})\in \Omega P^{-1}(M)\ .$$ Here 
$\cG^{Spin}$ is as in Lemma \ref{lem77}, 
 $\tilde \Phi_{\phi}(\nabla^{m})$ is defined below in (\ref{er10}),
and $H_{\alpha}\in \Omega^{3}(M)$ is the three-form given by the geometric string structure. 
 \end{theorem}
 \proof
\color{black} 
We fix an integer $\ell\ge n+1$.
For a sufficiently large   integer $m\in \nat$, $m\ge 3$ we can assume that
 the map $f$ has a factorization $$M\stackrel{f_{m}}{\to} BString(m) \stackrel{\iota_{m}}{\to} BString$$
 such that $\iota_{m}$ is $\ell$-connected.
 We have a fibre sequence
$$K(\Z,3)\to BString(m)\stackrel{p}{\to} BSpin(m)\to K(\Z,4)\ .$$
By \cite{MR0231369} the reduced $K$-theory group
$\tilde K^{*}(K(\Z,3))$   is \textcolor{black}{uniquely} divisible and consists of phantom classes, i.e. classes which vanish when pulled-back to finite $CW$-complexes. This suggests the following proposition which is probably well-known, but we could not find a reference for it.


 \begin{prop}\label{pinv}
 For $m\ge 3$
the pull-back  \begin{equation}\label{barp}p^{*}:\bar K^{*}(BSpin(m)) \to   \bar K^{*}(BString(m)) \end{equation} along the projection 
$p:BString(m)\to BSpin(m)$ represents $\bar K^{*}(BString(m))$ as a completion of $\bar K^{*}(BSpin(m))$
with respect to the topology induced on $\bar K^{*}(BSpin(m))$ via $p^{*}$ by the profinite topology of $\bar K^{*}(BString(m))$.
\end{prop}

We defer the   proof of this proposition  to Section \ref{hdklqwjdlqwdwqd}. 
and  continue with the construction of a geometrisation $\cG^{String}$ which is $\ell$-good. We choose an $\ell+1$-connected approximation $f_{u}:M_{u}\to BString(m)$ such that we can factorize
$f_{m}:M\to BString(m)$  over the closed embedding $h:M\to M_{u}$.  
As in Subsection \ref{sec801} we obtain
a natural refinement of $h$ to a $Spin^{c}$-map. Since $h$ is a closed embedding we can choose the
connections $\tilde \nabla^{u}$ on $P_{u}$ and $\nabla^{m,u}$ on $f_{u}^{*}Q^{Spin}_{m}$ such that
$h^{*}\tilde \nabla^{u}=\tilde \nabla^{TM}$ stably and $h^{*}\nabla^{m,u}=\nabla^{m}$.
Geometric string structures behave as flexible as connections and metrics \cite{2009arXiv0906.0117W}.
We can therefore assume that there is a geometric string structure $\alpha_{u}$ on $(f_{u}^{*}Q^{Spin}_{m},\nabla^{m,u})$ which restricts to the geometric string structure $\alpha$ on $M$.

In the proof of Lemma \ref{lem77} we have  constructed  a geometrisation
  $$\cG^{Spin(m)}_{u}:\bar K^{0}(BSpin(m))\to \hat K^{0}(M_{u})\ .$$ of $(M_{u}, p\circ f_{u},\nabla^{m,u})$.  
As a first approximation   we define
$$\cG^{String(m)}_{u,0}: \textcolor{black}{\bar K^{0}(BSpin(m))}  \to \hat K^{0}(M_{u})$$ by $$\cG^{String(m)}_{u,0}( \phi ):=\cG_{u}^{Spin(m)}(\phi)\ , \quad \phi\in \bar K^{0}(BSpin(m))\ .$$  
\color{black}In view of Proposition \ref{pinv}  the homomorphism
$\cG^{String(m)}_{u,0}$ is defined on a group which  defines $\bar K^{0}(BString(m))$ by completion with respect to a certain topology.  If $\cG^{String(m)}_{u,0}$  would be continuous with respect    to this topology it would descend to a continuous homomorphism defined on $\bar K^{0}(BString(m))$. Note that $\cG^{String(m)}_{u,0}$ is not continuous and 
 the idea is now to add a correction term  in order to improve the continuity.

Note that via $p^{*}$ we can identify 
 $$H\Q^{*}(BString(m);\Q)\cong \Q[[p_{2},\dots,p_{r_{m}}]]$$ 
 with the  quotient ring of $H^{*}(BSpin(m);\Q)$ given in  (\ref{eqf1}) by setting $p_{1}= 0$.
The problem with continuity comes from the contribution of $p_{1}(\tilde \nabla^{u})$ to the
curvature of $\cG_{u,0}^{String(m)}$.
The idea is now to kill this contribution by a correction term given by a geometric string structure
 $\alpha_{u}$ on $(\tilde f_{u}^{*}Q^{Spin}_{m},\nabla^{m,u})$. The geometric string structure provides the form $H_{\alpha_{u}}\in \Omega^{3}(M_{u})$ with the property that $2dH_{\alpha_{u}}=p_{1}(\nabla^{m,u})$ (see \ref{eq2303}).

For a formal power series 
$$\Lambda\in \Q[b,b^{-1}][[b^{2}p_{1},b^{4}p_{2},\dots]]$$ we define a new formal power series
\begin{equation}\label{cont2}\tilde \Lambda:=\frac{\Lambda-i_{p_{1}=0} \Lambda}{p_{1}}\in \Q[b,b^{-1}][[b^{2}p_{1},b^{4}p_{2},\dots]]\ .\end{equation}
In other words, the power series
$\tilde \Lambda$ is  $p_{1}^{-1}$ times the sum of those monomials of $\Lambda$ which contain $p_{1}$.
Since the periodic rational cohomology of any space $Y$ is complete, i.e.  we have
$HP\Q^{*}(Y)\cong \overline{HP\Q}^{*}(Y)$,  the Chern character factorizes over the completion of $K$-theory as
$\ch:\bar K^{0}(BSpin(m))\to HP\Q^{0}(BSpin(m))$. 
Let $\phi\in \bar K^{0}(BSpin(m))$. Then we define
\begin{equation}\label{er10v}\Phi_{\phi}:=\Td^{-1}\cup \ch(\phi)\in HP\Q^{0}(BSpin(m))\cong  \Q[[b^{2}p_{1},b^{4}p_{2},\dots,b^{2r_{m}}p_{r_{m}}]]\end{equation}
and obtain $\tilde \Phi_{\phi}$ as described above. We define the form
$$\nu_{u,\phi}:=\Td(\tilde \nabla^{u})^{-1}\wedge 2H_{\alpha_{u}}\wedge \tilde \Phi_{\phi}(\nabla^{m,u})\in \Omega P^{-1}(M_{u})\ ,$$
where we use the abbreviation
\begin{equation}\label{er10}\tilde \Phi_{\phi}(\nabla^{m,u}):= \tilde \Phi_{\phi}(p_{1}(\nabla^{m,u}),p_{2}(\nabla^{m,u}),\dots)\ .\end{equation}
We can now define the map 
$$\cG^{String(m)}_{u}: \textcolor{black}{\bar K^{0}(BSpin(m))}\to \hat K^{0}(M)$$ by the following prescription:
$$\cG_{u}^{String(m)}( \phi ):=\cG^{String(m)}_{u,0}( \phi) -a(\nu_{u,\phi})\ .$$
We further define
$$\cG^{String(m)}:=h^{*}\cG_{u}^{String(m)}\ .$$ 
Unfortunately we can not verify directly that $\cG_{u}^{String(m)}$ is continuous, but we have the following Lemma.
\begin{lem}\label{cont1}
The map
$$\cG^{String(m)}: \textcolor{black}{\bar K^{0}(BSpin(m))} \to \hat K^{0}(M)$$
extends by continuity to  $\bar K^{0}(BString(m))$.  \textcolor{black}{
The continuous extension (for which we use the same symbol $\cG^{String(m)}$)  admits the degree-preserving cohomological character given by $b^{2i}p_{i}\mapsto b^{2i}p_{i}(\nabla^{m})$ for all $i\in \{2,\dots,r_{m}\}$ and
is therefore a geometrisation of $(M,f,\tilde \nabla^{TM})$.}
 \end{lem}
 \proof
 \textcolor{black}{
 Since $BString(2m)$ has finite  skeleta  the profinite topology on $\bar K^{*}(BString(2m))$
 has a countable basis of neighbourhoods of zero so that we can check continuity of homomorphisms using sequences converging to zero.}
 Let us consider a sequence $(\phi_{k})$ in   $\bar K^{0}(BSpin(m))$  such that
$p^{*}\phi_{k}\to 0$ in the profinite topology of $\bar K^{0}(BString(m))$ as $k\to \infty$. We must show that there exists a $k_{0}\in \nat$ such that for all $k\ge k_{0}$ we have $\cG^{String(m)}( \phi_{k})=0$.
Let $t:N\to BString(m)$ be a compact $\dim(M_{u})\textcolor{black}{+5}$-connected approximation.
We can choose a $k_{0}\in \nat$ such that for all $k\ge k_{0}$ we have
$t^{*}p^{*}\phi_{k}=0$ and $t^{*}p^{*}(\Td^{-1}\wedge \ch(\phi_{k}))=0$. 
Since the pull-back
$t^{*}:H^{*}(BString(m);\Q)\to H^{*}(N;\Q)$ 
\color{black} is injective in degrees $\le \dim(M_{u})+4$
we see that $\Phi_{\phi_{k}}$ is a  formal power series of terms of  homogeneity   
$\ge \dim(M_{u})+5$ (here we count the topological degree, i.e. disregard the degree of $b$). 
Consequently, $i_{p_{1}=0} \Phi_{\phi_{k}}$ is a  formal power series of terms of homogeneity   
$\ge \dim(M_{u})+1$.

  We now calculate using (\ref{eq87100}) and (\ref{cont2}) that for $\phi\in \bar K^{0}(BSpin(m))$ 
 \begin{eqnarray}\Td(\tilde \nabla^{u})\wedge R(\cG_{u}^{String(m)}( \phi ))&=&
 \Td(\tilde \nabla^{u})\wedge R(\cG_{u}^{Spin(m)}(\phi))-p_{1}(\tilde \nabla^{m,u})\wedge \tilde \Phi_{\phi}(\nabla^{m,u})\nonumber\\
 &=&\Phi_{\phi}(\nabla^{m,u}) -p_{1}(\tilde \nabla^{m,u})\wedge \tilde \Phi_{\phi}(\nabla^{m,u})\nonumber\\
 &=&(i_{p_{1}=0} \Phi_{\phi})(\nabla^{m,u})\label{eqcont30}\ .
 \end{eqnarray}
We conclude that \textcolor{black}{
$\Td(\tilde \nabla^{u})\wedge R(\cG_{u}^{String(m)}( \phi_{k}  ))=0$ for $k\ge k_{0}$ since
$(i_{p_{1}=0} \Phi_{\phi})(\nabla^{m,u})$ would be a form on $M_{u}$ of a degree which    exceeds the dimension. Hence 
 for   $k\ge k_{0}$   the class}
$\cG^{String(m)}_{u}( \phi_{k} )$ is flat and in the kernel of $I:\hat K^{0}(M_{u})\to K^{0}(M_{u})$. We conclude that
$$\cG^{String(m)}_{u}( \phi_{k} )\in HP\R^{-1}(M_{u})/\im(\ch)\ .$$
An $\ell+1$-connected map induces an isomorphism in ordinary cohomology in degrees $\le \ell$.
Since $BString(m)$ is rationally even the odd-dimensional real  cohomology of the $\ell+1$-connected approximation $M_{u}$ is concentrated in degrees $\ge \ell+1$. 
Since $\dim(M)\le \ell$  the restriction  $h^{*}:HP\R^{-1}(M_{u})\to HP\R^{-1}(M)$ is trivial.
This implies that
$\cG^{String(m)}( \phi_{k} )=h^{*} \cG^{String(m)}_{u}( \phi_{k} )=0$ for all $k\ge k_{0}$.

The assertion about the cohomological character follows from  the relation
$$\Td(\tilde \nabla^{TM})\wedge R(\cG^{String(m)}(\phi))=(i_{p_{1}=0} \Phi_{ \phi })(\nabla^{m})$$
derived from (\ref{eqcont30}). This finishes the proof of Lemma \ref{cont1}. \hB

In order to show that this geometrisation $\cG^{String(m)}$   is  $\ell$-good we must show that
$ \cG_{u}^{String(m)}$ is continuous itself. To this end we argue similarly by representing this geometrisation
as a pull-back from a $\dim(M_{u})+1$-connected approximation of $BString(m)$. 

We now define   the geometrization $\cG^{String}$ of $[M,f,\tilde \nabla^{TM}]$ by
$$\cG^{String}:=\iota_{m,*} \cG^{String(m)}\ .$$  Since $\iota_{m}$ is $\ell$-connected, this geometrization is  $\ell$-good. If $m^{\prime}\in \nat $ satisfies $m\le m^{\prime}$, then it follows from the compatibilty of the family of geometrizations $(\cG_{u}^{Spin(m^{\prime})})_{m^{\prime}\ge m}$ that 
$$\cG^{String(m^{\prime})}(\phi)=\cG^{String(m)}(\phi_{|BString(m)})\ , \quad \phi\in \bar K^{0}(BString(m^{\prime}))\ .$$ 
Therefore $\cG^{String}$ does not depend on the choice of $m$.
In particular, it is $\ell$-good for all $\ell$. 
This finishes the proof of Theorem \ref{prop77}. \hB

\color{black}

We now specialise  Theorem
\ref{them2} in order to derive an intrinsic formula for $$b^{an}([M,f])=\bar W\circ \eta^{an}([M,f])\in T_{2m}\ .$$ The connection $\nabla^{k}$ on the $Spin(k)$-principal bundle $\tilde f^{*}Q^{Spin}_{k}\to M$ turns the real vector bundle $\tilde f^{*}\xi^{String}_{k}$ into a geometric bundle $\bN_{k}$. It is a geometric representative of the stable normal bundle of $M$, hence the notation. 
We have  $$R([TM]+1)=R(n+1+k-[\tilde f^{*} \xi^{String}_{k}])\in K[[q]]^{0}(M)\ .$$ Therefore we get an interpretation of $\bR(n+1+k-\bN_{k})$ as a virtual geometric representative of $R([TM]+1)$ (which differs from $\bR(TM+ 1)$ used in (\ref{eq1007}) since we work with the geometry on the normal bundle).  By construction we have
$$\cG^{String}(\lambda^{String}_{n+1})=[\bR(n+1+k-\bN_{k})]+a(\nu_{R(\lambda^{Spin}_{n+1})})\in \hat K^{0}(M)[[q]]\ .$$
In other words, the correction form for $[\lambda^{String}_{n+1}]\in K[[q]]^{0}(BString)$ is given by
$$\gamma_{R(\lambda^{String}_{n+1})}=\nu_{R(\lambda^{Spin}_{n+1})}=
\Td(\tilde \nabla^{TM})^{-1}\wedge 2H_{\alpha}\wedge \tilde  \Phi_{R(\lambda^{Spin}_{n+1})}(\nabla^{k})\in \Omega P^{-1}(M)[[q]]\ .$$
By  
Theorem
\ref{them2}
the composition $\bar W\circ \eta^{top}([M,f])\in T_{2m}$ is now represented by the formal power series
$$[- \int_{M}2H_{\alpha}\wedge \tilde \Phi_{R(\lambda^{Spin}_{n+1})}(\nabla^{k})]-\xi(\Dirac_{M}\otimes \bR(n+1+k-\bN_{k}))\in \R/\Z[[q]]\ .$$
This is the version of  (\ref{eq1007})  using the normal bundle geometry on the twisting bundles.

\subsection{The Crowley-Goette invariants} \label{subseccg}

In this Subsection we want to show how some of the  Eells-Kuiper or Kreck-Stolz type invariants from geometric topology can be understood from the point of view of the universal $\eta$-invariant. Our approach can be described schematically as follows. In a first step, called pre-construction, we translate the geometric topology data into elements of a suitable bordism group $MB$. We then apply the universal $\eta$-invariant. Finally, the desired invariant is obtained  by a suitable evaluation against classes in $K^{0}(MB)$.

In detail we will consider the example of the  Crowley-Goette invariant  recently introduced in \cite{2010arXiv1012.5237C} for $S^{3}$-principal bundles on certain $n=4m-1$-dimensional manifolds. We start with recalling the definitions from  \cite{2010arXiv1012.5237C}.  Since in the present paper we decided to work with $Spin^{c}$-bordism and complex Dirac operators we will define the variant $t^{\C}_{M}$ which coincides with the Crowley-Goette  invariant for even $m$ and is its double for odd $m$.
Let $S^{3}$ be the group of unit quaternions and $BS^{3}$ be its classifying space. The set of homotopy classes $[M,BS^{3}]$ is in natural bijection with the set of isomorphism classes of $S^{3}$-principal bundles
on $M$ denoted in \cite{2010arXiv1012.5237C} by $\Bun(M)$.

Let $M$ be a closed $n$-dimensional $Spin$-manifold
such that $H^{3}(M;\Q)=0$ and $H^{4}(M;\Q)=0$. Then the Crowley-Goette invariant is defined as a certain function
$$t_{M}:\Bun(M)\to \Q/\Z\ .$$
In the following we recall the intrinsic formula \cite[(1.9)]{2010arXiv1012.5237C} for $t_{M}$.
Note that $$HP\Q^{*}(BS^{3})\cong \Q[b,b^{-1}][[b^{2}c_{2}]]\ ,$$ and by the completion theorem \cite{MR0259946} we have the isomorphism
$$K^{0}(BS^{3})\cong R(S^{3})^{\hat{}}_{I_{S^{3}}}\ .$$
We let $\rho$ be the  representation of $S^{3}$ by left-multiplication on $\mathbb{H}\cong \C^{2}$.
The representation $\rho$ gives rise to a class
$[\rho]\in K^{0}(BS^{3})$ and a power series
$\ch([\rho])\in \Q[b,b^{-1}][[b^{2}c_{2}]]^{0}$ of total degree zero. There exists a unique power series
$\tilde \Phi\in  \Q[b,b^{-1}][[b^{2}c_{2}]]^{-4}$ of total degree $-4$ such that
$2-\ch([\rho])=c_{2}\ \tilde \Phi$.

Let $\tilde g\in \Bun(M)$ and $R\to M$ be a $S^{3}$-bundle classified by $\tilde g$.
We choose a connection $\nabla^{R}$ on $R$. For every unitary representation $(\lambda,V_{\lambda})$ of $S^{3}$ we let $E_{\lambda}:=P\times_{S^{3},\lambda}V_{\lambda}$ be the vector bundle associated to $R$ and $\lambda$. It comes with a natural hermitean metric $h^{E_{\lambda}}$. The connection $\nabla^{R}$ induces a connection $\nabla^{E_{\lambda}}$ which preserves  $h^{E_{\lambda}}$. In this way we get a geometric bundle $\bE_{\lambda}:=(E_{\lambda},h^{E_{\lambda}},\nabla^{E_{\lambda}})$. By our assumptions on the rational cohomology  of $M$ the Chern-Weyl representative  $c_{2}(\nabla^{R})$ of $c_{2}$ is exact, and there exists a unique element $\hat c_{2}(\nabla^{R})\in \Omega^{3}(M)/\im(d)$ such that $d\hat c_{2}(\nabla^{R})=c_{2}(\nabla^{R})$.
We define  $\tilde \Phi(\nabla^{R})\in \Omega P_{cl}^{-4}(M)$ by  replacing $c_{2}$ by  $c_{2}(\nabla^{R})$ in the power series $\tilde \Phi$.

We choose a Riemannian metric on $M$ which induces the Levi-Civita connection on $TM$.
Furthermore we choose the $Spin^{c}$-structure induced by the $Spin$-structure. We then get a natural
$Spin^{c}$-extension $\tilde \nabla^{TM}$ of the Levi-Civita connection.
The complex version
$$t_{M}^{\C}:\Bun(M)\to \R/\Z$$ of $t_{M}$ is now given by 
\cite[(1.9)]{2010arXiv1012.5237C}
\begin{equation}\label{cg2}t^{\C}_{M}(\tilde g):=[\int_{M}\Td(\tilde \nabla^{TM})\wedge \hat c_{2}(\nabla^{R})\wedge \tilde \Phi(\nabla^{R})]-2\xi(\Dirac_{M}) + \xi(\Dirac_{M}\otimes \bE_{\rho})\in \R/\Z\ .\end{equation}
 To be precise, the value of the integral belongs to $\R[b,b^{-1}]^{-4}$ which will be identified with $\R$ using the generator $b^{2}$.

In order to relate the Crowley-Goette invariant with the universal $\eta$-invariant we are led to consider a  bordism theory of $Spin^{c}$-manifolds with $S^{3}$-bundles with rationally trivial second Chern class. 
This bordism theory is constructed homotopy theoretically  as follows.
We have a fibre sequence 
$$ K(\Q,3)\to K(\Q/\Z,3) \stackrel{\partial}{\to} K(\Z,4)\to K(\Q,4)$$
of Eilenberg-MacLane spaces. We define the space $X$ by the following homotopy pull-back
\begin{equation}\label{cg1}\xymatrix{X\ar[d]^{q}\ar[r]&K(\Q/\Z,3)\ar[d]^{\partial}\\
BS^{3}\ar[r]^{c_{2}}&K(\Z,4)}\ .\end{equation}
Since $c_{2}$ is a rational isomorphism and $K(\Q/\Z,3)$ is rationally trivial we see that the space $X$ is rationally contractible. 

We consider the Thom spectrum $MB$ associated to the projection $$B:=BSpin^{c}\times X\to BSpin^{c}\ .$$  We conclude that for $n=4m-1$
$$\pi_{n}(MB)\otimes \Q\cong \pi_{n}(MSpin^{c}\wedge X_{+})\otimes \Q\cong \pi_{n}(MSpin^{c}\Q\wedge X_{+})\cong  \pi_{n}(MSpin^{c}\Q)\cong 0\ .$$
It follows that
$$\pi_{n}(MB)_{tors}=\pi_{n}(MB)$$ 
so that the universal $\eta$-invariant is defined on the whole $Spin^{c}$-bordism group of $X$: 
$$\eta^{top}=\eta^{an}:\pi_{n}(MB)\to Q_{n}(MB)\ .$$

  The pre-construction for the Crowley-Goette invariant is a map 
 $$s_{M}:\Bun(M)\to \pi_{n}(MB) $$
 given by the following Lemma.
 \begin{lem}\label{lemcg5}
A pair $((M,f),\tilde g)$ of  a cycle $(M,f)$ for $\pi_{n}(MSpin^{c})$     with $H^{3}(M;\Q)=0$ and $H^{4}(M;\Q)=0$,  and a map $\tilde g\in \Bun(M)$ gives naturally rise to a class
$s_{M}(\tilde g):=[M,(f,g)]\in \pi_{n}(MB)$. 
\end{lem}
\proof
The main point is to show that $\tilde g:M\to BS^{3}$ has a natural lift to $g:M\to X$ in the diagram
(\ref{cg1}).
The rationalisation of $\tilde g^{*}c_{2}$ vanishes so that there exists a class
$\hat c_{2}\in H^{3}(M;\Q/\Z)$ such that $\partial \hat c_{2}=\tilde g^{* }c_{2}$.
This lift is unique up to the image of a rational class of degree $3$, hence unique by our assumption.
The map
$\tilde g$ and the lift $\hat c_{2}:M\to K(\Q/\Z,3)$ together determine the lift $g:M\to X$.
 \hB

Our next task is to determine the element in $K^{0}(MB)$ at which want to evaluate.
To this end we calculate
the $K$-theory  $K^{*}(X)$.
We have a fibration
$$ X\stackrel{q}{\to} BS^{3}$$
with fibre $K(\Q,3)$. \textcolor{black}{Note that  $\Q$ is a countable abelian group.  Furthermore, the space $ BS^{3}$ has a $CW$-structure with 
 finite skeleta $(BS^{3}_{k})_{k\ge 0}$ and $\lim_{k}^{1} K^{*}(BS^{3}_{k})=0$ by  \cite{MR0259946}.
  We can apply  
  Proposition
 \ref{pinv1}  and see  that
\begin{equation}\label{cg15}q^{*}:\bar K^{*}( BS^{3} )\to \bar K^{*}( X)\end{equation}
represents $\bar K^{*}( X)$ as a completion
of $\bar K^{*}(  BS^{3} )$ with respect to the topology induced from the profinite topology on $\bar K^{*}( X)$ via $q^{*}$.}
 The domain of this map can be calculated the using the completion theorem \cite{MR0259946}. 
The element  
$2-\rho$ of the representation ring $R(S^{3})$ generates the dimension ideal $I_{S^{3}}$. If we let
$A:=2-[\rho]\in \bar K^{0}(BS^{3})$, then we have
$$\bar K^{*}(BS^{3})\cong \Z[[A]]\ .$$ 
Since $X$ is rationally contractible we have
$$\ch^{*}(q^{*}A)=q^{*}\ch(A)=0\ .$$
Let  $p:BSpin^{c}\times X\to X$ be the projection.  By Lemma \ref{lem1} the evaluation
$$\ev_{p^{*}q^{*}A}:Q_{n}(MB)\to \Q/\Z$$
is well-defined. We define
$$\varepsilon:= \ev_{p^{*}q^{*}A}\circ \eta^{top}:\pi_{n}(MB)\to \Q/\Z\ .$$

The following proposition clarifies the relation between $t_{M}^{\C}$ and the universal $\eta$-invariant.
\begin{prop}\label{cg78}
If $(M,f)$  is a cycle for $\pi_{n}(MSpin^{c})$ which satisfies    
$H^{3}(M;\Q)=0$ and $H^{4}(M;\Q)=0$, then we have the relation
$$t^{\C}_{M}=\varepsilon\circ s_{M}:\Bun(M)\to \Q/\Z\ .$$
\end{prop}
\proof
It is an instructive exercise in the use of geometrisations to
derive an intrinsic formula for the composition $\varepsilon\circ s_{M}$ which can be compared with the formula 
(\ref{cg2}) for $t_{M}^{\C}$. In a first step we must approximate the space $X$ by  spaces with finite skeleta. Note that we can write (compare with (\ref{eq4333}) for the connecting maps)
$$K(\Q/\Z,3):=\hocolim_{l}\: K(\Z/l\Z,3)\ .$$
If we define
$X_{l}$ by the pull-back 
\begin{equation}\label{cg133}\xymatrix{X_{l}\ar[d]^{q_{l}}\ar[r]&K(\Z/l\Z,3)\ar[d]^{\partial_{l}}&\\
BS^{3}\ar[r]^{c_{2}}&K(\Z,4)\ar[r]^{l}&K(\Z,4)}\ ,\end{equation}
then we get connecting maps $X_{l}\to X_{l^{\prime}}$ if $l|l^{\prime}$ and 
$$X\cong \hocolim_{l}\: X_{l}\ ,\quad  \pi_{n}(MB )=\colim_{l}\: \pi_{n}(MB_{l})\ ,$$ where the Thom spectrum $MB_{l}$ is associated to the projection 
$B_{l}:=BSpin^{c}\times X_{l}\to BSpin^{c}$.   The main advantage of $X_{l}$ is that it has finite skeleta. 

We consider a cycle $(M,f)$ for $\pi_{n}(MSpin^{c})$
 and an auxiliary map $g:M\to X$. We can assume that $g$ has a factorisation
$$g:M\stackrel{g_{l}}{\to} X_{l}\to X$$ for some $l$. We choose  a $Spin^{c}$-extension $\tilde \nabla^{TM}$ of the Levi-Civita connection. 
We are going to construct an $\ell$-good geometrisation for $(M,(f,g_{l}),\tilde \nabla^{TM})$ using similar ideas as in the $String$-bordism case Proposition \ref{prop77}.
We first take $m$ sufficiently large such that the canonical map $\iota_{m}:BSpin^{c}(m)\to BSpin^{c}$ is $\max\{\ell,n+1,4\}$-connected.
We choose a compact $ \max(\ell,n+1,4)$-connected approximation $(f_{u},g_{u}):M_{u}\to BSpin^{c}(m)\times X_{l}$
such that the map $(f,g_{l})$ factorizes over a closed embedding
$h:M\to M_{u}$. It is here where we use the property that $X_{l}$ has finite skeleta  which ensures that we can find a compact approximation $M_{u}$.

We choose the $Spin^{c}(m)$-connection $\tilde \nabla^{u}$ as in  Subsection \ref{sec801}. The map $h$ has a refinement to a $Spin^{c}$-map 
and we can assume that $h^{*}\tilde \nabla^{u}=\tilde \nabla^{TM}$ stably. 
 
The composition $q_{l}\circ g_{u}:M_{u}\to BS^{3}$ classifies an $S^{3}$-principal bundle $R_{u}\to M_{u}$ on which we choose a connection $\nabla^{R_{u}}$. We can assume that  $R\cong h^{*}R_{u}$ with connection $\nabla^{R}=h^{*}\nabla^{R_{u}}$.

We let $$\tilde \cG_{u}:K^{0}(BSpin^{c}(m)\times  BS^{3})\to \hat K^{0}(M_{u})$$ denote the
geometrisation of $(M_{u},(f_{u},q_{l}\circ g),\tilde \nabla^{u})$ which was constructed by the method of   Lemma \ref{cgl1} from a geometrization $\cG^{0}:=\cG^{Spin^{c}(m)}$ which is the analog of $\cG^{Spin(m)}$ in the proof of Proposition \ref{lem77}.

 We have a fibration
\begin{equation}\label{cg15551}BSpin^{c}(m)\times  X_{l}\stackrel{(\id,q_{l})}{\to}BSpin^{c}(m)\times  BS^{3}\end{equation}
with fibre $K(\Z,3)$.
\textcolor{black}{We can again apply Proposition \ref{pinv1} 
 and conclude  that
\begin{equation}\label{cg1555}(\id,q_{l})^{*}:\bar K^{*}(BSpin^{c}(m)\times  BS^{3})\to \bar K^{*}(BSpin^{c}(m)\times X_{l})\end{equation}  represents $\bar K^{*}(BSpin^{c}(m)\times X_{l})$ as a completion of $\bar K^{*}(BSpin^{c}(m)\times  BS^{3})$.}
We define 
$$\cG_{u,0}:\textcolor{black}{\bar K^{*}(BSpin^{c}(m)\times  BS^{3})   }   \to \hat K^{0}(M_{u})$$
by
$$\cG_{u,0}( \phi):=\tilde \cG_{u}(\phi)\in \hat K^{0}(M_{u})\ , \quad \phi\in \bar K^{*}(BSpin^{c}(m)\times  BS^{3})\ .$$
\textcolor{black}{This  map  needs a correction in order to be continuous with respect to the topology on
$\bar K^{*}(BSpin^{c}(m)\times  BS^{3})$ induced from the profinite topology of $\bar K^{*}(BSpin^{c}(m)\times X_{l})$.}
We must kill the contribution of $c_{2}(\nabla^{R_{u}})$ to the curvature of $\cG_{u,0}( \phi)$.
Note that $q_{l}^{*}c_{2}\in H^{4}(X_{l};\Z)$ is $l$-torsion.  Hence
we can choose a form $\alpha_{u}\in \Omega^{3}(M_{u})/\im(d)$ such that $d\alpha_{u}=c_{2}(\nabla^{R_{u}})$.
By an easy application of Serre's spectral sequence to the fibration (\ref{cg15551}) we see that
$$p^{*}:H^{*}(BSpin^{c}(m);\Q)\to H^{*}(BSpin^{c}(m)\times X_{l};\Q)$$ is an isomorphism.
Since $H^{*}(BSpin^{c}(m);\Q)$ is concentrated in even degrees
the odd-dimensional cohomology of $M_{u}$ is concentrated in degrees $\ge n+1$.
In particular we see that $\alpha_{u}$ is uniquely determined.
Moreover, the restriction 
\begin{equation}\label{okt2801-2012}h^{*}:HP\R^{-1}(M_{u})\to HP\R^{-1}(M)\end{equation} is trivial.

We have
\begin{eqnarray*}HP\Q^{*}(BSpin^{c}(m)\times BS^{3})&\cong& \Q[b,b^{-1}][[bc_{1},b^{2}p_{1},b^{4}p_{2}\dots,b^{2r_{m}}p_{r_{m}},b^{2}c_{2}]]\ ,\\ HP\Q^{*}(BSpin^{c}(m)\times X_{l})&\cong& \Q[b,b^{-1}][[bc_{1},b^{2}p_{1},b^{4}p_{2}\dots,b^{2r_{m}}p_{r_{m}}]]\ ,\end{eqnarray*}
where $c_{1}$ and the Pontrjagin classes come from $BSpin^{c}$, and $c_{2}$ is pulled back from $BS^{3}$.
The pull-back $(\id,q_{l})^{*}$ is the quotient map defined by setting $c_{2}=0$.
For $\phi\in K^{0}(BSpin^{c}(m)\times  BS^{3})$ we define the formal power series
$$\Phi_{\phi}:=\Td^{-1}\cup \ch(\phi)\in \Q[b,b^{-1}][[bc_{1},b^{2}p_{1},b^{4}p_{2}\dots,b^{2}c_{2}]]^{0}$$
and set
$$\tilde \Phi_{\phi}:=\frac{\Phi_{\phi}-i_{c_{2}=0}\Phi_{\phi}}{c_{2}}\in \Q[b,b^{-1}][[bc_{1},b^{2}p_{1},b^{4}p_{2}\dots,b^{2}c_{2}]]^{-4}\ .$$
For $\phi\in  \bar  K^{0}(BSpin^{c}(m)\times BS^{3})$ now define
$$\cG_{u}( \phi) :=\cG_{u,0}( \phi )-a(\alpha_{u}\wedge \Td(\tilde \nabla^{u})^{-1}\wedge \tilde \Phi_{\phi}(\nabla^{u},\nabla^{R_{u}}))\ ,$$
where $\tilde \Phi_{\phi}(\nabla^{u},\nabla^{R_{u}})\in \Omega P^{-4}(M)$ is obtained from $\tilde \Phi_{\phi}$ by replacing
the generators $c_{1},p_{i}$ and $c_{2}$ by their corresponding Chern-Weyl representatives
$c_{1}(\tilde \nabla^{u})$, $p_{i}(\tilde \nabla^{u})$, and $c_{2}(\nabla^{R_{u}})$.
We calculate similarly as in (\ref{eqcont30}) that
\begin{equation}\label{cg20} \Td(\tilde \nabla^{u})\wedge R(\cG_{u}(\phi))=i_{c_{2}=0}\Phi_{\phi}(\nabla^{u},\nabla^{R_{u}})\ .\end{equation}

We now define $$\cG: \textcolor{black}{\bar K^{*}(BSpin^{c}\times  BS^{3})   }\to \hat K^{0}(M)$$ by 
$$\cG( \phi ):=(\iota_{m}\times \id_{X_{l}})_{*}h^{*}\cG_{u}( \phi  )\ .$$
We claim that $\cG$ extends by continuity to a good geometrisation of $(M,(f,g_{l}),\tilde \nabla^{TM})$. 
The argument is very similar to that of Lemma \ref{cont1}.  We first show continuity. If $(\phi_{k})$ is a sequence in
$\bar K^{0}(BSpin^{c}\times BS^{3})$ with $(\id,q_{l})^{*}\phi_{k}\to 0$ as $k\to\infty$ in the profinite topology, 
then we can find a $k_{0}\in \nat $ such that
$$\cG_{u}(\phi_{k})\in HP\R^{-1}(M_{u})/\im(\ch)$$
for all $k\ge k_{0}$.
It follows from the vanishing of the map $h^{*}$ in  \eqref{okt2801-2012} that $h^{*}\cG_{u}(\phi_{k})=0$.

 Because of (\ref{cg20}) \textcolor{black}{a degree-preserving}
cohomological character of $\cG$ is given by 
$$bc_{1}\mapsto bc_{1}(\tilde \nabla^{TM})\ ,\quad  b^{2i}p_{i}\mapsto b^{2i}p_{i}(\tilde \nabla^{TM})\ .$$

It follows that $\cG$ is a geometrisation. In order to see that it is $\ell$-good we show that
$\cG_{u}$ itself is continuous using a similar argument based on a $\dim(M_{u})+1$-connected approximation of $BSpin^{c}(m)\times X_{l}$.

We can now apply  Theorem \ref{them2} in order to derive a formula for $\varepsilon([M,(f,g_{l})])\in \R/\Z$.
We can take $\hat c_{2}(\nabla^{R}):=h^{*}\alpha_{u}$ and have
$\tilde \Phi_{A}=\Td^{-1} \tilde \Phi$.
We have by construction
$$\cG((\id,q_{l})^{*}A)=[2-\bE_{\rho}]-a(\hat c_{2}(\nabla^{R})\wedge  \tilde \Phi(\nabla^{R}) )\ ,$$
hence the correction form (Definition  \ref{cgkl1}) is given by
$$\gamma_{(\id,q_{l})^{*}A}=-\hat c_{2}\wedge   \tilde \Phi(\nabla^{R}) \ .$$
It follows from \eqref{eq1000} that
$$\ev_{(\id,q_{l})^{*}A}(\eta^{an}([M,(f,g_{l})]))=[\int_{M}\Td(\tilde \nabla^{TM})\wedge \hat c_{2}(\nabla^{R})\wedge  \tilde \Phi(\nabla^{R})] -2\xi(\Dirac_{M})+\xi(\Dirac_{M}\otimes \bE)\ .$$
This is exactly the formula \eqref{cg2} for $t^{\C}_{M}(\tilde g)\in \R/\Z$. 
The Proposition \ref{cg78} now follows from Lemma \ref{lem1000} which gives the first equality in the chain
$$\ev_{(\id,q_{l})^{*}A}(\eta^{an}([M,(f,g_{l})]))=\ev_{(\id,q)^{*}A}(\eta^{an}([M,(f,g)]))=\varepsilon(s_{M}(\tilde g
))\ .$$
\hB 

The paper \cite{2010arXiv1012.5237C} provides a lot of interesting explicit calculations. 
Our general point of view is probably not of much help here. But it is useful to understand structural results like the relation with the Adams $e$-invariant \cite[Prop 1.11]{2010arXiv1012.5237C}. This is what we are going to explain now.
 We define the space $Y$ by extending the diagram  (\ref{cg1}) by another cartesian square as follows
\begin{equation}
\label{cg1v}
\xymatrix{Y\ar[r]^{H}\ar[d]^{r}&X\ar[d]^{q}\ar[r]&K(\Q/\Z,3)\ar[d]^{\partial} \\ 
S^{4}\ar[r]^{h}&BS^{3}\ar[r]^{c_{2}}&K(\Z,4)}   \  ,
\end{equation}
where $h$ generates $\pi_{4}(BS^{3})$
such that $h^{*}c_{2}\in H^{4}(S^{4};\Z)$ is the positive orientation class. 
We use the Serre spectral sequence in order to calculate the rational cohomology of $Y$:
$$H^{k}(Y;\Q)=\left\{
\begin{array}{cc}
\Q&k=0,7
\\0&\: k\not\in \{0,7\}\end{array}\right. \ .$$
This implies
\begin{equation}\label{cgh2}\pi_{4m-1}(S\wedge Y_{+})_{tors}=\pi_{4m-1}(S\wedge Y_{+})\end{equation}
for $m\ge 3$. 

From now on we assume that $m\ge 2$.
By Lemma \ref{lem1000}  we get the commutativity of the  squares (except of the lower right which will be explained below) of  the following  diagram:
\begin{equation}\label{cg100}
\xymatrix{\pi_{4m-1}(MSpin^{c}\wedge X_{+})\ar@/^1cm/@{.>}[rr]^{\varepsilon}\ar[r]^{\eta^{top}}&Q_{4m-1}(MSpin^{c}\wedge  X_{+})\ar[r]^{\ev_{p^{*}q^{*}A}}&\Q/\Z\\
\pi_{4m-1}(S\wedge Y_{+})_{tors}\ar[d]_{(\id\wedge r)_{*}}\ar[u]^{(\epsilon_{MSpin^{c}}\wedge H)_{*}}\ar[r]^{\eta^{top}}&Q_{4m-1}(S\wedge Y_{+})\ar[d]_{(\id, r)}\ar[u]^{(i_{*},H)}\ar[r]^{ev_{r^{*}h_{+}^{*}A}}&\Q/\Z\ar@{=}[u]\\
\pi_{4m-1}(S\wedge S^{4}_{+})\ar[d]_{w_{*}}\ar[r]^{\eta^{top}}&Q_{4m-1}(S\wedge S^{4}_{+})\ar[d]_{w_{*}}\ar[r]^{\ev_{h_{+}^{*}A}}&\Q/\Z\ar@{=}[u]\\
\pi_{4m-1}(S\wedge S^{4})\ar[d]_{\cong}\ar[r]^{\eta^{top}}&Q_{4m-1}( S\wedge S^{4})\ar[d]_{\cong}\ar[r]^{\ev_{h^{*}A}}&\Q/\Z\ar@{=}[u]\\
\pi_{4m-5}(S)\ar@/_1cm/@{.>}[rr]_{e^{Adams}_{\C}}\ar[r]^{\eta^{top}}&Q_{4m-5}(S)\ar[r]^{\ev_{1}}&\Q/\Z\ar@{=}[u]}\ .
\end{equation}
We need the condition $m\ge 2$ in order to have well-defined evaluations $\ev_{h_{+}^{*}A}$, $\ev_{h^{*}A}$ and $\ev_{1}$.
The map $w_{*}$ is induced by the map $w:S^{4}_{+}\to S^{4}$ which is the identity on $S^{4}$ and maps the extra base point to the base point of $S^{4}$. In terms of the canonical decomposition
$$K^{0}(S^{4}_{+})\cong  K^{0}(S^{4})\oplus \Z$$ we have $w^{*}=(\id_{K^{0}(S^{4})},\dim)$.
 This map induces 
$$w_{*}:Q_{4m-1}(S\wedge S^{4}_{+})\to Q_{4m-1}(S\wedge S^{4})$$
   We use the symbol $h_{+}:S^{4}_{+}\to BS^{3}$ for the map induced by $h$ which maps the extra  base point to a base point of $BS^{3}$. 
The lower left vertical map is the suspension isomorphism. The lower middle
vertical isomorphism is again induced by suspension and the Bott isomorphism
$$K^{0}(S^{4})\cong K^{-4}(S^{0})\stackrel{b^{-2}}{\to} K^{0}(S^{0})\ .$$
In order to see that the lower right square commutes note that this isomorphism maps
$h^{*}A$ to $1$. This follows from
$$\ch(2-A)=bc_{2}+O(b^{3})$$
and the fact that $c_{2}\in H^{4}(S^{4};\Z)$ is the suspension of $1\in H^{0}(*;\Z)$.  
The composition of the lower two arrows is the definition (\ref{eq1200}) of the complex version of the Adams $e$-invariant.

We conclude that
\begin{equation}\label{cg56}\varepsilon\circ (\epsilon_{MSpin^{c}}\wedge H)=e^{Adams}_{\C}\circ w_{*}\ .\end{equation}
The same argument as for Lemma \ref{lemcg5} gives
\begin{lem}\label{lemblav}
A pair $((M,f),\tilde g)$ of  a cycle for $\pi_{4m-1}(S)$ such that 
$H^{3}(M;\Q)=0$ and $H^{4}(M;\Q)=0$,  and a map $\tilde g\in [M,S^{4}]$ gives naturally rise to a class
$[M,(f,g)]\in \pi_{4m-1}(S\wedge Y_{+})$. 
\end{lem}
If $M$ satisfies the assumption of the Lemma, then we have a preconstruction map
$$\tilde s_{M}:[M,S^{4}]\to \pi_{4m-1}(S\wedge Y_{+})$$
and conclude from Proposition \ref{cg78}, (\ref{cg56}) and (\ref{cgh2}) that for $m\ge 3$ (or $m=2$ and $[M,(f,g)]\in \pi_{7}(S\wedge Y_{+})$ is a torsion class) 
$$t_{M}^{\C}=e^{Adams}_{\C}\circ w_{*}\circ \tilde s_{M}:[M,S^{4}]\to \Q/\Z\ .$$
This is \cite[Prop 1.11]{2010arXiv1012.5237C} if one takes the following
geometric description of the composition
$w_{*}\circ \tilde s_{M}(\tilde g)$ into account. First of all we have $\tilde s_{M}(\tilde g)=[M,(f,g)]$,
where $g:M\to Y$ is the lift of $\tilde g:M\to S^{4}$.
Then $w_{*}([M,(f,g)])=[M,(f,\tilde g)]-[M,const]\in \pi_{4m-1}(S\wedge S^{4})$. The geometric representative of the  $4$-fold desuspension
of this class is the stably normally framed manifold obtained by taking the  preimage $Y:=\tilde g^{-1}(s)$
of a regular point $s\in S^{4}$ of $\tilde g$.
\begin{kor}\label{intr2z} \cite[Prop 1.11]{2010arXiv1012.5237C}\footnote{The $e$-invariant used in the present paper is the negative of the $e$-invariant in the conventions of \cite[Prop 1.11]{2010arXiv1012.5237C}. This accounts for the different sign.}
We assume that $m\ge 2$. Let  $(M,(f,\tilde g))$  be as in Lemma \ref{lemblav}. 
If $m=2$, then in addition we assume that $[M,(f,g)]\in \pi_{7}(S\wedge Y_{+})$ is a torsion class.
Then we have
$$t^{\C}_{M}(h\circ \tilde g)=e^{Adams}_{\C}([Y,f^{\prime}])\ ,$$
where $Y$ is the  preimage $Y:=\tilde g^{-1}(s)$
of a regular point $s\in S^{4}$ of $\tilde g$ with its induced normal framing (and $f^{\prime}$ is the constant map).
\end{kor}

\begin{rem}\label{hiurfhrifhfx}{\rm In the following we consider the case $m=2$ and discuss what happens if we  drop the assumption that $[M,(f,g)]\in \pi_{4m-1}(S\wedge Y_{+})$ is a torsion element.
We consider the Hopf fibration $\tilde g:S^{7}\to S^{4}$. By Lemma \ref{lemblav} we get an element
$[S^{7},(f,g)]\in \pi_{7}(S\wedge Y_{+})$.  This element is not torsion.

 If it would be a torsion element,
then by Corollary \ref{intr2z} we would have $t^{\C}_{S^{7}}(h\circ \tilde g)=e_{\C}^{Adams}([Y,f^{\prime}])$,
where $(Y,f^{\prime})$ is a Hopf fibre with the induced framing.
It has been shown in \cite[Example 3.5]{2010arXiv1012.5237C} that
$t^{\C}_{S^{7}}(h\circ \tilde g)=0$. On the other hand, since
the Hopf fibration generates the stable homotopy group $\pi_{3}(S)\cong \Z/24\Z$ which is detected completely by $e^{Adams}$ we know that $e_{\C}^{Adams}([Y,f^{\prime}])\in \Q/\Z$
has order $12$, in particular is non-trivial.

Now $$\ev_{h_{+}^{*}A}(\eta^{top}((\id\wedge r)_{*}[S^{7},(f,g)]))  = e_{\C}^{Adams}([Y,f^{\prime}])\not=0$$  is a non-trivial torsion class of order $12$. 
On the other hand $$\ev_{p^{*}q^{*}A}(\eta^{top}((\epsilon_{MSpin^{c}}\wedge H)_{*}[S^{7},(f,g)]))=t^{\C}_{S^{7}}(h\circ \tilde g)=0\ .$$
We see that the upper half of the diagram \eqref{cg100} no longer commutes if we delete the subscript
 $(\dots)_{tors}$ in the second line.}\end{rem}

\section{The $K$-theory of $K(\pi,n)$-bundles}\label{hdklqwjdlqwdwqd}

\color{black}

The goal of this Section is to give a 
 proof of   Proposition \ref{pinv}. It uses some theory which we develop in greater generality for the purpose of applications at other places. The main result is Proposition \ref{pinv1}.  
 \bigskip

Let $Z$ be a space with an increasing filtration  $$\dots\subseteq Z_{k}\subseteq Z_{k+1}\subseteq\dots\ , \quad k\in \nat$$ such that $Z \simeq \hocolim_{k\in \nat}Z_{k}$. Then we consider the decreasing filtration
$$F^{k}K^{*}(Z):=\ker(K^{*}(Z)\to K^{*}(Z_{k-1}))$$
of the $K$-theory group $K^{*}(Z)$, and we write $\Gr^{k}_{F}(K^{*}(Z))$ for its associated 
subquotients. The filtration $(F^{k}K^{*}(Z))_{k\in \nat}$ induces a topology on $K^{*}(Z)$, and we write
${}^{F}\bar K^{*}(Z)$ for the Hausdorff completion of $K^{*}(Z)$ with respect to this topology. 

If $Z^{\prime}$ is a second space with  increasing filtration $(Z^{\prime}_{k})_{k\in \nat}$ and $Z\to Z^{\prime}$ is a filtration preserving map, then the pull-back
$ K^{*}(Z^{\prime})\to K^{*}(Z)$ is fitration preserving, continuous and induces a continuous map
${}^{F}\bar K^{*}(Z^{\prime})\to {}^{F}\bar K^{*}(Z)$.

  If the $Z_{k}$, $k\in \nat$, are finite $CW$-complexes, then the induced topology on $K^{*}(Z)$  is   the profinite topology.  
For a general filtration of $Z$ the profinite topology is always contained in the topology
associated to the filtration and we have a continuous map ${}^{F}\bar K^{*}(Z)\to \bar K^{*}(Z)$.



The filtration of $Z$ induces a spectral sequence
$(E_{r},d_{r})$ with $E_{1}^{s,t}\cong  K^{t+s}(Z_{s}/Z_{s-1})$.
If the filtration of $Z$ is the skeletal filtration, then this is the Atiyah-Hirzebruch spectral sequence whose second page is given by $E_{2}^{s,t}\cong H^{s}(Z,\pi_{-t}(K))$.

If we fix $(s,t) \in \nat\times \Z$, then for $r\ge s$ the group $E_{r}^{s,t}$ does not receive differentials. Hence we can consinder the sequence of inclusions of groups $$E_{\infty}^{s,t}:=\bigcap_{r\ge s} E_{r}^{s,t}\subseteq \dots\subseteq E_{r+2}^{s,t}\subseteq E_{r+1}^{s,t}\subseteq E_{r}^{s,t}\ .$$
For every $k\in \nat$ there is a natural map
\begin{equation}\label{mirette7}\sigma_{k}:\Gr_{F}^{k}K^{t+k}(Z)\to E_{\infty}^{k,t}\ . \end{equation}
We say that the spectral sequence converges strongly if
\eqref{mirette7} is an isomorphism for all $k\in \nat$ and $t\in \Z$. 
This is the case for example if the spectral sequence degenerates at a finite stage.

\begin{prop}\label{pinv1} Let $X$ be a $CW$-complex with an increasing filtration $(X_{k})_{k\in \nat}$ by finite subcomplexes $X_{k}$ such that $X=\bigcup_{k}X_{k}$.  We further fix an integer $n\ge 3$, a countable abelian group $\pi$, and consider a fibration
$p:Y\to X$ with fibre $K(\pi,n)$. Let $(Y_{k})_{k\in \nat}$ be the induced filtration of $Y$.  
\begin{enumerate}
\item
The projection
$p:Y\to X$ induces a continuous  injective   map \begin{equation}\label{barp111} \bar K^{*}(X) \to {}^{F} \bar K^{*}(Y)\ .\end{equation}
\item If   $\lim^{1}_{k}K^{*}(X_{k})=0$,     
then the image  of the composition \begin{equation}\label{barp1111}\bar K^{*}(X)\to {}^{F}\bar K^{*}(Y)\to \bar K^{*}(Y)\end{equation}  is dense.
\end{enumerate}
 \end{prop}
\begin{rem}{\rm If $\pi$ is a countable torsion group, then \cite[Prop. 4.7]{MR0231369} gives the stronger statement that
$p^{*}:K^{*}(X)\to K^{*}(Y)$ is an isomorphism. In this case we can even assume that $n\ge 2$, and we can drop the assumption that $X$ has finite skeleta.} \hB
\end{rem}

\proof[of Proposition \ref{pinv} assuming Proposition \ref{pinv1}]
Proposition \ref{pinv} follows from the combination of the two assertions of Proposition \ref{pinv1}
if we take $n:=3$, $\pi:=\Z$, and $p:BString(m)\to BSpin(m)$ (note that $m\ge 3$).  The homotopy type $BSpin(m)$ admits a cell structure with finite skeleta. Moreover, it has been shown in \cite[Sec. 2]{MR0259946} that
$\lim^{1}_{k} K^{*}(BSpin(m)_{k})=0$.
\hB

In order to show Proposition \ref{pinv1}
we first
need some preparations about divisible groups.
Let $A$ be some abelian group. Then we define its subgroup 
$$A_{div}:=\{a\in A\:|\: \forall n\in \nat\: \exists a^{\prime}\in A \:\:\mbox{s.t.}\:\: a=na^{\prime}\}$$
of divisible elements.
We consider the exact sequence
$$0\to A_{div}\to A\to \bar A\to 0\ .$$
Since a divisible group is injective this sequence is split. Hence we have a non-canonical decomposition
$$A\cong A_{div}\oplus \bar A\ .$$ This implies that
 $\bar A_{div}=0$.
 We now consider a short exact sequence of groups
 $$0\to A\to B\to C\to 0$$
 together with a homomorphism $B\to X$, where $X$ is a finitely generated abelian group.
 \begin{lem}\label{lemhilfe}
If $c\in C_{div}$,
then we can find a lift $b\in B$ of $c$ whose image in $X$ vanishes.
\end{lem}
\proof
We consider the diagram
 $$\xymatrix{0\ar[d]&0\ar[d]&0\ar[d]\\  A_{div}\ar[d]\ar[r]&B_{div}\ar[r]\ar[d]&C_{div} \ar[d] \\ A\ar[r]\ar[d]&B\ar[r]\ar[d]&C\ar[d]\\ \bar A\ar[r]\ar[d]&\bar B\ar[r]\ar[d]&\bar C\ar[d]\\0&0&0}$$
with vertical exact sequences. The Snake Lemma gives an isomorphism
$$\frac{C_{div}}{\im(B_{div}\to C_{div})}\cong \frac{\ker(\bar B\to \bar C)}{\im(\bar A\to \bar B)}\ . $$
This shows that the group on the right-hand side is divisible.
 Since any map from a divisible group to a finitely generated group is trivial we have a factorisation
$\bar B\to X$ of the map $B\to X$. We thus get a homomorphism
$$ \frac{\ker(\bar B\to \bar C)}{\im(\bar A\to \bar B)}\to \frac{X}{Y}\ ,$$
where $Y\subseteq X$ is the image of
$\bar A\to \bar B\to X$.  The quotient $X/Y$ is still finitely generated, and this implies that this map is trivial since its domain is divisible.

 We now  choose a preimage $b_{0}\in B$ of $c$.
Its image in $\bar b_{0}\in \bar B$ then vanishes when mapped to $\bar C$ so that it represents a class in $\frac{\ker(\bar B\to \bar C)}{\im(\bar A\to \bar B)}$.
The image of this class in $X/Y$ vanishes so that there exists $\bar a\in \bar A$ such that the image of $\bar b_{0}-\bar a$ in $X$ vanishes. We choose some lift $a\in A$ of $\bar a$. Then the image of $b:=b_{0}-a$ in $X$ vanishes. Moreover, $b$ can be taken as a lift of $c$, too. 
\hB

\begin{prop}\label{mirette4}
Assume that $X$ is a finite $CW$-complex and $p:Y\to X$ is a fibration with fibre
$K(\pi,n)$ for $n\ge 3$ and $\pi$ a countable abelian group.
Then the map \begin{equation}\label{mirette1}p^{*}:  K(X)\to   K^{*}(Y)\end{equation}  is injective and has dense range.
\end{prop}
\proof
Let $$\emptyset=X_{-1}\subseteq X_{0}\subseteq X_{1}\subseteq \dots\subseteq X_{\dim(X)}=X$$
be the filtration of $X$ given by the cellular structure. It induces a filtration
$$\emptyset=Y_{-1}\subseteq Y_{0}\subseteq Y_{1}\subseteq \dots\subseteq Y_{\dim(X)}=Y$$
by taking preimages. These filtrations induce spectral sequences
$E(\id)$ and $E(p)$ which both degenerate at the $\dim(X)$'th page and strongly converge. Here $E(\id)$ is the Atiyah-Hirzebruch spectral sequence for $X$, and $E(p)$ is the Serre spectral sequence
for the fibration $p$.

 We consider the map of fibrations
$$\xymatrix{Y\ar[r]\ar[d]^{p}&X\ar[d]^{\id}\\X\ar[r]&X}$$
which induces a morphism of  spectral sequences 
$$p^{*}:E(\id)\to E(p)\ .$$
\begin{lem}\label{mirette40}
The image $p^{*}E(\id)$ is a direct summand of $E(p)$. 
\end{lem}
\proof
We have  
 $$E^{s,t}_{1}(p)\cong C^{s}(X,K^{t-s}(K(\pi,n)))\ ,$$ 
 where $C^{s}(X,A)$ denotes the cellular cochains of $X$ with coefficients in the abelian group $A$.
  We have an exact sequence
$$0\to E_{1}^{s,t}(p)_{div}\to E_{1}^{s,t}(p)\to \overline E_{1}^{s,t}(p) \to 0\ .$$
Since $\tilde K^{*}(K(\pi,n))$ is uniquely divisible by
\cite[Thm. II]{MR0231369} the composition
$$E_{1}^{s,t}(\id)\stackrel{p^{*}}{\to}E_{1}^{s,t} (p)\to \overline E_{1}^{s,t}(p)$$ is an isomorphism. 
  Moreover, since $p^{*}$ is a chain map and there are no non-trivial maps from a divisible group to a finitely generated group we have a decomposition of chain complexes
$$E_{1}^{*,*}(p)\cong E_{1}^{*,*}(p)_{div}\oplus p^{*}E_{1}^{*,*}(\id)\ .$$
 
 Using that $E_{1}^{*,*}(p)_{div}$ is actually uniquely divisible and that a subquotient of a uniquely divisible group is again uniquely divisible we obtain an induced  decomposition of the
  higher pages and  a decomposition of the whole spectral sequence as
$$E(p)\cong E(p)_{div}\oplus p^{*}E(\id)\ .$$
\hB

%
%
%

\begin{lem}\label{mirette5}
The map \eqref{mirette1} is injective.
\end{lem}
\proof
Assume that $\phi\in  K^{*}(X)$ is such that $p^{*}(\phi)=0$. We are going to show that 
 $\phi\in F^{k}K^{*}(X) $ for all $k\ge -1$. For $k\ge \dim(X)+1$ this then implies that $\phi=0$.
 
We clearly have $\phi\in F^{-1}K^{*}(X)=K^{*}(X)$.  Assume now by induction that $\phi\in  F^{k}K^{*}(X)$.
Then $\phi$ is detected by an element $\sigma_{k}(\phi)\in E^{k,*}_{\infty}(\id)$.
We observe that  the image of $\sigma_{k}(\phi)$ in $E^{k,*}_{\infty}(p)$ 
 and in particular its component in $p^{*}E^{k,*}_{\infty}(\id)$ vanishes.  By Lemma \ref{mirette40}  we have  $\sigma_{k}(\phi)=0$. This implies $\phi\in F^{k+1}K^{*}(X) $.
\hB

\begin{lem}\label{mirette2}
The range of the map  \eqref{mirette1} is dense.
\end{lem}
\proof
We fix an element $\phi\in K^{*}(Y)$.
We must show that we can approximate this element in the profinite topology of $K^{*}(Y)$ by elements in the image of $p^{*}:K^{*}(X)\to K^{*}(Y)$.
 Let $t:T\to Y$ be a map from a finite $CW$-complex. Then $\ker(t^{*})\subseteq K^{*}(Y)$ is some neighbourhood of zero. We must  find an element
 $\psi\in K^{*}(X)$ such that 
 \begin{equation}\label{thegoal}\phi-p^{*}\psi\in \ker(t^{*})\ .\end{equation}
To this end we use the following Lemma
\begin{lem}\label{mirette3}
If $\hat \phi\in F^{k}K^{*}(Y)$, then there exists $\hat \psi\in K^{*}(X)$ and $\rho\in \ker(t^{*})$ such that
$\hat \phi-p^{*}\hat \psi-\rho\in F^{k+1}K^{*}(Y)$.
\end{lem}

Assuming this Lemma and using that $\phi\in F^{-1}K^{*}(Y)$ we obtain
the desired $\psi\in K^{*}(X)$ by a finite iteration. Thus Lemma \ref{mirette2}
follows from Lemma \ref{mirette3}. \hB 

\proof[of Lemma \ref{mirette3}]
The element $\hat\phi$ gives rise to an element 
$u:=\sigma_{k}(\hat \phi)\in E_{\infty}^{k,*}(p)$ which we can decompose as $u=v\oplus p^{*}w$ with $v\in E_{\infty}^{k,*}(p)_{div}$ and $w\in E^{k,*}_{\infty}(\id)$.
We let $\hat \psi\in F^{k}K^{*}(X)$ be an element which is detected by
 $w$.  
We apply Lemma \ref{lemhilfe} to the exact sequence
$$0\to F^{k+1}K^{*}(Y)\to F^{k}K^{*}(Y)\to E_{\infty}^{k,*}(p)\to 0$$
and the map $ t^{*}:F^{k}K^{*}(Y)\to K^{*}(T)$.
 By Lemma  \ref{lemhilfe}
we can find an element $\rho\in F^{k}K^{*}(Y)\cap \ker(t^{*})$ which represents
$v$.
Then we have
$\hat \phi-p^{*}\hat \psi-\rho\in F^{k+1}K^{*}(Y)$ as required.
\hB 

Proposition \ref{mirette4} now follows from Lemma \ref{mirette5} and Lemma \ref{mirette2}. \hB  

\proof[of Proposition \ref{pinv1}]

Let
$$\emptyset=Y_{-1}\subseteq Y_{0}\subseteq Y_{1}\subseteq Y_{2}\subseteq \dots$$ be the filtration of $Y$ induced by taking the preimages of the subcomplexes $X_{k}$ along $p:Y\to X$.
By construction the map $p$ preserves filtrations  and hence \eqref{barp111} is continuous.
%


\bigskip

We first show that \eqref{barp111} is injective. Note that the filtration 
$(F^{k}K^{*}(X))_{k\ge 0}$ induces the profinite topology.
Let $\phi\in K^{*}(X)$ be such that $p^{*}\phi\in \bigcap_{k\ge 0} F^{k}K^{*}(Y)$.
It suffices to show that this implies $\phi\in F^{k}K^{*}(X)$ for all $k\ge 0$ since then
$\phi$ represents the zero element in $\bar K^{*}(X)$.

We clearly have $\phi\in F^{-1}K^{*}(X)$. We now assume by induction that $\phi\in F^{k}K^{*}(X)$. We can apply Proposition \ref{mirette4} to the fibration $p_{k}:Y_{k} \to X_{k}$.
Since $$p_{k}^{*}(\phi_{|X_{k}})=(p^{*}\phi)_{|Y_{k}}=0$$ we conclude that $\phi_{|X_{k}}=0$ and hence $\phi\in F^{k+1}K^{*}(X)$.

\bigskip

We now show that the image of \eqref{barp1111} is dense.
Let $\phi\in K^{*}(Y)$. We must approximate $\phi$ by elements in the image of \eqref{barp1111}.
 Let $t:T\to Y$ be a map from a finite $CW$-complex. Then $\ker(t^{*})\subseteq K^{*}(Y)$ is some neighbourhood of zero. We must   find an element
 $\psi\in K^{*}(X)$ such that 
 \begin{equation}\label{thegoal1111}\phi-p^{*}\psi\in \ker(t^{*})\ .\end{equation}
  Since $T$ is finite there exists a number  $k \in \nat $ such that there is a factorisation $$t:T\stackrel{  t_{k}}{\to} Y_{k }\stackrel{\iota_{k}}{\to} Y\ .$$ 
 It suffices to find  an element $\psi\in K^{*}(X)$ such that
  \begin{equation}\label{thegoal111}\iota_{k}^{*}\phi-\iota_{k}^{*}p^{*}\psi\in \ker(   t_{k}^{*})\ .\end{equation}
  Then indeed $\phi-p^{*}\psi\in \ker(t^{*})$.

  For every $r\in \nat$  we apply  Proposition \ref{mirette4} to the fibration $p_{k+r }:Y_{k+r }\to X_{k+r }$
  in order to see that $p_{k+r}^{*}$ induces an isomorphism
$$\frac{K^{*}(X_{k+r})}{\ker((p_{k+r}\circ t_{k+r})^{*})}\stackrel{\sim}{\to}  \frac{K^{*}(Y_{k+r})}{\ker(t^{*}_{k+r})}\ .$$

We have a projective system of exact sequences indexed by $r\in \nat$
$$0\to \ker((p_{k+r}\circ t_{k+r})^{*})\to K^{*}(X_{k+r}) \to  
\frac{K^{*}(Y_{k+r})}{\ker(t^{*}_{k+r})}\to 0\ .$$ If we take the limit, 
identify
$$\frac{\lim_{r} K^{*}(X_{k+r})}{\lim_{r} \ker((p_{k+r}\circ t_{k+r})^{*})}\cong  \frac{\bar K^{*}(X)}{\ker((p\circ t)^{*})}\ ,$$
and use the assumption that $\lim^{1}_{r} K^{*}(X_{k+r})=0$, then
we get the exact sequence
$$0\to  \frac{\bar K^{*}(X)}{\ker((p\circ t)^{*})}  \to \lim_{r}\frac{K^{*}(Y_{k+r})}{\ker(t^{*}_{k+r})}
\to \lim_{r}^{1}  \ker((p_{k+r}\circ t_{k+r})^{*}) \to 0\ .$$
Note note that the graded groups
$$\lim_{r}\frac{K^{*}(Y_{k+r})}{\ker(t^{*}_{k+r})}\cong \bigcap_{r\in \nat} \im(t_{k+r}^{*}) \ , \quad   \ker((p_{k+r}\circ t_{k+r})^{*})$$  are of  finite type, since they are subgroups of the groups of finite type $K^{*}(T)$ and $K^{*}(X_{k+r})$, respectively.
\begin{fact}[\cite{gray}] If $(G_{r})_{r\in \nat}$ is a projective system of countable abelian groups, then either
$\lim_{r}^{1} G_{r}=0$ or $\lim_{r}^{1}G_{r}$ is uncountable.  
\end{fact} 
Since 
$\lim_{r}^{1}  \ker((p_{k+r}\circ t_{k+r})^{*})$ is a quotient of a countable group 
we conclude that  this group is actually trivial.
The pull-back along the family of maps $(p\circ i_{k+r})_{r\in \nat}$ induces an  isomorphism
 $$\frac{\bar K^{*}(X)}{\ker((p\circ t)^{*})}  \stackrel{\sim}{\to}  \lim_{r}\frac{K^{*}(Y_{k+r})}{\ker(t^{*}_{k+r})}\ .$$
 The element $\phi$ represents a class  
$[\phi]\in \lim_{r}K^{*}(Y_{k+r})/\ker(t^{*}_{k+r})$. Hence we can find a class
$[\psi ]\in \bar K^{*}(X)/\ker((p\circ t)^{*})$ which is mapped to $[\phi]$ under the isomorphism above. The elements $\phi$ and $\psi$ satisfy the relation \eqref{thegoal111}.
 This finishes the proof of Proposition \ref{pinv1}.
\hB

%
%
%
%
%


\color{black}

\bibliographystyle{plain}
\bibliography{universaleta}

\end{document}